\documentclass[11 pt]{article}
\usepackage{fullpage, latexsym,amsmath,amsfonts,amscd,amsthm,upref,graphics,epsfig}

\newtheorem{numbered}{chapter}[section]
\newtheorem{theorem}[numbered]{Theorem}

\newtheorem{nothing*}[numbered]{}
\newtheorem{remark}[numbered]{Remark}

\newtheorem{definition}[numbered]{Definition}

\newtheorem{lemma}[numbered]{Lemma}
\newtheorem{proposition}[numbered]{Proposition}
\newtheorem{corollary}[numbered]{Corollary}
\newtheorem{example}[numbered]{Example}

\numberwithin{equation}{section}
\numberwithin{figure}{section}

\title{Transition from Rotating Waves to Modulated Rotating Waves on the Sphere}
\author{Adela N. Comanici $^{a}$\\
$^{a}$\emph{\small{Department of Mathematics, University of Houston}}\\
\emph{\small{651 Philip G. Hoffman Hall}}\\
\emph{\small{Houston, Texas, 77204-3008, USA}}\\
\emph{\small{adela@math.uh.edu}}}

\begin{document}
\maketitle

\begin{abstract}
In this article, we study parameter-dependent systems of reaction-diffusion equations on the sphere,
which are equivariant under the group $SO(3)$ of all rigid rotations on a sphere. Two main types of spatial-temporal
patterns that can appear in such systems are rotating waves (equilibrium in a co-rotating frame) and modulated rotating
waves (periodic solution in a co-rotating frame). The transition from rotating
waves to modulated rotating waves on spherical domains is explained via a supercritical Hopf bifurcation from a rotating
wave, $SO(3)$-symmetry and finite-dimensional equivariant center manifold reduction. The Baker-Campbell-Haussdorff
formula in the Lie algebra $so(3)$ is used to get reduced differential equations
on $so(3)$, a formula for a primary frequency vector, as well as a formula for the periodic part associated to any
modulated rotating wave obtained by a supercritical Hopf bifurcation from a rotating wave. As a consequence, there are
three types of motions for the tips of these modulated rotating waves.
In the resonant case (with two parameters), we obtain that the primary frequency vectors of a branch of these modulated
rotating waves are generically orthogonal to the frequency vector of the initial rotating wave undergoing Hopf
bifurcation.\\
\end{abstract}

\emph{\small{Keywords $\colon$}} equivariant center manifold, rotating wave, modulated rotating wave, Hopf bifurcation

\section{Introduction}

The main motivation of this article is the presence of spiral waves in excitable media, especially in cardiac tissue. Spiral waves arise as stable spatio-temporal
patterns in various chemical, physical systems and biological systems, as well as numerical simulations of
reaction-diffusion systems on excitable media with various geometries. Excitable media are extended non-equilibrium systems
having a uniform rest state that is linearly stable but susceptible to finite perturbations.
Spiral waves have been observed experimentally, for instance, in catalysis of platinum surfaces \cite{NORE}, Belousov-Zhabotinsky chemical
reactions \cite{JSW,LOPS}, Rayleigh-Bernard convection \cite{PB}, slime-mould cells \cite{Ge} and the most important, cardiac tissue \cite{DPSBJ}.
Numerical simulations of the spiral waves have been done by \cite{Ba,Ba1}, for example.\\
It is now believed that spiral and scroll waves that appear in the heart muscle can lead to cardiac arrhythmias (abnormal rhythms in the heart), giving rise to atrial fluttering or ventricular
fibrillation. In normal hearts cardiac arrhythmias are rare, but in diseased hearts cardiac arrhythmias can become more
common. For example, if chambers of the heart become abnormally large, they are susceptible to serious arrhythmias in
which waves are believed to circulate in a fashion that is similar to the circulation of the Belousov-Zhabotinsky waves
in a chemical medium. Real human hearts are enormously complex three-dimensional structures. In this article, we assume
that the geometry of the excitable media is a sphere, which in the case of cardiac tissue is clearly an approximation.\\
In the planar case, a rigidly rotating spiral wave is an example of wave pattern rotating around a center and being well
approximated by an Archimedean spiral wave far from rotation center. Near the rotation center, there is a core region of
the spiral wave, where the front of the wave has a tip, whose structure is considered to be the most important in
understanding the behavior of the whole spiral wave \cite{KS,Ke}. Barkley \cite{Ba} was the first who performed
a numerical linear stability analysis for the basic-time periodic spiral wave solution in a reaction-diffusion system on
the unbounded plane and showed evidence of a Hopf bifurcation. In particular, a simple pair of eigenvalues was shown to cross the imaginary axis while
three neutral eigenvalues lie on the imaginary axis and the remainder of the spectrum is bounded into the left-half plane.
Later, using an ad hoc model, Barkley \cite{Ba1} was the first to realize the key importance of the group $SE(2)$ of all
planar translations and rotations in describing  the dynamics and bifurcations of planar spiral waves.\\
It is well known now that the tip of the spiral wave rotates steadily or meanders or linearly drift in plane
\cite{GLM1, GLM2}.
From a mathematical point of view, rigidly rotating spiral waves are examples of rotating waves, meandering spiral waves
are examples of modulated rotating waves and linearly drifting spiral waves are examples of modulated rotating waves.\\
The first rigorous mathematical theory of the planar spiral waves was done by Wulff \cite{Wu1}. In her thesis, Wulff
studied the external periodic forcing of rotating waves which leads to modulated rotating waves or modulated travelling
waves, and also proved an $SE(2)$-equivariant Hopf theorem for the bifurcation from rotating waves to modulated
(rotating or travelling) waves in autonomous systems. The external periodic forcing of rotating waves was
studied using a contraction mapping theorem on scales of Banach spaces, and the proof for the $SE(2)$-equivariant Hopf
theorem for rotating waves was based on Liapunov-Schmidt reduction on scales of Banach spaces. In both cases, it is shown
that modulated travelling waves emanate if the rotation frequency is a multiple of the external frequency, respectively
the modulus of the eigenvalue leading to Hopf bifurcation. The main difficulty comes from the fact that the group $SE(2)$
is noncompact and the action of $SE(2)$ on the usual spaces of functions is not smooth. The proofs in \cite{Wu1} are
based on the basic assumption that the linearization at the rotating wave in the co-rotating frame does not exhibit continuous spectrum near the imaginary axis.\\
Later, Sandstede, Scheel and Wulff \cite{SSW} proved a finite-dimensional center bundle reduction theorem near a relative
equilibrium $Gu_{0}$ of an infinite-dimensional vector field on a Banach space $X$ on which acts a finite-dimensional Lie group (not
necessarily compact). This generalizes Krupa's results on bifurcation from relative equilibria \cite{Kr}, from compact
groups to noncompact groups and from finite dimensions to infinite dimensions. Using the results of \cite{SSW}, the Hopf bifurcation
from one-armed rotating spiral wave to meandering waves can be studied \cite{GLM1, SSW}.\\
In the case of a one-armed rotating spiral wave, the relative equilibrium $SE(2)u_{0}$ is diffeomorphic to $SE(2)$,
which is diffeomorphic to $\mathbb{C} \times \mathbf{S}^{1}$. Near Hopf bifurcation, the reduced differential equations
on the center bundle $SE(2) \times \mathbb{C}$ are given by
\begin{equation}\label{E:intro0}
\begin{array}{lll}
\dot{p} &=& e^{i\phi}f(q,\lambda),\\
\dot{\phi} &=& F^{\phi}(q,\lambda),\\
\dot{q} &=& F^{q}(q,\lambda),\\
\end{array}
\end{equation} where $F^{\phi}(0,0)=\omega_{rot}$, $F^{q}(0,0)=0$, $f(0,0)=0$. The rotating wave $u_{0}$ corresponds to
$q= 0$ at $\lambda= 0$ and to the solution $(0,\omega_{rot}t)$ for the first two differential equations in (\ref{E:intro0}).\\
In case of a supercritical Hopf bifurcation we have $d_{q}F^{q}(0,0)= i\omega_{2}$, $Re \, ( d_{q}F^{q}_{\lambda}(0,0)
) \neq 0$. It follows that for any small $\lambda > 0$ there is a meandering spiral wave which becomes a drifting
linear wave if $\omega_{1}(\lambda)= k \omega_{2}(\lambda)$, for some $k \in \mathbb{Z}$, where $\omega_{2}(\lambda)$ is
the frequency that appears due to the Hopf bifurcation and $\omega_{1}(\lambda)= \frac{\left| \omega_{2}(\lambda)\right
|}{2\pi}\int_{0}^{\frac{2\pi}{\left | \omega_{2}(\lambda)\right |}} F^{\phi}(q(t,\lambda),\lambda) \, dt$, and
$q(t,\lambda)$ is the $\frac{2\pi}{\left | \omega_{2}(\lambda) \right |}$ periodic solution of the third differential equation in
(\ref{E:intro0}) that appears due to the supercritical Hopf bifurcation for $q= 0$ at $\lambda= 0$. Similar results for
m-armed spiral waves were also obtained.\\
In \cite{Wu2}, a $G$-equivariant semilinear system of parabolic equations (where $G$ is a finite-dimensional possibly
non-compact Lie group) is studied. In particular, the periodic forcing of relative equilibria and resonant periodic forcing of relative
equilibria to relative periodic orbits, as well as Hopf bifurcation from relative equilibria to relative periodic orbits
are treated using Liapunov-Schmidt reduction. Resonant drift phenomena are also studied.
Then, these results are applied to the planar spiral waves.\\
For the skew-product finite-dimensional system of differential equations on the center manifold near a relative
equilibrium \cite{FSSW}, the normal form method which further simplifies the system is presented in \cite{FT}.
The normal form for the case $G=SE(N)$ is obtained. Then, the known results regarding the meandering and drifting of
planar spiral waves are recovered, as well as new results regarding the relative homoclinic and
heteroclinic trajectories to relative equilibria of $SE(2)$-actions.\\
In \cite{Sc}, Scheel shows that for a large class of reaction-diffusion systems on the plane, m-armed spiral waves
bifurcate from a homogeneous equilibrium when the latter undergoes a Hopf bifurcation. This was done using spatial
dynamics. \\
All previous results are valid for planar spiral waves. The interest to consider spiral waves on non-planar surfaces is
motivated by the applicability to problems in physiology, biology and chemistry. Therefore, the study of spiral waves
by experiments and numerical simulations of reaction-diffusion systems on a sphere have recently been undertaken. In the case of spiral waves
on a sphere, the dynamics is expected to be quite different because any spiral wave starting from a rotating center
cannot end at a point. The number of tips of a wave front cannot be odd, and therefore, the dynamics of spiral waves may
acquire a new feature qualitatively different from the planar case.\\
The dynamics of spiral waves in an excitable reaction-diffusion systems on a sphere was numerically investigated by
\cite{GA, ZM, ZM1}, and \cite{YMY} who employ a spectral method using spherical harmonics as basis
functions. Maselko \cite{Ma}, as well as Maselko and Showalter \cite{MS} performed experiments with Belousov-Zhabotinsky
chemical waves propagating on the surface of a sphere. They observed that a spiral wave winds outward from a meandering
source at the north pole and undergoes self-annihilation as it winds into itself at the south pole.\\
In \cite{ZM}, the evolution of spiral waves on a circular domain and on a spherical surface is studied by numerical
integration of a reaction-diffusion system.
Two different asymptotic regimes are observed for both domains. The first regime is a rigid rotation of an excitation
wave around  the symmetry axis of the domain. The second one is a compound rotation including a drift of the rotation
center of the spiral wave either along the boundary of the disk or along the equator of the sphere. In this case the
shape of the wave and its rotation velocity are periodically changing in time. Simplified analytical estimates are
presented to describe the rigid rotation.\\
In \cite{GA}, numerical integration of an excitable reaction-diffusion system on a sphere is presented.
The evolution of counter-rotating double spiral waves on the sphere is studied and it is shown that the tips of the
spiral can either perform a meandering motion or rigidly rotate around a fixed center, depending on the system control
parameter. It is observed that the rotation of the spiral wave on spherical surface is similar to that obtain  on the
planar surface, except that in the absence of the boundary on a spherical surface some parts of the wave can undergo
self-annihilation in contrast to the spiral wave behavior on bounded planar surfaces.\\
In \cite{ZM1}, the evolution of spiral waves on a spherical surface is studied by integration of a reaction-diffusion
system with a global feedback. It is shown that depending on intensity, sign, and/or time delay in the feedback loop a
global coupling can be effectively used either to stabilize the rigid motion of a spiral wave or to completely destroy
spiral waves and to suppress self-sustained activity in a confined domain of an excitable medium.\\
In \cite{YMY}, the dynamics of chemical spiral waves in an excitable reaction-diffusion system on a sphere
is numerically investigated employing a spectral method using spherical harmonics as basis functions. Different types
of spiral waves --symmetric or antisymmetric source-source or nearly antisymmetric source-nonsource-- have been obtained
depending on whether the medium is homogeneous or inhomogeneous, and it has been observed that the tips can either rotate
steadily or change their shapes.\\
The influence of the topological constraints and the inhomogeneity in the excitability on the geometry and dynamics of spiral
waves on a thin spherical shell of excitable media are presented in \cite{DGK}.
Also, rigidly rotating waves on spherical domains have been studied using kinematical theory by \cite{GH, GG, MG}.\\
In \cite{GH}, the geometrical stability of the symmetric counter-rotating spiral waves propagating on the unit sphere is
studied. By the use of the eikonal equations, it is demonstrated that these solutions are stable under small perturbations
normal to the wave front lying on a unit sphere.\\
In \cite{GG}, the authors showed that stationary rotating solutions on a sphere of the eikonal equation
under some boundary conditions on tips must be symmetric with respect to the equator with spiral winding out from
source tips at polar points.\\
In \cite{MG}, using the eikonal approximation to a reaction-diffusion system on a sphere, the authors prove
existence of a class of counter-rotating double spiral solutions which are highly asymmetric with respect to the equator.
They also derive a power law, linking the angular rotation of the spiral waves with the velocity of plane waves in
medium.\\
In \cite{Re}, Renardy considered bifurcations from rotating waves of semi-linear equations that are equivariant under a
general compact Lie group and applied his results to the Laser equations. His results do not cover the resonance case,
which will be covered in our paper. The theorems were proved using a generalized implicit function theorem on scales of
Banach spaces.\\
In \cite{Ra}, Rand examined modulated rotating waves in rotating fluids and applied his results to the Taylor-Couette
problem, where modulated rotating waves, so called modulated wavy vortices occur.\\
In this article we have studied only the trivial isotropy case, because as far as we are aware, there have been observed no m-spiral
waves ($m > 1$) on spherical surfaces. Also, in \cite{DGK} it was numerically verified that there is a critical size of
the sphere below which self-sustained spiral waves cannot exist. Therefore, we use a sphere of an arbitrary,
but fixed radius $r$. \\
In this article, we study non-resonant and resonant supercritical Hopf bifurcation from a rotating wave. The relation
between spiral waves and rotating waves, respectively modulated rotating waves will be presented somewhere else. While we have
concentrated on the case of Hopf bifurcation, it should be mentioned that we would obtain similar results for
\emph{ periodic forcing } of  a rotating wave.\\
In Section \ref{S:RDS}, we define the representation of the Lie group $SO(3)$ on the usual function space.
It is a smooth unitary representation on the fractional spaces relative to the appropriate sectorial operator. We also associate an
$SO(3)$-equivariant semiflow to the reaction-diffusion system on a sphere, and we recall the notion of tip position
function. In Section \ref{S:RW_MRW}, we recall the notions
of relative equilibrium, relative periodic orbit, rotating wave and modulated rotating wave. We also define the notions of
frequency vector of a rotating wave and primary frequency vector of a modulated rotating wave, and present some of their
properties. Section \ref{S:HopfBif} deals with a supercritical Hopf bifurcation from a rotating wave. Decomposition of a
modulated rotating wave obtained by Hopf bifurcation into a primary frequency part and the associated periodic part is
proved. In Section \ref{S:Frequency}, the form of the reduced differential equations on $so(3)$,
as well as formulae for a primary frequency vectors and the associated periodic part are obtained using
Baker-Campbell-Hausdorff formula in the Lie algebra $so(3)$, the properties of the exponential map $e \colon so(3) \rightarrow
SO(3)$ and the properties of the adjoint representations of $SO(3)$ and $so(3)$. As a consequence, there are three
types of motions for the tips of the modulated rotating waves obtained by a Hopf bifurcation from a rotating wave. Namely,
the tip can either quasi-periodically meander on the sphere, such that the primary frequency vector is of order $O(1)$
near the frequency vector of the initial rotating wave,
or quasi-periodically meander very slowly (order $O(\sqrt{\lambda})$) about the primary frequency vector, which can be
orthogonal or not to the frequency vector of the initial rotating wave undergoing Hopf bifurcation. In the orthogonal case,
the tip motion can be approximated with a very slow drift (order $O(\sqrt{\lambda})$) along an equator of the sphere.
Section \ref{S:Examples_Transition} presents some examples that illustrate the theoretical results obtained in
Sections \ref{S:HopfBif} and \ref{S:Frequency}.
In Section \ref{S:Res_drift}, the resonant Hopf bifurcation from a rotating wave with two parameters is studied.
Near the resonant Hopf bifurcation, we proved that there exists generically a unique branch of modulated rotating waves with primary
frequency vectors orthogonal to the frequency vector of the rotating wave undergoing Hopf bifurcation.
Then, Section \ref{S:Numerical} deals with some numerical results that illustrate the theoretical results obtained in
Sections \ref{S:HopfBif} and \ref{S:Frequency}. The proofs of all theorems are presented in Section
\ref{S:Proof_Transition}. Three appendices are included, presenting $BCH$ formula in $so(3)$, some properties of the
exponential map $e \colon so(3) \rightarrow SO(3)$ and of the adjoint representations of $SO(3)$ and $so(3)$, the
equivariant center manifold reduction theorem and some computations done in Section \ref{S:Numerical}.\\

\section{Reaction-diffusion Systems on Sphere $r\mathbf{S^{2}}$}\label{S:RDS}

Let $r > 0$ and $\mathbf{S^{2}}$ be the unit sphere in $\mathbb{R}^{3}$. We consider a reaction-diffusion system of the
form
\begin{equation}\label{E:RD}
\frac{\partial{u}}{\partial{t}}(t,x)= D\Delta_{S}{u}(t,x)+F(u(t,x)) \mbox{ on } r\mathbf{S^{2}},
\end{equation}
$u= (u_{1}$,~$u_{2}$, \ldots, ~$u_{N}) \colon \mathbb{R} \times r\mathbf{{S}^{2}} \rightarrow \mathbb{R}^{N}$ with $N \geq 1$,\\
$D= \left ( \begin{array}{ccc}
            d_{1}  & \dots  & 0 \\
             \vdots & \ddots & \vdots\\
            0      & \dots  & d_{N}\\
            \end{array}
    \right )$ with $d_{i} \geq 0$ for $i=1$,~$2$, \ldots, ~$N$ the diffusion coefficients, $\Delta_{S}$ is
Laplace-Beltrami operator on $r\mathbf{S^{2}}$ and $F= (F_{1}$,~$F_{2}$, \ldots, $F_{N}) \colon \mathbb{R}^{N} \rightarrow  \mathbb{R}^{N}$.
Also, $D \Delta_{S} u  =\left (\begin{array}{c}
                             d_{1}\Delta_{S} u_{1} \\
                             d_{2}\Delta_{S} u_{2} \\
                             \vdots\\
                             d_{N}\Delta_{S} u_{N}\\
                               \end{array}
                        \right )$. \\
We study the reaction-diffusion system $(\ref{E:RD})$ on the function space
\[\mathbf{Y}= \begin{cases}
               \mathbf{L^{2}}(r\mathbf{S^{2}}, \mathbb{R}^{N})&\text{if  $d_{i}> 0$ for $i=1$,~$2$, \ldots, ~$N$;}\\
               \mathbf{H^{2}}(r\mathbf{S^{2}}, \mathbb{R}^{N}) &\text{if there exists $i \in \{ 1$,~$2$ \ldots, ~$N \}$ such that $d_{i}= 0$.}
              \end{cases}
\]
Let $\alpha \in (\frac{1}{2},1)$ and $\mathbf{Y^{\alpha}}=D((I-D\Delta_{S})^{\alpha})$ be the fractional spaces of $\mathbf{Y}$ relative
to the sectorial operator $-D\Delta_{S}$. Using \cite{Fr, Heb, He, Pa, Ta, Wu1}, to the
reaction-diffusion system (\ref{E:RD}) we associate a local semiflow on an open set
of $\mathbf{Y}$. Let $F \colon \mathbb{R}^{N} \rightarrow \mathbb{R}^{N}$ be a $C^{k+2}$ function such that $F(0)= 0$, where
$0 \leq k \leq \infty $. Then the reaction-diffusion system (\ref{E:RD}) defines a sufficiently smooth local semiflow
$\mathbf{\Phi}$ on the function space $\mathbf{Y^{\alpha}}$. The flow is $C^{k+2}$ if $\mathbf{Y}= \mathbf{L^{2}}(r\mathbf{S^{2}},
\mathbb{R}^{N})$ and $C^{k}$ if $\mathbf{Y}=\mathbf{H^{2}}(r\mathbf{S^{2}}, \mathbb{R}^{N})$. Namely \cite{Pa}, for any $u_{0} \in \mathbf{Y}^{\alpha}$, let $u(t, u_{0})$
be the sufficiently smooth solution ($C^{k+2}$ or $C^{k}$  respectively) of the initial value problem given by the
reaction-diffusion system (\ref{E:RD}) and by the initial condition $u(0)= u_{0}$, defined on the maximal interval of
existence $I(u_{0})=[0,t_{0}(u_{0}))$. Let $W=\{(t,u_{0}) \in [0,\infty) \times \mathbf{Y}^{\alpha} \mid t \in I(u_{0})\}$.
Then, the local semiflow $\mathbf{\Phi} \colon W \rightarrow \mathbf{Y}^{\alpha}$ is defined by $\mathbf{\Phi}(t,u_{0})=u(t,u_{0})$, for any $(t, u_{0}) \in W$.
Any reaction-diffusion system of the above type has the property that is $SO(3)$-equivariant, where the precise meaning of
$SO(3)$-equivariance is defined below.
\begin{definition}\label{def:Action1}
The \emph{representation} $T$ of $SO(3)$ on $\mathbf{Y}$ is the function $T \colon SO(3) \rightarrow GL(\mathbf{Y})$, defined by
\begin{equation}\label{E:Rep}
T(A)u(x)= u(A^{-1}x) \mbox{ where } A \in SO(3), u \in \mathbf{Y}, x \in r\mathbf{S^{2}}.
\end{equation}
\end{definition}
The linear action $\theta \colon SO(3) \times \mathbf{Y}\rightarrow \mathbf{Y}$ associated to the unitary representation
$T$ (\cite{Heb, He, LM, Ta, Wu1}) is defined by
\begin{equation}\label{E:Action}
\theta(A,u)(x)= u(A^{-1}x) \mbox{ where } A \in SO(3), u \in \mathbf{Y} \mbox{ and } x \in r\mathbf{S^{2}}.
\end{equation}
We denote $\theta(A,u)=Au$. From now on, when we talk about $SO(3)$-equivariance, we mean equivariance with respect to the
action $\theta$ defined in (\ref{E:Action}).
Using \cite{DM, Su, Te, Va1, Va2}, we obtain the following two propositions:
\begin{proposition}\cite{Co} \label{thm:Smoothness}
The restriction of the representation $T$ to $\mathbf{Y^{\alpha}}$, where $\alpha \in (\frac{1}{2},1)$, is smooth.
\end{proposition}
\begin{proposition}\cite{Co}\label{thm:Equivariant}
Let $F \colon \mathbb{R}^{N} \rightarrow \mathbb{R}^{N}$ be a $C^{k+2}$ function such that $F(0)= 0$ and
$1 \leq k \leq \infty$. The local semiflow $\mathbf{\Phi}$  is $SO(3)$-equivariant with respect to the
action $\theta$ restricted to $\mathbf{Y^{\alpha}}$.
\end{proposition}
For any $u_{0} \in \mathbf{Y}^{\alpha}$ we follow the time evolution of $u_{0}$ under the local semiflow $\mathbf{\Phi}$
using the concept of tip motion of $u_{0}$.
\begin{definition}\label{def:tip}
For \, $0 \leq k \leq \infty$, a $C^{k}$ function $x_{tip} \colon \mathbf{Y^{\alpha}} \rightarrow r\mathbf{S}^{2}$ that
is $SO(3)$-equivariant is called a \emph{tip position} function. By the \emph {tip} of $u \in \mathbf{Y^{\alpha}}$, we understand the point
$x_{tip}(u) \in r\mathbf{S^{2}}$. Also, for any $u_{0} \in \mathbf{Y^{\alpha}}$, $x_{tip}(\mathbf{\Phi}(t, u_{0}))$
with $t \in [0, \infty)$ is called the \emph{tip motion} of $u_{0}$. If the function $x_{tip}$ is defined only on an
open subset of $\mathbf{Y^{\alpha}}$, then the function $x_{tip}$ is called a \emph{local tip position} function.
\end{definition}
The tip position function is usually used to follow the time evolution of $u_{0}$ under the local semiflow
$\mathbf{\Phi}$, when $\mathbf{\Phi}(t, u_{0})$ is a rotating wave or a modulated rotating wave.

\section{Rotating Waves and Modulated Rotating Waves on $r\mathbf{S^{2}}$}\label{S:RW_MRW}

We consider that north hemisphere of $\mathbf{S}^{2}$ is the set $N=\{(x,y,z) \in \mathbf{S}^{2}
\mid z>0 \mbox{ or } z=0, x \in [-1,1), y \in [0,1]\}$, and south hemisphere is the set $S=\mathbf{S}^{2}\setminus N$.\\
For any $G$-equivariant dynamical system, where $G$ is a Lie group, we can talk about the concepts of relative equilibrium
and relative periodic orbit. If $SO(n) \subset G$ for some integer $ n \geq 2$, we can also talk about rotating waves
and modulated rotating waves.
\begin{definition} \cite{Fi1, Kr, Wu1} \label{def:RE}
Let $u_{0} \in \mathbf{Y^{\alpha}}$ be such that the stabilizer of $u_{0}$ is $\Sigma_{u_{0}}= {I_{3}}$.
The orbit group of $u_{0}$ is called a \emph{relative equilibrium} for (\ref{E:RD}) if there exists a matrix $X_{0} \in so(3)$ such that
\begin{equation}\label{E:RE}
\mathbf{ \Phi}(t, u_{0})= e^{X_{0}t}u_{0} \mbox{ for any } t \geq 0.
\end{equation}
Sometimes, when $SO(3)u_{0}$ is a relative equilibrium, we call $u_{0}$ a \emph{relative equilibrium}.
If $X_{0} \neq O_{3}$, then any solution of the reaction-diffusion system (\ref{E:RD}) of the form
$A\mathbf{ \Phi}(., u_{0})$, where $A \in SO(3)$ is called a \emph{rotating wave} for (\ref{E:RD}).
\end{definition}
Because the action $\theta$ restricted to $\mathbf{Y^{\alpha}}$ is smooth, then the relative equilibrium from Definition
\ref{def:RE} is a sufficiently smooth manifold in $\mathbf{Y^{\alpha}}$ diffeomorphic to $SO(3)$ since
$\Sigma_{u_{0}}=I_{3}$.
\begin{definition} \label{def:frequency_RW}
Let $\mathbf{\Phi}(., Au_{0})$ be a rotating wave as in Definition \ref{def:RE}, where $A \in SO(3)$.
Then, the vector $\overrightarrow {AX_{0}A^{-1}}$ is called the \emph{frequency vector} of the rotating wave.
If the vector $\frac{1}{\left | X_{0} \right | }\overrightarrow {AX_{0}A^{-1}}$ is in the north hemisphere,
then $\omega_{0}= \left | X_{0} \right |$ is called the \emph{frequency} of the rotating wave. Otherwise, $-\omega_{0}$ is called the \emph{frequency} of the rotating wave.
\end{definition}
Then, the rotating waves in $SO(3)u_{0}$ have their anti-symmetric matrices associated with the frequency vectors on the same adjoint orbit.
The following property shows the importance of the frequency vector associated to any rotating wave:
\begin{proposition}\cite{Co}\label{prop:TipmRW}
For a rotating wave $\mathbf{\Phi}(t, u_{0})$, the tip motion $x_{tip}(\mathbf{\Phi}(t, u_{0}))$ is a circle on the
sphere $r\mathbf{S}^{2}$ with the center on the line having the direction of the frequency vector of $\mathbf{\Phi}(t, u_{0})$, and this is independent of the choice of the tip position function.
\end{proposition}
\begin{definition} \cite{Fi1, Kr, Wu1} \label{def:RPO}
Let $u_{0} \in \mathbf{Y^{\alpha}}$ be such that the stabilizer of $u_{0}$ is $\Sigma_{u_{0}}= {I_{3}}$.
The set defined by $\{A \mathbf{\Phi}(t, u_{0}) \mid A \in SO(3), t \in [0, \infty)\}$ is called a
\emph{relative periodic orbit} for $(\ref{E:RD})$ if it is not a relative equilibrium and there exist a number $T > 0$ and a matrix $X_{0} \in so(3)$ such that
\begin{equation}\label{E:RPO}
\mathbf{ \Phi}(T, u_{0})= e^{X_{0}T}u_{0} \mbox { and } \mathbf{\Phi}(t, u_{0}) \notin SO(3)u_{0} \mbox { for any }
t \in (0,T).
\end{equation}
If $\left | X_{0} \right | T \neq 2k\pi$ for any $k \in \mathbb{Z}$, then any solution of the reaction-diffusion
system (\ref{E:RD}) of the form $A\mathbf{\Phi}(., u_{0})$, with $A \in SO(3)$ is called a \emph{modulated rotating wave} for (\ref{E:RD}).
\end{definition}
If $\left | X_{0} \right | = \frac{2k\pi}{T}$ for some $ k \in \mathbb{Z}$, then $\mathbf{ \Phi}(t, u_{0})$ is a
$T$-periodic solution of the reaction-diffusion (\ref{E:RD}).
Because the action $\theta$ of $SO(3)$ on $\mathbf{Y^{\alpha}}$ is smooth, then the relative periodic orbit defined in
Definition \ref{def:RPO} is a sufficiently smooth manifold in $\mathbf{Y^{\alpha}}$.
There are two properties related to modulated rotating waves, namely:
\begin{remark}\cite{Co} \label{rem:rot}
\begin{enumerate}
\item $\mathbf{ \Phi}(T_{1}, u_{0})= e^{X_{0}T_{1}}u_{0}$ if and only if $T_{1} \in T\mathbb{Z}$.
\item There are periodic solutions of the reaction-diffusion system (\ref{E:RD}) that are modulated rotating waves in the sense of Definition \ref{def:RPO}. Then their period is an integer multiple of the corresponding $T$ from Definition \ref{def:RPO} and $\left | X_{0} \right | T \in 2\pi \mathbb{Q}$.
\end{enumerate}
\end{remark}
\begin{definition}\label{def:frequency_MRW}
Let $\mathbf{\Phi}(., Au_{0})$ be a modulated rotating wave as in Definition \ref{def:RPO}, where $A \in SO(3)$.
Then, the vector $\overrightarrow {AX_{0}A^{-1}}$ is called a \emph{primary frequency vector} of the modulated rotating wave and the positive number $T$ is called the \emph{relative period} of the modulated rotating wave. If the vector $\frac{1}{\left | X_{0} \right |  }\overrightarrow {AX_{0}A^{-1}}$ is in the north hemisphere, then
$\omega_{0}= \left | X_{0} \right | $ is called the \emph{primary frequency} of the modulated rotating wave. Otherwise, $-\omega_{0}$ is called the \emph{primary frequency} of the modulated rotating wave.
\end{definition}
A primary frequency vector $\overrightarrow{X_{1}}$ of a modulated rotating wave $\mathbf{\Phi}(t, Au_{0})$ is unique up
to an integer multiple of $\frac{2\pi}{\left | X_{1} \right | T}\overrightarrow{X_{1}}$. Then the modulated rotating
waves in $SO(3)\mathbf{\Phi}(.,u_{0})$ have the anti-symmetric matrices associated with primary frequency vectors on
the same adjoint orbit (if we consider primary frequency vectors to be unique up to a multiple of
$\frac{2k\pi}{T}$, $k \in \mathbb{Z}$ of the corresponding unit primary frequency vector).
If $\mathbf{\Phi}(t, Au_{0})$ is a $T$-periodic solution of the reaction-diffusion system (\ref{E:RD}) such that
$\mathbf{ \Phi}(T, Au_{0})= Au_{0}$ and $\mathbf{ \Phi}(t, Au_{0}) \notin SO(3)u_{0}$ for all $t \in (0,T)$,
then a primary frequency vector of $\mathbf{\Phi}(t, Au_{0})$ can be any vector $\overrightarrow{X_{1}} \in
\mathbb{R}^{3}$ with ${\left | X_{1} \right | }= \frac{2k\pi}{T}$ for $k \in \mathbb{Z}$.
The following property shows the importance of the primary frequency vector associated to any modulated rotating wave:
\begin{proposition}\cite{Co} \label{prop:TipmMRW}
For a modulated rotating wave $\mathbf{\Phi}(t,u_{0})$, the tip motion $x_{tip}(\mathbf{\Phi}(t, u_{0}))$  has
the property that $x_{tip}(\mathbf{\Phi}(kT, u_{0}))$ for $k \in \mathbb{Z}$ are points of  a circle on the sphere
$r\mathbf{S}^{2}$ with the center on the line having the direction of the primary frequency vector of $\mathbf{\Phi}(t, u_{0})$; this is independent of the choice of the tip position function.
\end{proposition}

\section{Hopf Bifurcation from Rotating Waves to Modulated Rotating Waves on $r\mathbf{S^{2}}$}\label{S:HopfBif}

\begin{figure}
\begin{center}
\psfig{file=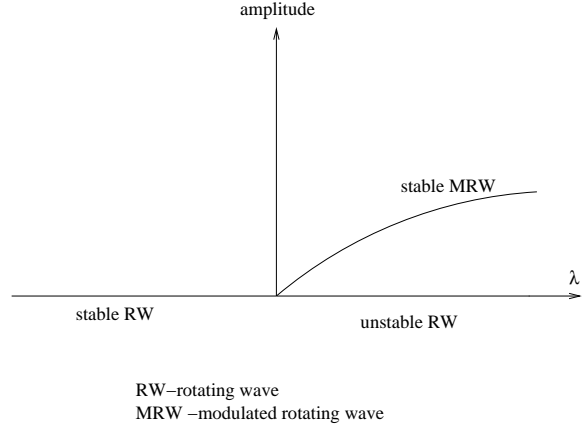,width=3in}
\end{center}
\caption{Bifurcation Diagram for Nonresonant Hopf Bifurcation}
\label{figure:figura1}
\end{figure}

In this section, we will consider a  supercritical Hopf bifurcation from an equilibrium. Because the bifurcating periodic solution  has amplitude of order $\sqrt{\lambda}$, the following definition will simplify the language.
\begin{definition}
Let $M$ be a smooth manifold, $X$ a normed space or the empty set, $p \geq 0$ an integer and
$Y \colon X \times [0, \lambda_{0}) \times \mathbb{R}^{p} \rightarrow M$ for $\lambda_{0} >0$ small. We say that the function $Y$ is \emph{$CS$} on $X \times [0, \lambda_{0}) \times \mathbb{R}^{p}$ if the function $Z \colon X \times [0, \epsilon_{0}) \times \mathbb{R}^{p} \rightarrow M$  defined by \[Z(x,\epsilon, \mu)= Y(x, \epsilon^{2}, \mu),\] is smooth on  $X \times [0,\epsilon_{0}) \times \mathbb{R}^{p}$, where $\epsilon_{0}=\sqrt{\lambda_{0}}$.\\
We say that $Y$ is \emph{$C^{k}$-$CS$} if $Z$ is $C^{k}$, where $k \in \mathbb{Z}$, $ k \geq 1$. We say that
$Y$ is \emph{sufficiently $CS$} if $Z$ is sufficiently smooth.
\end{definition}
Let us consider the following reaction-diffusion system
\begin{equation}\label{E:RDparam}
\frac{\partial{u}}{\partial{t}}(t,x)= D\Delta_{S}{u}(t,x)+F(u(t,x),\lambda) \mbox{ on } r\mathbf{S^{2}},
\end{equation}
$u= (u_{1}$,~$u_{2}$, \ldots, ~$u_{N}) \colon \mathbb{R} \times r\mathbf{{S}^{2}} \rightarrow \mathbb{R}^{N}$ with $N \geq 1 $, $D= \left(\begin{array}{ccc}
                   d_{1}  & \dots  & 0 \\
                   \vdots & \ddots & \vdots\\
                    0      & \dots  & d_{N}
                    \end{array}
                 \right)$, $d_{i} \geq 0$, $i= 1$,~$2$, \ldots, $N$ the diffusion coefficients and
$F= (F_{1}$,~$F_{2}$, \ldots, ~$F_{N}) \colon \mathbb{R}^{N}\times \mathbb{R} \rightarrow \mathbb{R}^{N}$ a $C^{k+2}$ function with $3 \leq k \leq \infty $. We study the reaction-diffusion system (\ref{E:RDparam}) on the function space $\mathbf{Y}$ defined in Section \ref{S:RDS}.
Let $u_{0} \in \mathbf{Y^{\alpha}}$ be a relative equilibrium that is not an equilibrium for (\ref{E:RDparam}) at $\lambda= 0$ and such that the stabilizer of $u_{0}$ is $\Sigma_{u_{0}}= {I_{3}}$.  Let $\mathbf{\Phi}(t, u_{0}, 0)= e^{X_{0}t}u_{0}$. Consider $L$ the linearization of the right-hand side of (\ref{E:RDparam}) with respect to the rotating wave $\mathbf{\Phi}(t, u_{0}, 0)= e^{X_{0}t}u_{0}$ at $\lambda= 0$ in the co-rotating frame, that is
\[L= D\Delta_{S}+D_{u}F(u_{0},0)-X_{0}.\]
Suppose that:
\begin{enumerate}
\item $\sigma(L) \cap \{z \in \mathbb{C} \mid Re \, (z) \geq 0 \}$ is a spectral set with spectral projection $P_{*}$,
and dim$(R(P_{*})) < \infty $;
\item the semigroup $e^{Lt}$ satisfies $\left |e^{Lt}|_{R(1-P_{*})} \right | \leq Ce^{-\beta_{0}t}$ for
some $\beta_{0} > 0$ and $C > 0$.
\end{enumerate}
Theorems \ref{thm:CMR_REparam} (see Appendix \ref{S:app2_CMR}) can be applied.
Then, there exist sufficiently smooth functions $X_{G} \colon V_{*} \times \mathbb{R} \rightarrow so(3)$ and $X_{N} \colon V_{*} \times \mathbb{R} \rightarrow V_{*}$ such that any solution of
\begin{equation}\label{E:reduceddiffeqparam1}
    \begin{array}{lll}
    \dot{A} &=& AX_{G}(q, \lambda),\\
    \dot{q} &=& X_{N}(q, \lambda),\\
    \end{array}
\end{equation}
on $SO(3) \times V_{*}$ corresponds to a solution of the reaction-diffusion system (\ref{E:RDparam}) on
$M_{u_{0}}^{cu}(\lambda)$ under the diffeomorphic identification for $\left | \lambda  \right | $ small.
Also, $X_{G}(0, 0)= X_{0}$, $X_{N}(0, 0)= 0$ and $\sigma(D_{u}X_{N}(0,0))= \sigma(Q_{*}L|_{V_{*}})$, where $Q_{*}$
is the projection onto $V_{*}$ along $T_{u_{0}}(SO(3)u_{0})$.
An example of tip position function is $x_{tip} \colon M_{u_{0}}^{cu}(\lambda) \times \mathbb{R} \rightarrow
r\mathbf{S}^{2}$ defined by $x_{tip}(A\mathbf{\Psi}(q), \lambda)=Ax_{0}$, where $x_{0} \in r\mathbf{S}^{2}$ is fixed.
Also, we have the following result. Suppose we omit the parameter $\lambda$. Let $\mathbf{\Phi}(., u_{1})$ be a
modulated rotating wave as defined in Definition \ref{def:RPO} with $u_{1} \in M^{cu}_{u_{0}}$ and $\Sigma_{u_{1}}= I_{3}$,
where $M^{cu}_{u_{0}}$ is defined in Theorem \ref{thm:CMR_REparam}. Suppose that $\mathbf{\Phi}(., u_{0})$ corresponds to
$0$ and $\mathbf{\Phi}(., u_{1})$ corresponds to $q_{1}(t)$, where $0$ and $q_{1}(t)$ are solutions of the second
differential equation in (\ref{E:reduceddiffeqparam1}). Using \cite{FSSW, GSS, SSW}, we have:
\begin{proposition} \cite{Co}\label{prop:Stability_RW}
\begin{enumerate}
\item $\mathbf{\Phi}(., u_{0})$ is orbitally stable (respectively unstable) if $0$ is stable (respectively unstable) in the second differential equation of (\ref{E:reduceddiffeqparam1}).
\item $\mathbf{\Phi}(., u_{1})$ is orbitally stable (respectively unstable) if $q_{1}(t)$ is stable (respectively unstable) in the second differential equation of (\ref{E:reduceddiffeqparam1}).
\end{enumerate}
\end{proposition}
Let dim$(R(P_{*}))=5$. Suppose that a supercritical Hopf bifurcation with eigenvalues $\pm i \omega_{bif}$ takes place
in the second differential equation of (\ref{E:reduceddiffeqparam1}) in $V_{*}$ at $q= 0$ for $\lambda= 0$
\cite{GS, HS, Wi}. Namely,
\begin{enumerate}
\item $X_{N}(0, 0)= 0$;
\item $D_{q}X_{N}(0, 0)$ has eigenvalues $\pm i \omega_{bif}$; without loss of generality, we assume that
$X_{N}(0, \lambda)= 0$ for $\left | \lambda \right | $ small;
\item $D_{q}X_{N}(0, \lambda)$ has the eigenvalues $\alpha(\lambda) \pm i(\omega_{bif}+\beta(\lambda))$
with $\alpha(0)= \beta(0)= 0$ such that $\alpha^{'}(0) > 0$;
\item the branch of periodic solutions $q(t,\lambda)$ bifurcating from $q= 0$ satisfies $q(t, \lambda)=
O(\sqrt{\lambda})$ for $\lambda \geq 0$ small;
\end{enumerate}
Let $q(t, 0)=0$ for all $t \in [0, \infty)$. Let $T(\lambda)= \frac{2\pi}{\left | \omega_{\lambda} \right |}$ be
the period of the solution $q(t,\lambda)$ near $q= 0$ that appears for $\lambda > 0$ small due to the
supercritical Hopf bifurcation, where $\omega_{\lambda}= \omega_{bif}+O(\lambda)$ for $\lambda \geq 0$ small.
Let us define
\begin{equation}\label{E:xg}
X^{G}(t,\lambda)=X_{G}(t, \lambda)=
               \begin{cases}
                       X_{G}(q(t, \lambda), \lambda)
                        &\text{if $\lambda > 0$ and $t \in [0, \infty),$}\\
                       X_{0} &\text{if $\lambda= 0$ and $t \in [0, \infty),$}
               \end{cases}
\end{equation}
where $X_{G}(q, \lambda)$ is defined in Theorem \ref{thm:CMR_REparam}(see Appendix \ref{S:app2_CMR}).\\
Since $q(t,\lambda)=\sqrt{\lambda}r(t,\lambda)$ is sufficiently $CS$ and $X_{G}(q,\lambda)$ is sufficiently smooth,
it follows that the function $X^{G}$ is sufficiently $CS$.\\
Writing $q(t, \lambda)=\sqrt{\lambda}r(t,\lambda)$, we have
\begin{equation}\label{E:InitialData}
\begin{array}{lll}
X^{G}(t, \lambda) &=& x_{0}(t,\lambda)X_{0}^{1}+\sqrt{\lambda}x_{1}(t, \lambda)X_{1}
+\sqrt{\lambda}x_{2}(t, \lambda)X_{2},\\
x_{0}(t, \lambda) &=& \left | X_{0} \right |+\sqrt{\lambda}x_{01}(t)+ \lambda x_{02}(t, \lambda)
\end{array}
\end{equation}
for $t \in [0,\infty)$ and $\lambda \geq 0$ small. The functions $X^{G}(t,\lambda)$, $x_{0}(t,\lambda)$,
$x_{1}(t,\lambda)$, $x_{2}(t, \lambda)$ are $\frac{2\pi}{\left | \omega_{\lambda}\right |}$-periodic in $t$
for $\lambda \geq 0$ small.\\
Let $A(t, \lambda)$ for $t \in [0,\infty)$ and $\lambda \geq 0$ small be the solution of the initial value
problem
\begin{equation}\label{E:AuxFreq}
\begin{array}{lll}
\dot{A} &=& AX^{G}(t, \lambda),\\
A(0) &=& I_{3},
\end{array}
\end{equation}
where $X^{G}(t, \lambda)$ is defined in (\ref{E:xg}).
\begin{lemma}[Decomposition of $A(t, \lambda)$]\label{lem:Hopf_group}
Suppose that the assumptions made in this section hold. Consider the initial value problem
\begin{equation}\label{E:auxHopf}
\begin{array}{lll}
\dot{A} &=& AX^{G}(t, \lambda),\\
A(0) &=& I_{3},\\
\end{array}
\end{equation}
where $X^{G}(t, \lambda )$ defined by (\ref{E:xg}) is sufficiently $CS$ and
$\frac{2\pi}{\left | \omega_{\lambda}\right |}$-periodic in $t$ for $\lambda > 0$ small. Then, there exists a
sufficiently $CS$ solution $A(t, \lambda)= e^{X(\lambda)t}B(t, \lambda)$ of the initial value problem
(\ref{E:auxHopf}), where $X(\lambda) \in so(3)$ and $B(t, \lambda)$ is a
$\frac{2\pi}{\left | \omega_{\lambda} \right | }$-periodic function such that $B(0,\lambda)= I_{3}$.
\end{lemma}
We call $B(t, \lambda)$ the \emph{periodic part} of $\mathbf{\Phi}(t,u_{\lambda},\lambda)$ associated to
$\overrightarrow{X(\lambda)}$. Using Lemma \ref{lem:Hopf_group} we get the following result:
\begin{theorem}[Hopf Bifurcation Theorem for Rotating Waves on a Sphere]\label{thm:Hopf_RW}
Suppose that the assumptions made in this section hold. Then, there exists a sufficiently $CS$ branch $\mathbf{\Phi}(t,u_{\lambda}, \lambda)$
of solutions for the reaction-diffusion system (\ref{E:RDparam}) such that $\mathbf{\Phi}(t, u_{0},0)=e^{X_{0}t}u_{0}$
and for $\lambda > 0$ small, $\mathbf{\Phi}(t,u_{\lambda}, \lambda)$ is either an orbitally stable modulated rotating
wave with a primary frequency vector $\overrightarrow{X(\lambda)}$ and the secondary frequency $\omega_{\lambda}$, or an
orbitally stable periodic solution with the period $\frac{2\pi}{\left | \omega_{\lambda} \right | }$.
\end{theorem}
A slightly modified version of this theorem was proved by \cite{Re, Wu2}. They used the Liapunov-Schmidt reduction
on scales of Banach spaces. We use the equivariant manifold reduction and Lemma \ref{lem:Hopf_group}.\\
The bifurcation diagram is illustrated in Figure \ref{figure:figura1}.

\section{The Primary Frequency Vector Formula for $\mathbf{\Phi}(t,u_{\lambda},\lambda)$}\label{S:Frequency}

Let $X_{0}^{1}= \frac{1}{\left | X_{0} \right |} X_{0}$. There exist $X_{1},X_{2} \in so(3)$ such that the set
$\{\overrightarrow X_{0}^{1}, \overrightarrow X_{1}, \overrightarrow X_{2} \}$ is an orthonormal basis in
$\mathbb{R}^{3}$ satisfying
\[\overrightarrow X_{0}^{1} \times \overrightarrow X_{1}= \overrightarrow X_{2},
\overrightarrow X_{1} \times \overrightarrow X_{2}= \overrightarrow X_{0}^{1} \mbox{ and }
\overrightarrow X_{2} \times \overrightarrow X_{0}^{1}=\overrightarrow X_{1}.\] Since $so(3)$ is isomorphic to
$\mathbb{R}^{3}$ and $\overrightarrow {[X,Y]}=\overrightarrow X \times \overrightarrow Y$,
$\{X_{0}^{1}, X_{1}, X_{2}\}$ is a basis of the Lie algebra $so(3)$ such that
\[ [X_{0}^{1}, X_{1}]=  X_{2}, [X_{1},  X_{2}]= X_{0}^{1} \mbox{ and } [ X_{2}, X_{0}^{1}]= X_{1}.\]
Let us consider the reaction-diffusion system (\ref{E:RDparam}).\\
Throughout this section, we will suppose that $u_{0} \in \mathbf{Y^{\alpha}}$ is a relative equilibrium that
is not an equilibrium for (\ref{E:RDparam}) at $\lambda= 0$ and such that the stabilizer of $u_{0}$ is
$\Sigma_{u_{0}}= {I_{3}}$. Also, let $\mathbf{\Phi}(t, u_{0}, 0)= e^{X_{0}t}u_{0}$.\\
Suppose that Theorem \ref{thm:Hopf_RW} holds. Recall that a primary frequency vector $\overrightarrow{X(\lambda)}$
of $\mathbf{\Phi}(t, u_{\lambda}, \lambda)$ is given by $A(\frac{2\pi}{\left | \omega_{\lambda}\right |}, \lambda)=
e^{X(\lambda)\frac{2\pi}{\left | \omega_{\lambda} \right | }}$ for $\lambda > 0$ small.

\subsection{Existence of a \emph{sufficiently $CS$} branch $X(\lambda)$}\label{SS:Existence}

The theorems stated in this subsection can be proved without using $BCH$ formula in $so(3)$ \cite{Co}.\\
Let us define the smooth following function
\begin{equation}\label{E:qiuinverse}
q([\overrightarrow{Y}])= \begin{cases}
        Y &\text{if $\frac{1}{\left | Y \right |}\overrightarrow{Y}$,
        $\frac{1}{\left | X_{0} \right |}\overrightarrow{X_{0}}$
            are in the same hemisphere or $Y=O_{3}$}\\
        (-\frac{2\pi}{\left | Y \right |}+1)Y
        &\text{if $\frac{1}{\left | Y \right |}\overrightarrow{Y}$,
           $\frac{1}{\left | X_{0} \right |}\overrightarrow{X_{0}}$
           are in different hemispheres and $ Y \neq O_{3}$},\\
        \end{cases}
\end{equation}
where points $0$ and $2\pi$ are identified, and $[\overrightarrow{Y}]$ is an equivalence class of $D$ (see Appendix
\ref{S:app1_BCH} for definition of $D$). Then $e^{q([Y])}=e^{[Y]}$ for any $[Y] \in D$.

\begin{theorem}\cite{Co}\label{thm:FrequencyGeneral1}
Suppose that the hypotheses of Theorem \ref{thm:Hopf_RW} are satisfied. Then, there exists a unique sufficiently $CS$
branch $X(\lambda)$ such that $\frac{2\pi}{\left | \omega_{\lambda}\right |}X(\lambda) \in q(D)$ for $\lambda \geq 0$,
$\overrightarrow {X(\lambda)}$ is a primary frequency vector of $\mathbf{\Phi}(t, u_{\lambda},\lambda)$ for $\lambda > 0$
small and $\mathbf{\Phi}(\frac{2\pi}{\left | \omega_{bif}\right |}, u_{0},0)= e^{X(0)\frac{2\pi}{\left | \omega_{bif}
\right |}}u_{0}$.
\end{theorem}
The branch $X(\lambda)$ does not have the property that $X(0)=X_{0}$. The next corollary shows that we can construct a
branch of primary frequency vectors $\overrightarrow{X^{f}(\lambda)}$ such that $X^{f}(0)=X_{0}$.
\begin{corollary}\cite{Co}\label{cor:FrequencyGeneral2}
Suppose that the hypotheses of Theorem \ref{thm:Hopf_RW} are satisfied. Then, there exists a branch $X^{f}(\lambda)$
for $\lambda \geq 0$ such that $\overrightarrow {X^{f}(\lambda)}$ is a primary frequency vector of
$\mathbf{\Phi}(t, u_{\lambda}, \lambda)$ for $\lambda > 0$ small and $X^{f}(0)= X_{0}$. Moreover,
\begin{enumerate}
\item if  $\left | X_{0} \right | \neq k\omega_{bif}$ for all $k \in \mathbb{Z}$, then the branch $ {X^{f}(\lambda)}$
for $\lambda \geq 0$ is sufficiently $CS$ and
\begin{equation}\label{E:freqTaylor01}
X^{f}(\lambda)= \left [\left | X_{0} \right|+O(\sqrt{\lambda})\right ]X_{0}^{1}+O(\sqrt{\lambda})X_{1}+
O(\sqrt{\lambda})X_{2} \mbox{ for small } \lambda \geq 0.
\end{equation}
(The branch $\mathbf{\Phi}(t,u_{\lambda},\lambda)$ contains only modulated rotating waves for $\lambda \geq 0$ small).
\item if $\left | X_{0} \right | = k \omega_{bif}$ for some $k \in \mathbb{Z}$, $k \neq 0$, then
$\left | X^{f}(\lambda) \right | $ for $\lambda \geq 0$ is continuous and
\begin{equation}\label{E:freqTaylor02}
\left | X^{f}(\lambda) \right | = \left | X_{0} \right | +O(\lambda^{\frac{1}{4}}) \mbox{ for small } \lambda \geq 0.
\end{equation}
\end{enumerate}
\end{corollary}
Therefore, we can not have frequency -locking phenomena for the modulated rotating waves  $\mathbf{\Phi}(t,u_{\lambda},\lambda)$.
If $\left | X_{0} \right | = k \omega_{bif}$ for some $ k \in \mathbb{Z}$, $k \neq 0$, it may be possible that
the branch $X^{f}(\lambda)$ be discontinuous at $\lambda= 0$ or/and at any $\lambda > 0$ such that
$\left | X(\lambda) \right | =0$.

\subsection{Construction of the branches $X(\lambda)$ and $X^{f}(\lambda)$}\label{SS:Construction}

In the rest of the section we give a way of constructing $X(\lambda)$ and $X^{f}(\lambda)$ for $\lambda \geq 0$
using the $BCH$ formula in $so(3)$ presented in Theorem \ref{thm:BCH_formula}.
Theorems \ref{thm:Aux1}, \ref{cor:Aux2} and \ref{thm:GlobalZeta} are used to construct the branch $X(\lambda)$.
Recall that a primary frequency vector $\overrightarrow{X(\lambda)}$ of $\mathbf{\Phi}(t, u_{\lambda}, \lambda)$ is
given by $A(\frac{2\pi}{\left | \omega_{\lambda}\right |}, \lambda)=e^{X(\lambda)\frac{2\pi}{\left | \omega_{\lambda}
\right | }}$ for $\lambda > 0$ small, we show that there exists a sufficiently $CS$ function $Z(t, \lambda) \in q(D)$
such that $A(t,\lambda)= e^{Z(t, \lambda)}$ for at least $t \in \left [0,\frac{2\pi}{\left |
\omega_{\lambda}\right | } \right ]$ and $\lambda \geq 0$ small.
\begin{theorem}\label{thm:Aux1}
Suppose that the hypotheses of Theorem \ref{thm:Hopf_RW} are satisfied.
Let us consider the following initial value problem
\begin{equation}\label{E:ODEs_frequency}
\begin{array}{rll}
\overrightarrow {\dot{Z}} &=& \left [I_{3}+\frac{1}{2}Z+ \left (\frac{1}{\left | Z \right | ^{2}}-
\frac{\cos \frac{\left | Z \right |}{2}} {2 \sin\frac{\left | Z \right |}{2}\left| Z \right |} \right )Z^{2} \right ]
\overrightarrow {X^{G}(t, \lambda)},\\
Z(0) &=& O_{3}.
\end{array}
\end{equation}
Then,
\begin{enumerate}
\item there exists a positive integer $n$ independent of $\lambda $ such that the initial value problem
$(\ref{E:ODEs_frequency})$ has a unique sufficiently $CS$ solution $Z_{1}(t, \lambda)$ on $t \in
[0, \frac{2\pi}{n\left | \omega_{\lambda}\right | }]$ and $\lambda \geq 0$ small;
\item
\begin{equation}\label{E:Zeta0}
Z_{1}(t,\lambda)= \left ( \left | X_{0} \right |t+ \sqrt{\lambda}
\int_{0}^{t}x_{01}(s)\, ds+ \lambda x_{03}(t, \lambda) \right)X_{0}^{1}+
\sqrt{\lambda} z_{1}(t, \lambda)X_{1}+\sqrt{\lambda} z_{2}(t, \lambda)X_{2}
\end{equation}
for any $t \in [0, \frac{2\pi}{n\left | \omega_{\lambda}\right |}]$ and $\lambda \geq 0$ small;
\end{enumerate}
\end{theorem}
The integer $n$ that appears in Theorem \ref{thm:Aux1} can be found and, in general, it is not $1$. Therefore,
in general, the function $Z_{1}(t, \lambda)$ is not defined on the entire interval
$[0, \frac{2\pi}{\left | \omega_{\lambda} \right |}]$, for $\lambda \geq 0$ small.
\begin{corollary}\label{cor:Aux2}
Suppose that the hypotheses of Theorem \ref{thm:Hopf_RW} are satisfied and the positive integer $n$ is the one obtained in Theorem \ref{thm:Aux1}.
Then, for any $i= 1$,~$2$, \ldots,~$n-1$ the initial value problem
\begin{equation}\label{E:ODEs_frequency1}
\begin{array}{rll}
\overrightarrow {\dot{Z}} &=& \left [I_{3}+\frac{1}{2}Z+ \left(\frac{1}{\left | Z \right | ^{2}}-\frac{\cos \frac{\left | Z\right |}{2}} {2\sin \frac{\left | Z \right |}{2}\left| Z \right|} \right )Z^{2}\right ]
\overrightarrow {X^{G}(t, \lambda)},\\
Z(i\frac{2\pi}{n \left | \omega_{\lambda}\right | }) &=& O_{3}
\end{array}
\end{equation}
has a unique sufficiently $CS$ solution $Z_{i+1}(t, \lambda)$ on $[i\frac{2\pi}{n\left |\omega_{\lambda}\right | },
(i+1)\frac{2\pi}{n\left |\omega_{\lambda}\right | }]$ and $\lambda \geq 0$ small.
\end{corollary}
Combining Theorem \ref{thm:Aux1} and Corollary \ref{cor:Aux2}, we obtain the following result:
\begin{theorem}\label{thm:GlobalZeta}
Suppose that the hypotheses of Theorem \ref{thm:Hopf_RW} are satisfied. Then, there exists a sufficiently $CS$ function
$Z(t, \lambda) \in q(D)$ such that $A(t,\lambda)= e^{Z(t, \lambda)}$ for $t \in [0,\infty)$
and $\lambda \geq 0$ small, where $A(t,\lambda)$ is the solution of the initial value problem (\ref{E:AuxFreq}).
In fact, the function $Z(t, \lambda)$ satisfies the following initial value problem for $t \in [0,\infty)$ and $\lambda \geq 0$ small :
\begin{equation}\label{E:ODEs_frequency_final}
\begin{array}{rll}
 \overrightarrow {\dot{Z}} &=& \left [I_{3}+\frac{1}{2}Z+ \left (\frac{1}{\left | Z \right | ^{2}}-
\frac{\cos \frac{\left | Z \right |}{2}} {2 \sin\frac{\left | Z \right |}{2}\left| Z \right |} \right )Z^{2} \right ]
\overrightarrow {X^{G}(t, \lambda)} (\mbox{ mod } 2\pi),\\
\overrightarrow{ Z(0)} &=& \overrightarrow{O_{3}}.
\end{array}
\end{equation}
\end{theorem}
\begin{remark}
In fact, the function $Z(t, \lambda)$ obtained in Theorem \ref{thm:GlobalZeta} is constructed as follows:
\begin{enumerate}
\item On the interval $[0, \frac{2\pi}{\omega_{\lambda}}]$:
\begin{equation}\label{E:ZetaFormula00}
Z(t, \lambda)= Z^{0}(t, \lambda)= q(BCH(Z_{1}(\frac{2\pi}{n\left | \omega_{\lambda}\right |}, \lambda),
Z_{2}(2\frac{2\pi}{n \left | \omega_{\lambda}\right | }, \lambda), \ldots, Z_{i}(t, \lambda)))
\end{equation}
for all $t \in [(i-1)\frac{2\pi}{n\left | \omega_{\lambda}\right |}, i\frac{2\pi}{n\left | \omega_{\lambda}\right | }]$
and $\lambda \geq 0$ small, where $i= 1$,~$2$, \ldots, ~$n$ and $Z_{i}$ is the solution of the initial
value problem (\ref{E:ODEs_frequency1}) on the interval $[(i-1)\frac{2\pi}{n\left | \omega_{\lambda}\right |},
i\frac{2\pi}{n\left | \omega_{\lambda}\right | }]$.
\item On the interval $[i\frac{2\pi}{\left |\omega_{\lambda}\right | },
(i+1)\frac{2\pi}{\left |\omega_{\lambda}\right | }]$ for any integer $i \geq 1$ and $\lambda \geq 0$ small :
\begin{equation}\label{E:ZetaFormula1}
Z(t,\lambda)= Z^{i}(t, \lambda)= q(BCH(Z^{0}(\frac{2\pi} {\left| \omega_{\lambda}\right | }, \lambda),
Z^{i-1}(t-\frac{2\pi}{\left | \omega_{\lambda}\right |},\lambda))).
\end{equation}
\end{enumerate}
\end{remark}


Using the function $Z(t, \lambda)$ obtained in Theorem \ref{thm:GlobalZeta}, we construct the branches
$X(\lambda)$ ($\lambda \geq 0$) as follows:
\begin{proposition}[Primary Frequency Vectors Associated to $\mathbf{\Phi}(t,u_{\lambda},\lambda)$]
\label{prop:Frequency1}
Suppose that the hypotheses of Theorem \ref{thm:Hopf_RW} are satisfied. Then, the branch $X(\lambda)$ for
$\lambda \geq 0$ from Theorem \ref{thm:FrequencyGeneral1} is defined by:
\begin{equation}\label{E:FrequencyFormula}
X(\lambda) =\frac{\left | \omega_{\lambda} \right |}{2\pi}Z(\frac{2\pi}{\left | \omega_{\lambda} \right | }, \lambda)
\mbox{ for small } \lambda \geq 0,
\end{equation}
where $Z(t, \lambda)$ is given in Theorem $\ref{thm:GlobalZeta}$;
\end{proposition}
\begin{corollary}[Primary Frequency Vectors Associated to $\mathbf{\Phi}(t,u_{\lambda},\lambda)$]\label{cor:Frequency2}
Then, the branch $X^{f}(\lambda)$ for $\lambda \geq 0$ from Corollary \ref{cor:FrequencyGeneral2} is
defined by:
\begin{enumerate}
\item if $\left | X_{0} \right | \neq k\omega_{bif}$ for all $k \in \mathbb{Z}$, then
we define $X^{f}(\lambda)=\frac{\left | X(\lambda) \right |+k\left |\omega_{\lambda}\right | }{\left | X(\lambda)
\right | }X(\lambda)$ for $\lambda > 0$ small.\\
\item if $\left | X_{0}\right | = k\omega_{bif}$ for some $ k \in \mathbb{Z}$, $k \neq 0$, then we define:\\
$X^{f}(0)=X_{0}$;\\
if $\left | X(\lambda) \right | \neq 0$, then we define $X^{f}(\lambda)= \frac{\left | X(\lambda) \right |+
 k\left |\omega_{\lambda}\right |}{\left | X(\lambda) \right | }X(\lambda)$ for $\lambda \geq 0$ small;\\
if $\left | X(\lambda) \right | = 0$, then we define $X^{f}(\lambda)= X(\lambda)+ k\left | \omega_{\lambda}\right |
Q(\lambda)$, where $Q(\lambda) \in so(3)$, $\left | Q(\lambda)\right | =1$ for $\lambda > 0$ small.
\end{enumerate}
\end{corollary}
Using Lemma \ref{lem:Hopf_group} for $\lambda >0$ small, we can associate  $A(t, \lambda)= e^{X(\lambda)t}B(t, \lambda)=
e^{X^{f}(\lambda)t}B^{f}(t, \lambda)$ to each modulated rotating wave or periodic solution of period $\frac{2\pi}
{\left | \omega_{\lambda}\right |}$ obtained by Theorem \ref{thm:Hopf_RW}, where $B(t, \lambda)$, $B^{f}(t, \lambda)$
are $\frac{2\pi}{\left | \omega_{\lambda} \right | }$-periodic in $t$ for
$\lambda > 0$ small, as well as $B(0,\lambda)=B^{f}(0, \lambda)=I_{3}$. Also, let $B(t,0)$ and $B^{f}(t,0)$ be such that $A(t,0)=e^{X(0)t}B(t, 0)= e^{X^{f}(0)t}B^{f}(t, 0)$, that is $B^{f}(t,0)=I_{3}$.
\begin{proposition}[Periodic Parts Associated to $\mathbf{\Phi}(t,u_{\lambda},\lambda)$]\label{thm:Periodic_part}
Suppose that the hypotheses of Theorem \ref{thm:Hopf_RW} are satisfied. Then,
\begin{enumerate}
\item if $\left | X_{0} \right | \neq k\omega_{bif}$ for all $k \in \mathbb{Z}$, we have
\begin{equation}\label{E:PeriodicPartFormula1}
B^{f}(t, \lambda)= e^{\sqrt{\lambda}Y(t,\lambda)} \mbox{ for  } t \in [0, \infty) \mbox { and } \lambda \geq 0 \mbox{ small, }
\end{equation}
where $Y(t, \lambda)$ is $\frac{2\pi}{\left |\omega_{\lambda} \right | }$-periodic in $t$ and sufficiently $CS$ for $t \in [0, \infty)$ and $\lambda \geq 0$ small.\\
\item if $\left | X_{0} \right | = k\omega_{bif}$ for some $k \in \mathbb{Z}$, $k \neq 0$, we have
\begin{equation}\label{E:PeriodicPartFormula2}
B(t, \lambda)= e^{X_{0}t+\sqrt{\lambda}H(t,\lambda)} \mbox { for } t \in [0, \infty) \mbox { and } \lambda \geq 0 \mbox{ small, }
\end{equation}
where $H(t, \lambda)$ is sufficiently $CS$ for $t \in [0, \infty)$ and $\lambda \geq 0$ small, and
$e^{X_{0}\frac{2\pi}{\left | \omega_{\lambda}\right |}+\sqrt{\lambda}H(\frac{2\pi}{\left | \omega_{\lambda}\right |},\lambda)}=I_{3}$ for $\lambda \geq 0$ small.
\end{enumerate}
\end{proposition}
\begin{remark}
In fact, $B^{f}(t, \lambda)= e^{q(BCH(-X^{f}(\lambda)t, Z(t,\lambda)))}$ (respectively
$B(t, \lambda)= e^{q(BCH(-X(\lambda)t, Z(t, \lambda)))}$) for $t \in [0, \infty)$ and $\lambda \geq 0$ small, where $Z(t, \lambda)$ is given in Theorem \ref{thm:GlobalZeta}.
\end{remark}

Using Theorems \ref{prop:Frequency1}, \ref{cor:Frequency2} and \ref{thm:Periodic_part} we get the following result:
\begin{theorem}\label{thm:MRW}
Suppose that the hypotheses of Theorem \ref{thm:Hopf_RW} are satisfied. Then,
\begin{enumerate}
\item
if $\left | X_{0} \right | \neq k\omega_{bif}$ for all $k \in \mathbb{Z}$, we have
\[ \mathbf{\Phi}(t,u_{\lambda},\lambda)= e^{\left [ (\left | X_{0} \right|+O(\sqrt{\lambda}))X_{0}^{1}+
O(\sqrt{\lambda})X_{1}+O(\sqrt{\lambda})X_{2}\right ]t}e^{\sqrt{\lambda}Y(t,\lambda)}\Psi(\sqrt{\lambda}r(t,\lambda)),\]
where $Y(t,\lambda)$ and $r(t, \lambda)$ are sufficiently $CS$ and $\frac{2\pi}{\left | \omega_{\lambda} \right | }$-periodic in $t$ for $\lambda \geq 0$ small.
\item
if $\left | X_{0} \right | = k\omega_{bif}$ for some $k \in \mathbb{Z}$, $k \neq 0$, we have
\[ \mathbf{\Phi}(t,u_{\lambda},\lambda)=e^{\left [ O(\sqrt{\lambda})X_{0}^{1}+O(\sqrt{\lambda})X_{1}+
O(\sqrt{\lambda})X_{2}\right ] t}e^{X_{0}t+\sqrt{\lambda}H(t,\lambda)}\Psi(\sqrt{\lambda}r(t,\lambda)),\]
where $r(t, \lambda)$ and $H(t, \lambda)$ are sufficiently $CS$ and $r(t, \lambda)$ is
$\frac{2\pi}{\left | \omega_{\lambda} \right | }$-periodic in $t$ for $\lambda \geq 0$ small, as well as
$e^{X_{0}\frac{2\pi}{\left | \omega_{\lambda} \right | }+\sqrt{\lambda}H(\frac{2\pi}{\left | \omega_{\lambda}\right |},
\lambda)}=I_{3}$ for $\lambda \geq 0$ small.
\end{enumerate}
\end{theorem}
In fact, we have that the initial value problem (\ref{E:reduceddiffeqparam1}) is equivalent with the following initial value problem
\begin{equation}\label{E:Zeta_final}
\begin{array}{rll}
\overrightarrow {\dot{Z}} &=& \left [I_{3}+\frac{1}{2}Z+ \left (\frac{1}{\left | Z \right | ^{2}}-
\frac{\cos \frac{\left | Z \right |}{2}} {2 \sin\frac{\left | Z \right |}{2}\left| Z \right |} \right )Z^{2} \right ]
\overrightarrow {X_{G}(q, \lambda)} \left ( \mbox{ mod } 2\pi \right ),\\
\dot{q} &=& X_{N}(q, \lambda),\\
\overrightarrow{Z(0)} &=& \overrightarrow{O_{3}}.
\end{array}
\end{equation}
We could have considered that the differential equation $\dot{q}= X_{N}(q, \lambda)$ has a
$\frac{2\pi}{\left | \omega_{\lambda}\right | }$-periodic solution $q(t, \lambda)=O(\sqrt{\lambda})$
instead of assuming that a supercritical Hopf bifurcation takes places at $q=0$ for $\lambda=0$. All the results remain
valid.

\section{Examples}\label{S:Examples_Transition}

Throughout this section, we let $\epsilon= \sqrt{\lambda}$. Here we present some examples of functions
$X^{G}(t, \lambda)$ for which we can get the closed form solution $A(t, \lambda)$ of the initial value problem
(\ref{E:AuxFreq}). They obey Theorem \ref{thm:MRW}.
\begin{example}[ Hopf bifurcation to modulated rotating waves]\label{ex:exHoph1}
Let $g(t,\lambda)$ be a sufficiently smooth periodic function of period $\frac{2\pi}{\left | \omega_{bif}+\lambda \right | }$ such that $g(0, \lambda)= 0$ for $\lambda \geq 0$ small. For any $t \in [0,\infty)$ and $\lambda
\geq 0$ small, we define
\[X^{G}(t,\lambda)= (2 \epsilon X_{1} + 2 \epsilon X_{2}+2 \epsilon X_{0}^{1})\dot{g}(t,\lambda)+
e^{-(2\epsilon X_{1}+2\epsilon X_{2}+ \epsilon X_{0}^{1})g(t,\lambda)}(X_{0}+
\epsilon X_{1})e^{(2\epsilon X_{1}+2\epsilon X_{2}+2\epsilon X_{0}^{1})g(t, \lambda)}.\]
It is clear that $X^{G}(t,\lambda)$ is a sufficiently $CS$, $\frac{2\pi}{\left | \omega_{bif}+\lambda \right |}$-periodic function such that $X^{G}(t,0)= X_{0}$. Then, the initial value problem (\ref{E:AuxFreq}) has the solution
\[A(t,\lambda)= e^{(X_{0}+\epsilon X_{1})t}e^{(2 \epsilon X_{1}+2 \epsilon X_{2}+2 \epsilon X_{0}^{1})g(t,\lambda)}.\]
Clearly, $\overrightarrow X_{0}+\epsilon \overrightarrow X_{1}$ is not orthogonal to $\overrightarrow X_{0}$.
\end{example}

\begin{example}[Example 1. Resonant drift phenomena for modulated rotating waves]\label{ex:exHoph2}
Let $\omega_{bif}= \left | X_{0} \right | $ and $g(t,\lambda)$ be a sufficiently smooth periodic function of period
$\frac{2\pi}{\left | \omega_{bif}+\lambda \right | }$ such that $g(0,\lambda)=0$ for $\lambda \geq 0$ small. \\
For any $t \in [0,\infty)$ and $\lambda \geq 0$ small, we define
\begin{multline}
X^{G}(t, \lambda)= (X_{0}+\epsilon X_{2})\left (\frac{\left |\omega_{bif}+\lambda \right |}{\sqrt{\left | X _{0} \right | ^{2}+\lambda}}+\lambda \dot{g}(t,\lambda)\right )+ \epsilon e^{-(X_{0}+\epsilon X_{2})\left (\frac{\left | \omega_{bif}+ \lambda \right |}{\sqrt{\left | X _{0} \right | ^{2}+\lambda}}t+ \lambda g(t,\lambda)\right )}X_{1}\\
\cdot e^{(X_{0}+\epsilon X_{2})\left (\frac{\left | \omega_{bif}+\lambda \right |}{\sqrt{\left | X _{0} \right | ^{2}+\lambda}}t+\lambda g(t,\lambda)\right )}.
\end{multline}
It is clear that $X^{G}(t,\lambda)$ is a sufficiently $CS$, $\frac{2\pi}{\left | \omega_{bif}+\lambda  \right |}$- periodic function such that $X^{G}(t,0)= X_{0}$. Then, the initial value problem (\ref{E:AuxFreq}) has the solution
\begin{equation}
A(t,\lambda)= e^{\epsilon X_{1}t}e^{(X_{0}+\epsilon X_{2})\left (\frac{\left | \omega_{bif}+\lambda \right |} {\sqrt{\left | X_{0} \right | ^{2}+ \lambda}}t+\lambda g(t,\lambda)\right )}.
\end{equation}
Clearly, $\epsilon \overrightarrow X_{1}$ is orthogonal to $\overrightarrow X_{0}$.
\end{example}

\begin{example}[Example 2. Resonant drift phenomena for modulated rotating waves]\label{ex:exHoph3}
Let $\omega_{bif}= \left | X_{0} \right | $ and $g(t,\lambda)$ be a sufficiently smooth periodic function of period
$\frac{2\pi}{\left | \omega_{bif}+\lambda \right | }$ such that $g(0, \lambda)=0$ for $\lambda \geq 0$ small. \\
For any $t \in [0,\infty)$ and $\lambda \geq 0$ small, we define
\begin{multline}
X^{G}(t, \lambda)= (X_{0}+\epsilon X_{2})\left (\frac{\left | \omega_{bif}+\lambda \right |}{\sqrt{\left | X _{0} \right |^{2}+\lambda}}+ \lambda \dot{g}(t,\lambda)\right)+ \epsilon e^{-(X_{0}+\epsilon X_{2})\left (\frac{\left |\omega_{bif}+\lambda \right |}{\sqrt{\left | X _{0} \right | ^{2}+\lambda}}t+\lambda g(t,\lambda)\right )}\\
\cdot (X_{0}+X_{1})e^{(X_{0}+\epsilon X_{2})\left (\frac{\left | \omega_{bif}+ \lambda \right |}{\sqrt{\left | X _{0} \right |^{2}+ \lambda}}t+\lambda g(t,\lambda)\right )}.
\end{multline}
It is clear that $X^{G}(t,\lambda)$ is a sufficiently $CS$, $\frac{2\pi}{\left | \omega_{bif}+\lambda \right |}$-periodic function such that $X^{G}(t,0)= X_{0}$. Then, the initial value problem (\ref{E:AuxFreq}) has the solution
\begin{equation}
A(t,\lambda)= e^{\epsilon(X_{0}+X_{1})t}e^{(X_{0}+\epsilon X_{2})\left (\frac{\left | \omega_{bif}+\lambda \right |}
{\sqrt{\left | X _{0} \right | ^{2}+\lambda}}t+ \lambda g(t,\lambda)\right )}.
\end{equation}
Clearly, $\epsilon(\overrightarrow X_{0}+ \overrightarrow X_{1})$ is not orthogonal to $\overrightarrow X_{0}$.
\end{example}

\begin{example}[Example 3. Resonant drift phenomena for modulated rotating waves]\label{ex:exHoph4}
Let $\omega_{bif}= \frac{1}{k}\left | X_{0} \right |$ for $ k \in \mathbb{Z}$, $k \neq 0$ and $g(t, \lambda, \mu)$ be a sufficiently smooth periodic function of period $\frac{2\pi}{\left | \omega_{bif}+\lambda \right | }$ such that $g(0, \lambda, \mu )=0$ for $ \lambda \geq 0$, $\left | \mu \right | $ small. \\
For any $t \in [0,\infty)$, $\lambda \geq 0$ small and $\left | \mu \right | $ small, we define
\begin{multline}
X^{G}(t, \lambda,\mu)=(X_{0}+\mu X_{1})\left (\frac{k \left | \omega_{bif}+\lambda \right |}{\sqrt{\left | X _{0} \right | ^{2}+\mu^{2}}}+\lambda \dot{g}(t,\lambda,\mu)\right ) \\
+\epsilon^{k} e^{-(X_{0}+\mu X_{1}) \left (\frac{k \left | \omega_{bif}+\lambda \right |} {\sqrt{\left | X _{0} \right | ^{2}+\mu^{2}}}t+\lambda g(t,\lambda,\mu)\right )}((\epsilon-\mu)X_{0}+X_{1}+ X_{2})e^{(X_{0}+ \mu X_{1})
\left (\frac{k \left | \omega_{bif}+\lambda \right |}{\sqrt{\left | X _{0} \right | ^{2}+ \mu^{2}}}t+\lambda g(t,\lambda,\mu)\right )}.
\end{multline}
It is clear that $X^{G}(t,\lambda,\mu)$ is a sufficiently $CS_{\lambda}$, $\frac{2\pi}{\left | \omega_{bif}+\lambda \right |}$ -periodic function such that $X^{G}(t,0,\mu)= X_{0}+\mu X_{1}$. Then, the initial value problem (\ref{E:AuxFreq}) for two parameters $\lambda$, $\mu$ has the solution
\begin{equation}
A(t, \lambda, \mu)= e^{\epsilon^{k}((\epsilon-\mu)X_{0}+X_{1}+ X_{2})t}e^{(X_{0}+\mu X_{1})\left (\frac{k\left | \omega_{bif}+\lambda \right |}{\sqrt{\left | X _{0} \right | ^{2}+\mu^{2}}}t+\lambda g(t, \lambda, \mu)\right )}.
\end{equation}
Clearly, for $ \mu=\sqrt{\lambda}$, we have orthogonality, that is $\epsilon ^{k}(\overrightarrow X_{1}+\overrightarrow X_{2})$ is orthogonal to $\overrightarrow X_{0}$.
\end{example}

\begin{example}[Nonuniform rigidly rotation]\label{ex:exHopf5}
Let $\omega_{bif}= \left | X_{0} \right | $ and $g(t,\lambda)$ be a sufficiently smooth periodic function of period
$\frac{2\pi}{\left | \omega_{bif}+\lambda \right | }$ such that $g(0, \lambda)=0$ for $\lambda \geq 0$ small. \\
For any $t \in [0,\infty)$ and $\lambda \geq 0$ small, we define
\begin{equation}
X^{G}(t, \lambda)= (X_{0}+\epsilon X_{1})(1+\epsilon \dot{g}(t, \lambda)).
\end{equation}
It is clear that $X^{G}(t,\lambda)$ is a sufficiently $CS$, $\frac{2\pi}{\left | \omega_{bif}+\lambda \right |}$-periodic function such that $X^{G}(t,0)= X_{0}$. Then, the initial value problem (\ref{E:AuxFreq}) has the solution
\begin{equation}
A(t, \lambda)=e^{(X_{0}+\epsilon X_{1})(t+g(t, \lambda))}.
\end{equation}
This clearly represents a non-uniformly rigid rotation about the line containing the vector
$\overrightarrow{X_{0}+\epsilon X_{1}}$. This is a degenerate situation for the case of one parameter.
\end{example}
Also, the case of a $\frac{2\pi}{\left | \omega_{\lambda} \right |}$-periodic solution is a degenerate situation for
one parameter. Therefore, we do not present it here.\\
These examples will be used in Section \ref{S:Numerical}.
As we have seen in Examples \ref{ex:exHoph2} and \ref{ex:exHoph3}, if we have a single parameter $\lambda$ and $\left | X_{0} \right | = k \omega_{bif}$ for some $k \in \mathbb{Z}$, we can have three different types of solutions for the initial value problem (\ref{E:AuxFreq}). In Example(\ref{ex:exHoph4}), we have two parameters $\lambda \geq 0$ and $\mu$, as well as $\left | X_{0} \right | = k \omega_{bif}$ for some $k \in \mathbb{Z}$ and we have a branch
$\mu=\mu(\lambda)$ for which $\overrightarrow{X(\lambda, \mu(\lambda))}$ is orthogonal to $\overrightarrow{X_{0}}$.
In next section, we prove that this is generically valid.

\section{Resonant Drift Phenomena for Modulated Rotating Waves on $r\mathbf{S^{2}}$ }\label{S:Res_drift}

We consider the following reaction-diffusion system
\begin{equation}\label{E:RDparams}
\frac{\partial{u}}{\partial{t}}(t,x)= D\Delta_{S}{u}(t,x)+F(u(t,x),\lambda, \mu) \mbox{ on } r\mathbf{S^{2}},
\end{equation}
where $u=(u_{1}$,~$u_{2}$, \ldots, ~$u_{N})\colon \mathbb{R} \times r\mathbf{S^{2}} \rightarrow \mathbb{R}^{N}$ with $N \geq 1$, $D=\left(\begin{array}{ccc}
            d_{1}  & \dots  & 0 \\
            \vdots & \ddots & \vdots\\
            0      & \dots  & d_{N}
             \end{array}
       \right)$ with  $d_{i} \geq 0$ for $i= 1$,~$2$, \ldots, ~$N$ are the diffusion coefficients and
$F=(F_{1}$,~$F_{2}$, \ldots,~$F_{N})\colon \mathbb{R}^{N} \times \mathbb{R} \times \mathbb{R}\rightarrow \mathbb{R}^{N}$ are sufficiently smooth functions such that $F(0,\lambda,\mu)=0$ for $\left | \lambda \right | $, $\left | \mu \right | $ small.\\
We study the reaction-diffusion system (\ref{E:RDparams}) on the function space $\mathbf{Y}$ defined in Section
\ref{S:RDS}.
Let $\mathbf{\Phi}(t,u,\lambda,\mu)$ be the $SO(3)$-equivariant sufficiently smooth local semiflow defined as in
Section \ref{S:RDS}. Let $u_{0} \in \mathbf{Y}^{\alpha}$ be a relative equilibrium that is not an equilibrium for (\ref{E:RDparams}) at $(\lambda, \mu)= (0, 0)$ and  such that the stabilizer of $u_{0}$ is $\Sigma_{u_{0}}= {I_{3}}$. consider $L$ the linearization of the right-hand side of (\ref{E:RDparams}) with respect to the rotating wave $\mathbf{\Phi}(t, u_{0}, 0, 0)= e^{X_{0}t}u_{0}$ at $(\lambda,\mu)=(0, 0)$ in the co-rotating frame, that is
\[L =D\Delta_{S}+D_{u}F(u_{0}, 0, 0)-X_{0}.\]
Suppose that:
\begin{enumerate}
\item $\sigma(L) \cap \{z \in \mathbb{C} \mid Re \, (z) \geq 0 \}$ is a spectral set with spectral projection $P_{*}$, and dim$(R(P_{*})) < \infty;$ \item the semigroup $e^{Lt}$ satisfies
$\left | e^{Lt}|_{R(1-P_{*})} \right | \leq Ce^{-\beta_{0} t}$ for some $\beta_{0} > 0$ and $C > 0$.
\end{enumerate}
Theorems \ref{thm:CMR_REparam} with parameters $\lambda$, $ \mu$ can be applied.\\
There exist sufficiently smooth functions $X_{G} \colon V_{*} \times \mathbb{R} \times \mathbb{R} \rightarrow so(3)$ and $X_{N} \colon V_{*} \times \mathbb{R} \times \mathbb{R} \rightarrow V_{*}$ such that any solution of
\begin{equation}\label{E:reduceddiffeqparam2}
    \begin{array}{lll}
    \dot{A} &=& AX_{G}(q, \lambda, \mu),\\
    \dot{q} &=& X_{N}(q, \lambda, \mu),\\
    \end{array}
\end{equation}
on $SO(3) \times V_{*}$ corresponds to a solution of the reaction-diffusion system (\ref{E:RDparam}) on
$M_{u_{0}}^{cu}(\lambda, \mu)$ under the diffeomorphic identification for $\left | \lambda  \right |$ and
$\left | \mu  \right |$ small. Also, $X_{G}(0, 0, 0)= X_{0}$, $X_{N}(0, 0, 0)= 0$ and $\sigma(D_{u}X_{N}(0, 0, 0))= \sigma(Q_{*}L|_{V_{*}})$, where $Q_{*}$ is the projection onto $V_{*}$ along $T_{u_{0}}(SO(3)u_{0})$.
Let dim$R(P_{*})=5$. Suppose that a supercritical Hopf bifurcation with eigenvalues $\pm i \omega_{bif}$ takes place
in the second differential equation of (\ref{E:reduceddiffeqparam2}) in $V_{*}$ at $q= 0$ for $(\lambda, \mu)= (0,0)$,
that is:
\begin{enumerate}
\item $X_{N}(0, 0, 0)= 0$;
\item $D_{q}X_{N}(0, 0, 0)$ has eigenvalues $\pm i \omega_{bif}$; without loss of generality, we
assume that $X_{N}(0, \lambda, \mu)= 0$ for $\left | \lambda \right | $, $\left | \mu \right | $ small;
\item $D_{q}X_{N}(0,\lambda, \mu)$ has the eigenvalues $\alpha(\lambda, \mu) \pm i(\omega_{bif}+
\beta(\lambda, \mu)$ with $\alpha(0, 0)= \beta(0,0)= 0$) and \, $\alpha_{\lambda}(0,0)> 0$.
This implies that $\alpha(\lambda_{H}(\mu), \mu)= 0$ for some sufficiently smooth curve
$\lambda= \lambda_{H}(\mu)$ with $\lambda_{H}(0)= 0$. This curve represents the Hopf points and without
loss of generality, we suppose that $\lambda_{H}(\mu)= 0$ for $\left | \mu \right | $ small;
\item the branch of periodic solutions $q(t, \lambda, \mu)$ bifurcating from $q=0$ generically satisfies
$q(t, \lambda, \mu)=\sqrt(\lambda)r(t, \lambda, \mu)$.
\end{enumerate}

\begin{figure}
\begin{center}
\psfig{file=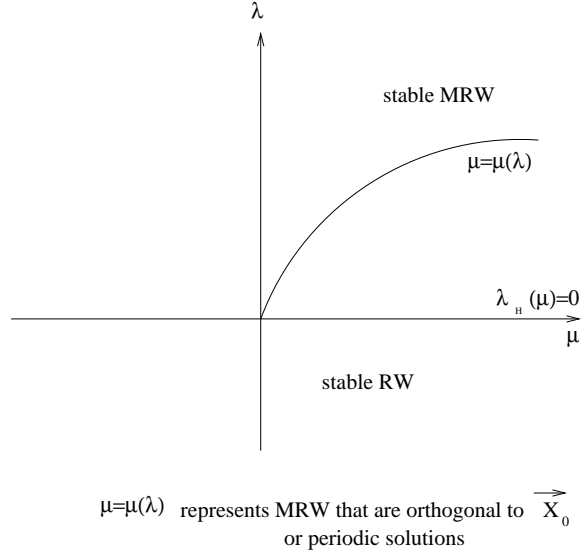,width=3in}
\end{center}
\caption{Parameter Space for Resonant Hopf Bifurcation}
\label{figure:figura2}
\end{figure}

For $\lambda >0$ small and $\left | \mu \right | $ small, let  $T(\lambda, \mu)= \frac{2\pi}{\left | \omega_{\lambda,\mu}
\right | }$ be the period of the periodic solution $q(t, \lambda,\mu)$ near $q=0$, that appears due to the supercritical Hopf bifurcation, where $\omega_{\lambda, \mu}= \omega_{bif}+O(\mu)+\lambda s( \lambda, \mu)$   for $\lambda \geq 0$ small and $\left | \mu\right |$ small.\\
By Theorem \ref{thm:Hopf_RW}, there exists a sufficiently $CS$ branch $\mathbf{\Phi}(t,u_{\lambda,\mu},\lambda,\mu)$ for $\lambda \geq 0$ small and $\left |\mu \right | $ small such that for $\left | \mu \right | $ small,
$\mathbf{\Phi}(t,u_{0,\mu},0,\mu)=e^{X_{G}(0,0,\mu)t}u_{0,\mu}$ and such that for $ \lambda > 0$
small and $\left | \mu \right | $ small, $\mathbf{\Phi}(t,u_{\lambda,\mu},\lambda,\mu)$ is an orbitally
stable modulated rotating wave or periodic solution of period $\frac{2\pi}{\left | \omega_{\lambda,\mu}\right | }$. Let $\overrightarrow {X(\lambda, \mu)}$ be the sufficiently $CS$ branch such that for  $\left | \mu \right | $ small,
$e^{X(0,\mu)\frac{2\pi}{\left | \omega_{0,\mu} \right |}}= e^{X_{G}(0,0, \mu)\frac{2\pi}{\left | \omega_{0,\mu} \right |}}$, and for $\lambda > 0$ small and $\left | \mu \right | $ small, $\overrightarrow{X(\lambda,\mu)}$ is a
primary frequency vector corresponding to $\mathbf{\Phi}(t,u_{\lambda,\mu},\lambda,\mu)$.
We have the following two results:
\begin{theorem}[Resonance Case]\label{thm:Resonant_drift1}
Suppose the previous assumptions hold and that $\left | X_{0} \right | = k\omega_{bif}$ for some
$k \in \mathbb{Z}$, $k \neq 0$. Let \[X_{G}(q, \lambda, \mu)= x_{0}(q, \lambda,\mu)X_{0}^{1}+x_{1}(q, \lambda, \mu)X_{1}
+ x_{2}(q, \lambda,\mu)X_{2}.\]
Then, if $(x_{0})_{\mu}(0,0,0) \neq k(\omega_{0,\mu})^{'}|_{\mu=0}$, there exists a sufficiently $CS$ curve
$\mu= \mu(\lambda)$ for $\lambda \geq 0$ small such that $\mu(0)=0$, and there exists a sufficiently $CS$
branch of orbitally stable modulated rotating waves $\mathbf{\Phi}(t, u_{\lambda, \mu(\lambda)}, \lambda, \mu(\lambda))$
for $\lambda > 0$ small, with a primary frequency vector $\overrightarrow {X(\lambda, \mu(\lambda))}$ orthogonal to
$\overrightarrow {X_{0}}$.
\end{theorem}
\begin{corollary}\label{cor: Resonant_drift2}
Suppose that Theorem \ref{thm:Resonant_drift1} hold. Then,
\begin{enumerate}
\item the branch $X(\lambda, \mu(\lambda))$ is sufficiently $CS$ for $\lambda \geq 0$ small and
$ X(\lambda, \mu(\lambda))= O(\lambda)^{\frac{k}{2}})X_{1}+O(\lambda ^{\frac{k}{2}})X_{2}$ for
$\lambda \geq 0$ small.
\item $B(t, \lambda,\mu(\lambda))=e^{(O(\sqrt{\lambda})X_{0}^{1}+O(\sqrt{\lambda})X_{1}+
O(\sqrt{\lambda})X_{2})t}e^{X_{0}t+\sqrt{\lambda}H(t,\lambda)}$, for $ t \in [0,\infty)$
and $\lambda \geq 0$ small, where $H(t,\lambda)$ is a sufficiently $CS$ function such that
$e^{X_{0}\frac{2\pi}{\left | \omega_{\lambda} \right | }+\sqrt{\lambda}H(\frac{2\pi}
{\left | \omega_{\lambda}\right |}, \lambda)}=I_{3}$ for $\lambda \geq 0$ small.
\item there exists a branch $X^{f}(\lambda, \mu(\lambda))$ with $\lambda \geq 0$ that is discontinuous at
$\lambda= 0$, but $\left | X^{f}(\lambda, \mu(\lambda)) \right |$ is continuous for $\lambda \geq 0$ small and
$\left | X^{f}(\lambda, \mu(\lambda)) \right | = \left | X_{0} \right | + O(\lambda ^{\frac{k}{4}})$ for $\lambda \geq 0$ small.
\end{enumerate}
\end{corollary}
The branch $\mu= \mu(\lambda)$ for $\lambda \geq 0$ small can be found using the $BCH$ formula in $so(3)$.
The parameter space in the case of the resonant Hopf bifurcation is illustrated in Figure (\ref{figure:figura2}).

The results in Sections \ref{S:HopfBif}, \ref{S:Frequency} and \ref{S:Res_drift} show that a rotating wave on
a sphere generically undergoes a transition to quasi-periodic meandering (modulated rotating wave) at a
Hopf bifurcation. Furthermore, the primary frequency vector of the modulated rotating wave is determined by
both the critical Hopf eigenvalues and the frequency of the rotating wave undergoing the bifurcation.
In particular, resonances between the critical Hopf eigenvalue and the frequency of the rotating wave undergoing
the bifurcation lead to orthogonal drift. While we have concentrated on the case of Hopf bifurcation,
it should be mentioned that we would obtain similar results for \emph{ periodic forcing } of  a rotating wave: equation
(\ref{E:AuxFreq}) represents the dynamical equations of a periodically forced rotating wave.
In a spherical heart, for example, one could alter the dynamics of a rotating spiral wave simply by applying a
(weak) periodic forcing (e.g. pacing), thereby inducing meandering. By an appropriate choice of the
forcing frequency, one could then control the direction of the primary frequency vector (i.e. steer the spiral
tip to another location). Of course, this assumes perfect spherical symmetry. The case when geometrical imperfections
appear can also influence the dynamics of rotating waves and it will be presented somewhere else.

\section{Some Numerical Results}\label{S:Numerical}

Theorem \ref{thm:MRW} shows that, generically, by a supercritical Hopf bifurcation of a rotating wave, not only we get modulated rotating waves (that is a quasi-periodic tip motion), but also we get a quasi-periodic meandering tip motion. It is possible that this meandering motion is not of epicycle-type.
We illustrate this using the examples presented in Section \ref{S:Examples_Transition}.\\
Let $L_{x}=\left(\begin{array}{ccc}
              0 & 0 & 0\\
              0 & 0 & -1\\
              0 & 1 & 0\\
           \end{array} \right )$,
$L_{y}=\left ( \begin{array}{ccc}
                     0 & 0 & 1\\
                     0 & 0 & 0\\
                    -1 & 0 & 0\\
                  \end{array} \right )$,
$L_{z}=\left ( \begin{array}{ccc}
                 0 & -1 & 0\\
                 1 & 0 & 0\\
                 0 & 0 & 0\\
               \end{array} \right )$,\\
$R_{x}(\theta)= e^{L_{x}\theta}=\left ( \begin{array}{ccc}
                                           1 & 0 & 0\\
                                           0 & \cos \theta & -\sin \theta\\
                                           0 & \sin \theta & \cos \theta
                                     \end{array} \right )$,
$R_{y}(\theta)= e^{L_{y}\theta}=\left ( \begin{array}{ccc}
                                        \cos \theta & 0 & \sin \theta\\
                                                  0 & 1 & 0\\
                                       -\sin \theta & 0 & \cos \theta\\
                                       \end{array} \right )$,\\
$R_{z}(\theta)=e^{L_{z}\theta}=\left ( \begin{array}{ccc}
                                      \cos \theta & -\sin \theta & 0\\
                                      \sin \theta & \cos \theta & 0\\
                                                0 & 0 & 1\\
                                        \end{array} \right )$.\\

\begin{figure}[h]
\begin{center}
$\begin{array}{c@{\hspace{1in}}c}
\psfig{file=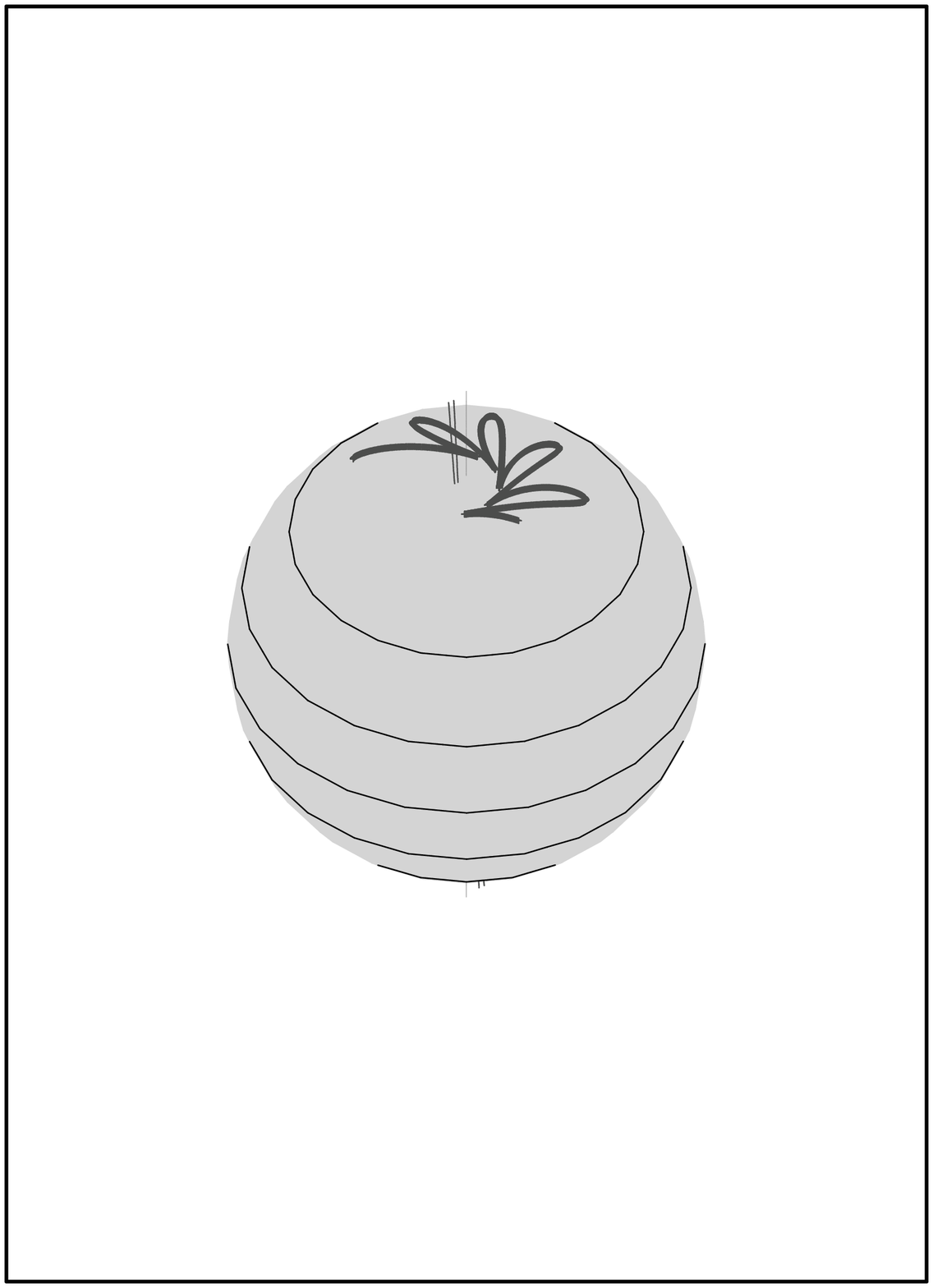,width=2in } &  \psfig{file=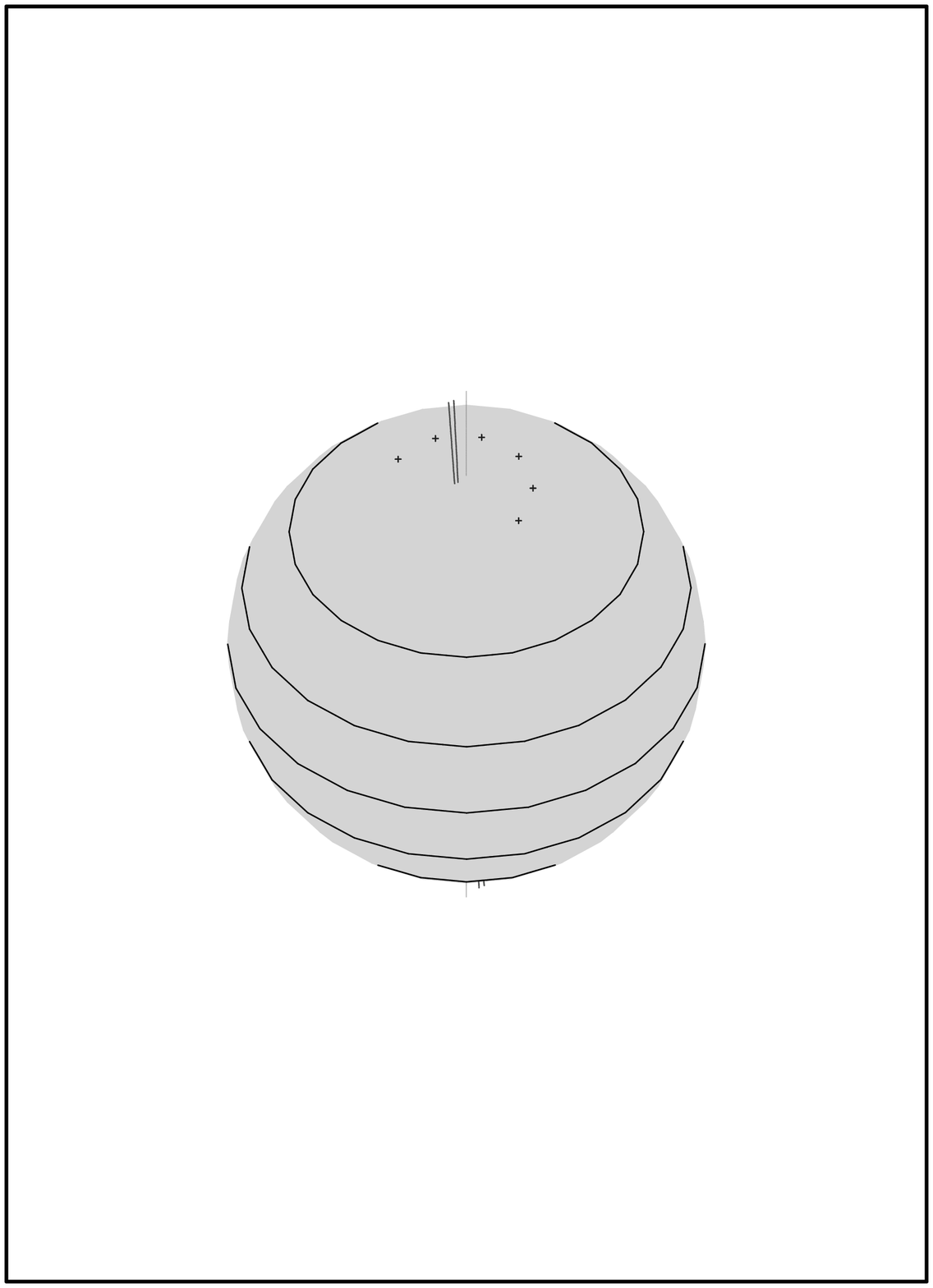,
width=2in }\\ [0.4cm] \mbox{$x_{tip}(\mathbf{\Phi}(t,u_{0.01},0.01))$} &
\mbox{$x_{tip}(\mathbf{\Phi}(i\frac{2\pi}{20.01},u_{0.01},0.01))$,$i=1,\ldots, 5$}
\end{array}$
\end{center}
\caption{Case 1, $\lambda=0.01$ }\label{figure:figtest-fig1}
\end{figure}

\begin{figure}[h]
\begin{center}
$\begin{array}{c} \psfig{file=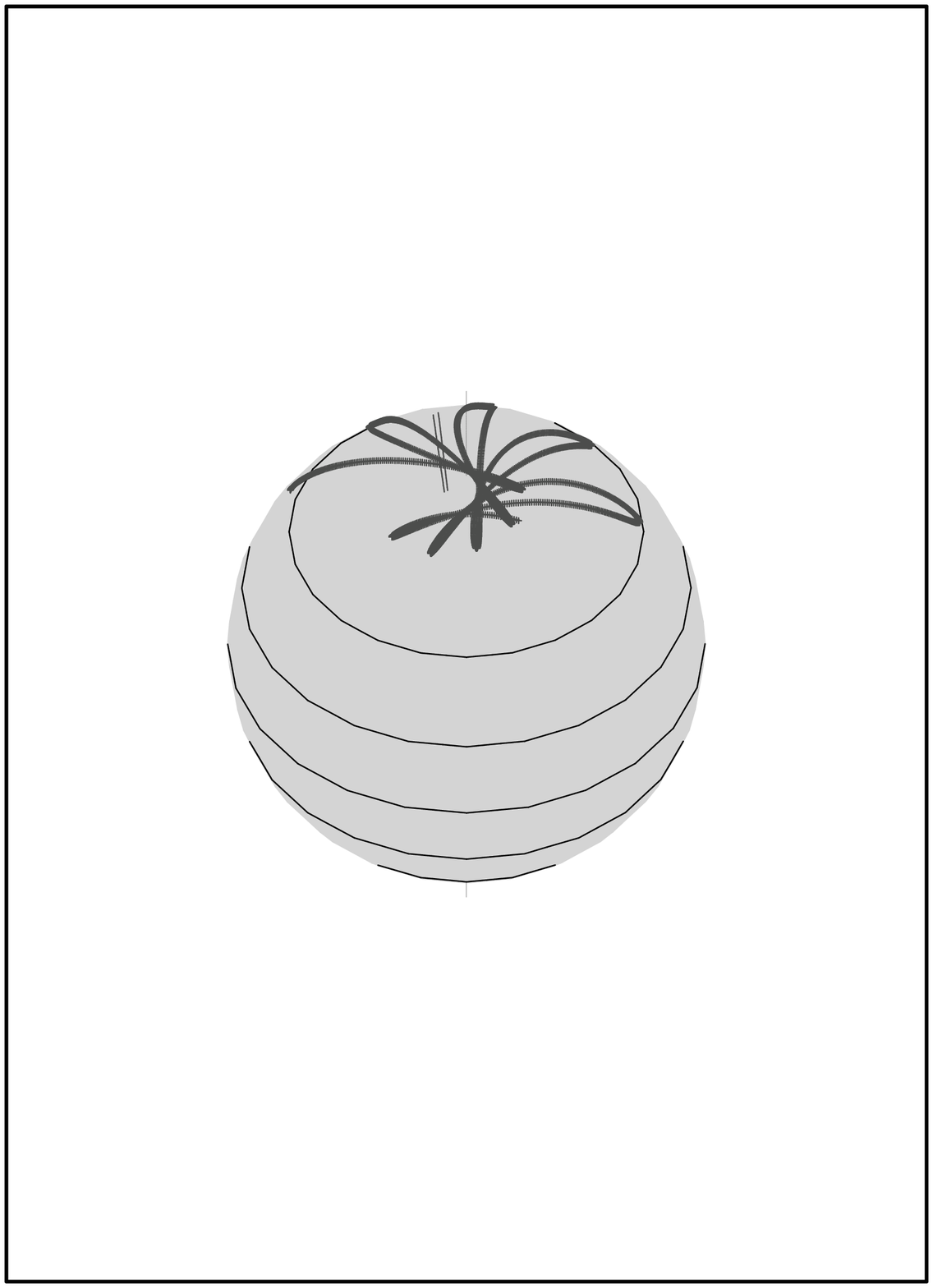,width=2in }   \\
[0.4cm] \mbox{$x_{tip}(\mathbf{\Phi}(t,u_{0.05},0.05))$}
\end{array}$
\end{center}
\caption{Case 1, $\lambda=0.05$ }\label{figure:figtest-fig2}
\end{figure}

\begin{figure}[h]
\begin{center}
$\begin{array}{c@{\hspace{1in}}c}
\psfig{file=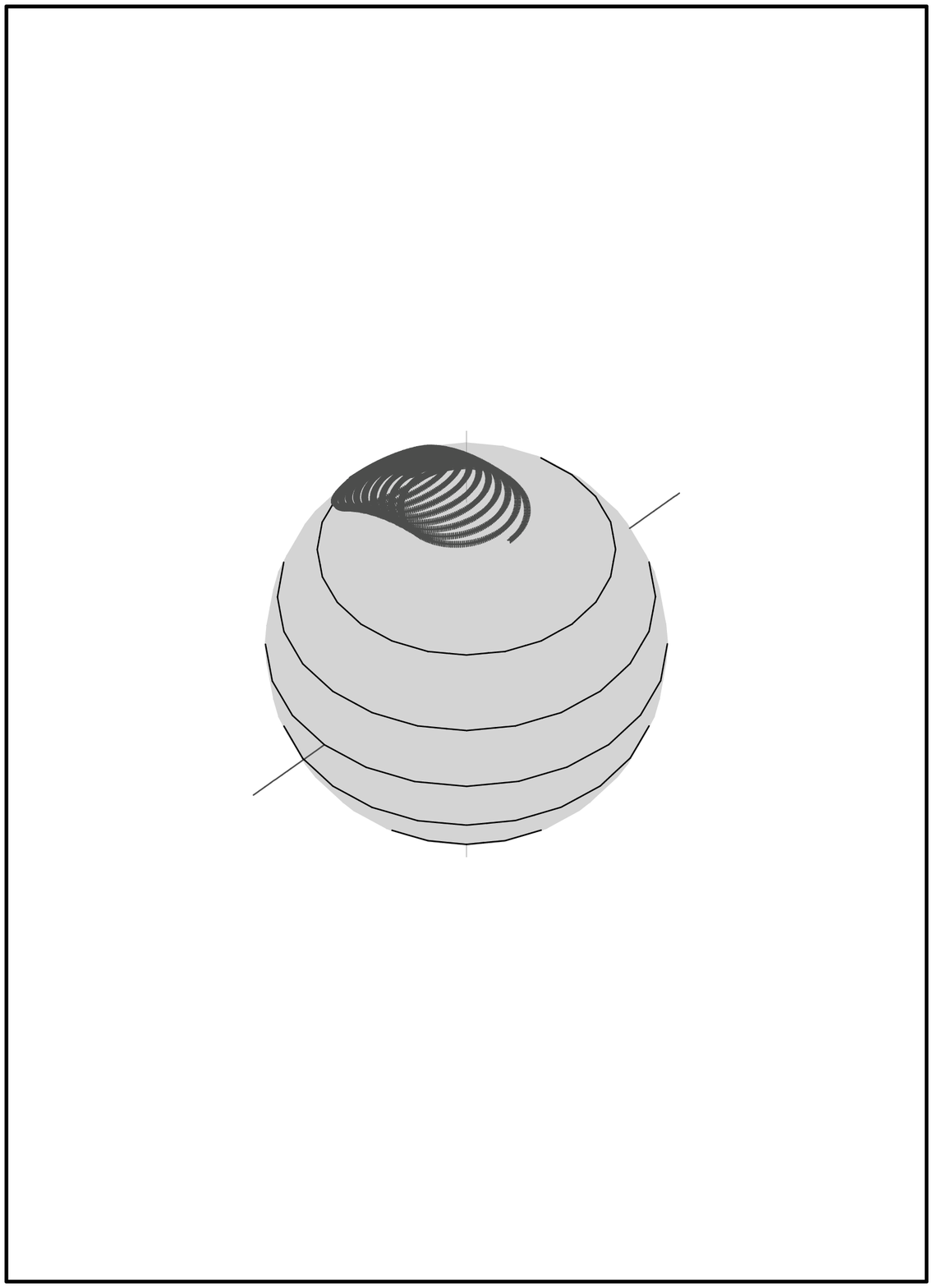,width=2in } &  \psfig{file=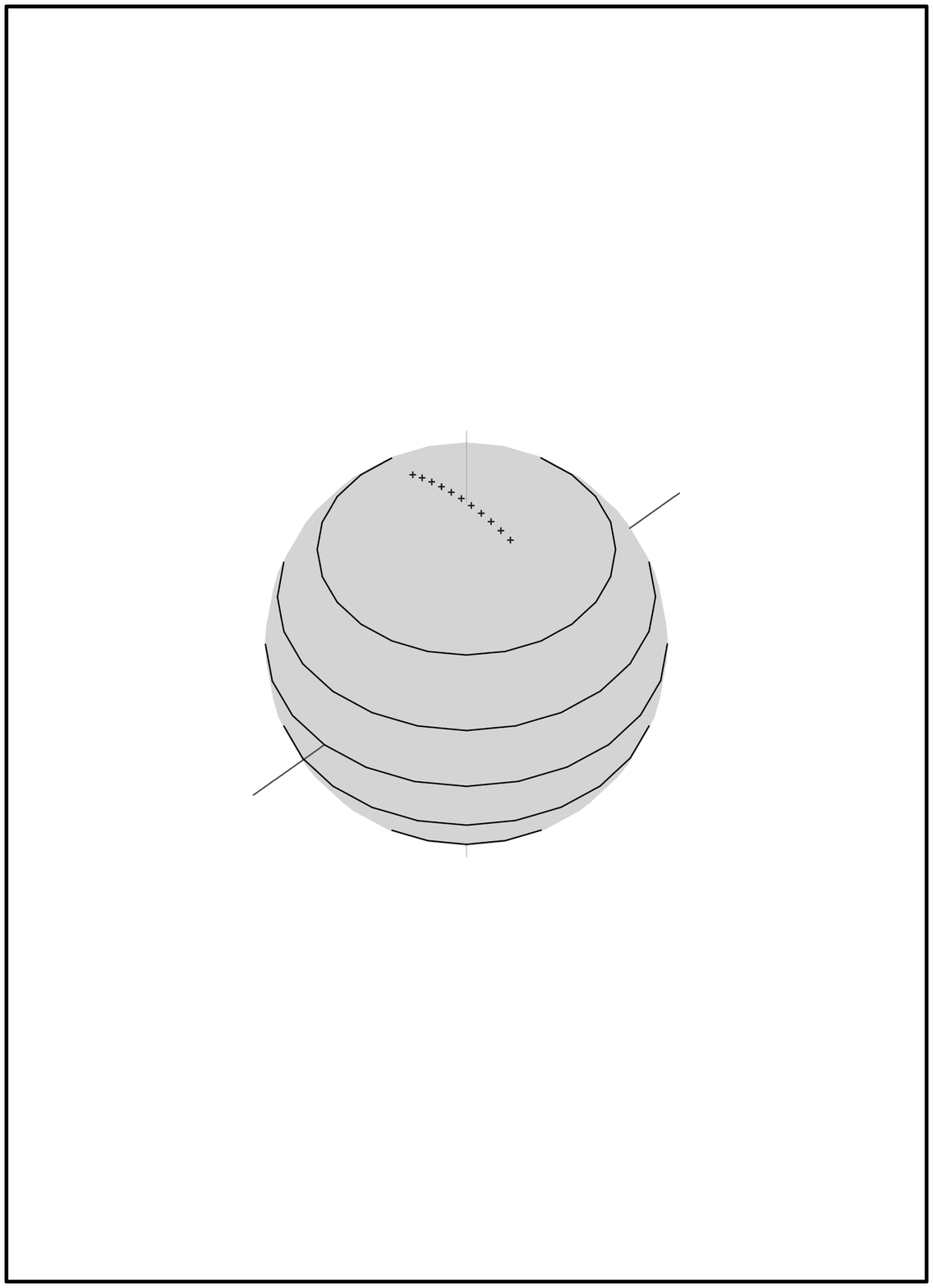,
width=2in }\\ [0.4cm] \mbox{$x_{tip}(\mathbf{\Phi}(t,u_{0.05},0.05))$} &
\mbox{$x_{tip}(\mathbf{\Phi}(i\frac{2\pi}{20.05},u_{0.05},0.05))$,$i=1,\ldots, 5$}
\end{array}$
\end{center}
\caption{Case 2, $\lambda=0.05$ }\label{figure:figtest-fig3}
\end{figure}

\begin{figure}[h]
\begin{center}
$\begin{array}{c} \psfig{file=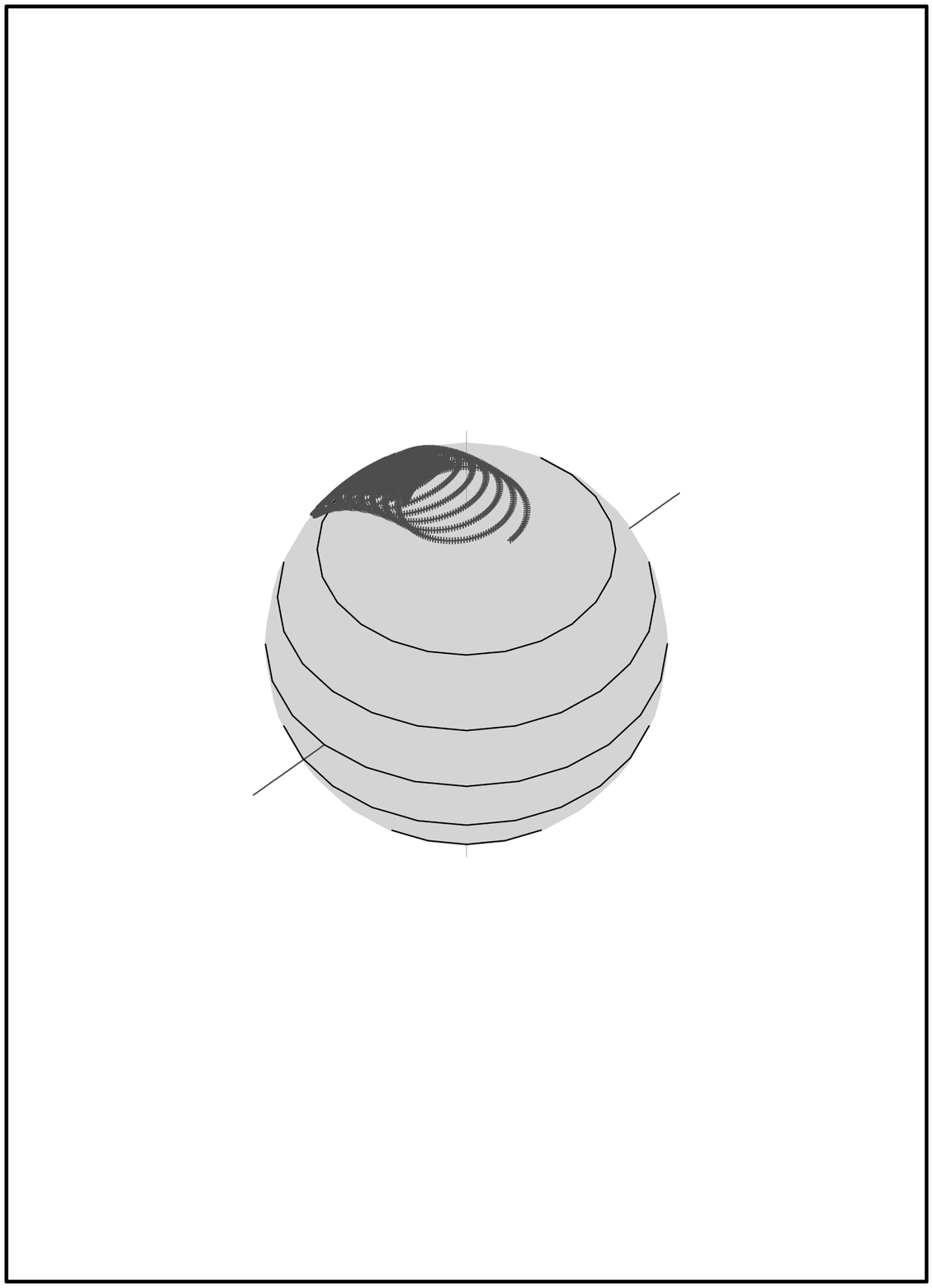,width=2in }   \\
[0.4cm] \mbox{$x_{tip}(\mathbf{\Phi}(t,u_{0.1},0.1))$}
\end{array}$
\end{center}
\caption{Case 2, $\lambda=0.1$ }\label{figure:figtest-fig4}
\end{figure}

\begin{figure}[h]
\begin{center}
$\begin{array}{c@{\hspace{1in}}c}
\psfig{file=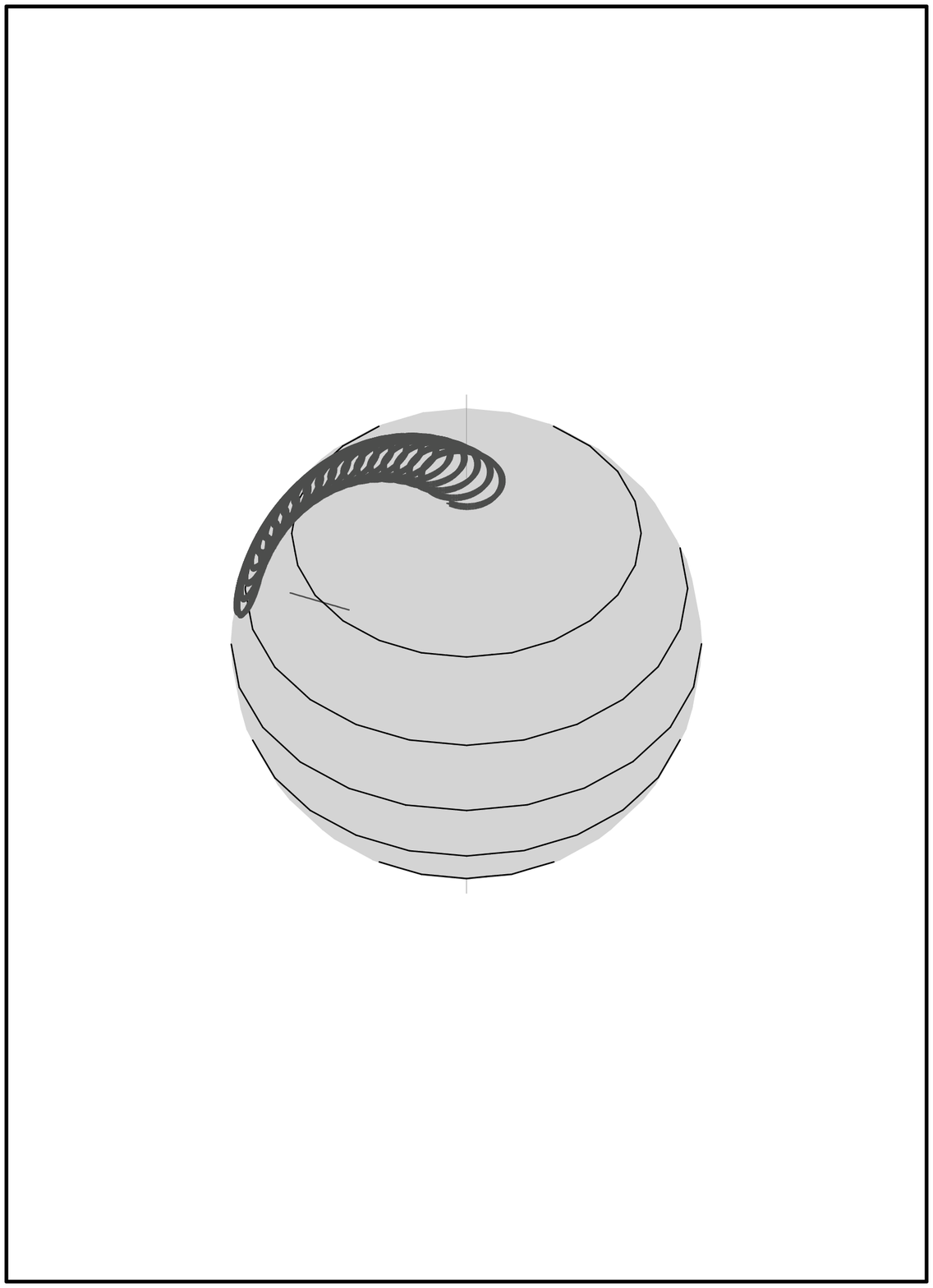,width=2in } &  \psfig{file=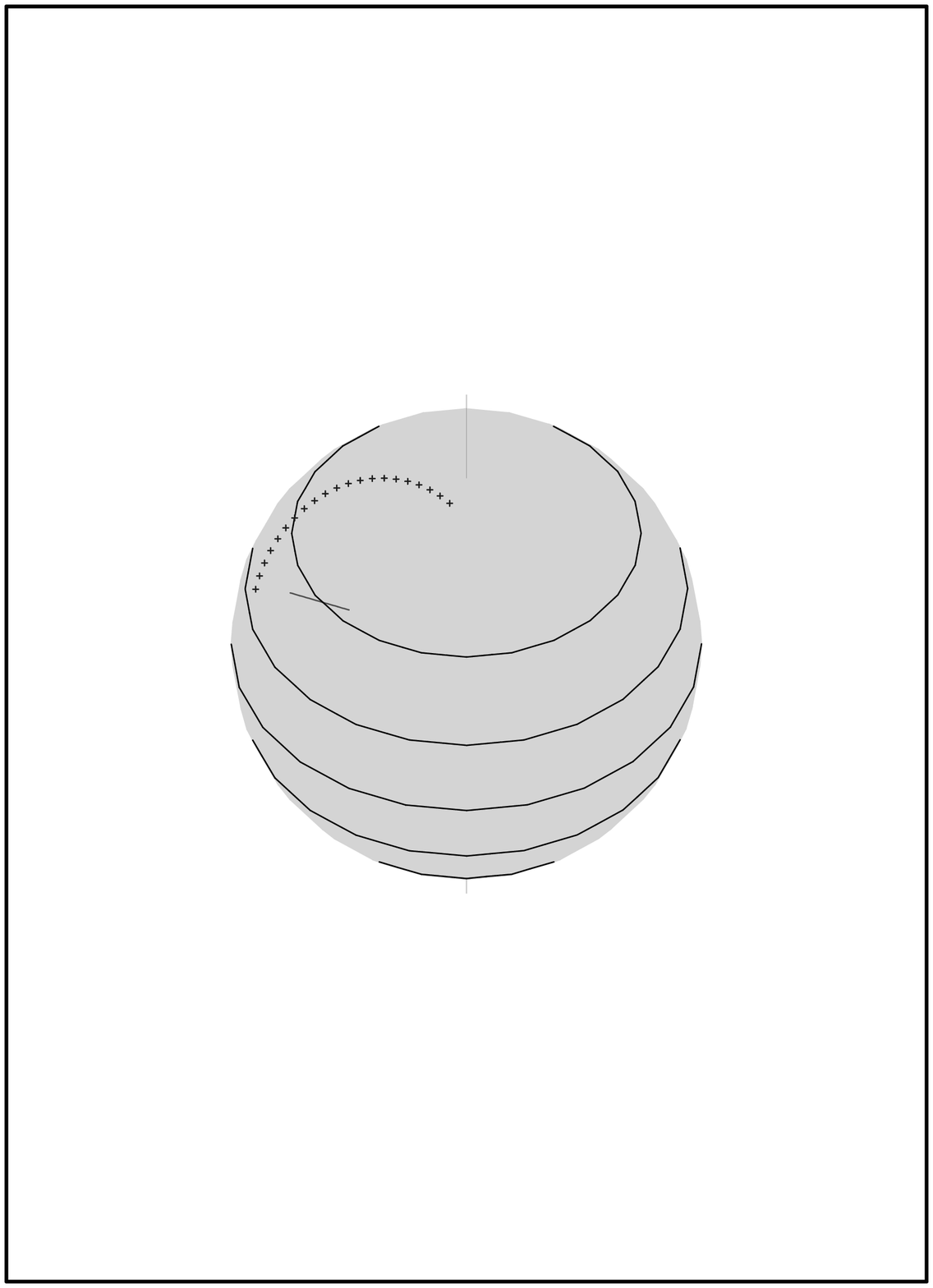,
width=2in }\\ [0.4cm] \mbox{$x_{tip}(\mathbf{\Phi}(t,u_{0.05},0.05))$} &
\mbox{$x_{tip}(\mathbf{\Phi}(i\frac{2\pi}{20.05},u_{0.05},0.015))$, $i=1,
\ldots, 5$}
\end{array}$
\end{center}
\caption{Case 3, $\lambda=0.05$ }\label{figure:figtest-fig5}
\end{figure}

\begin{figure}[h]
\begin{center}
$\begin{array}{c} \psfig{file=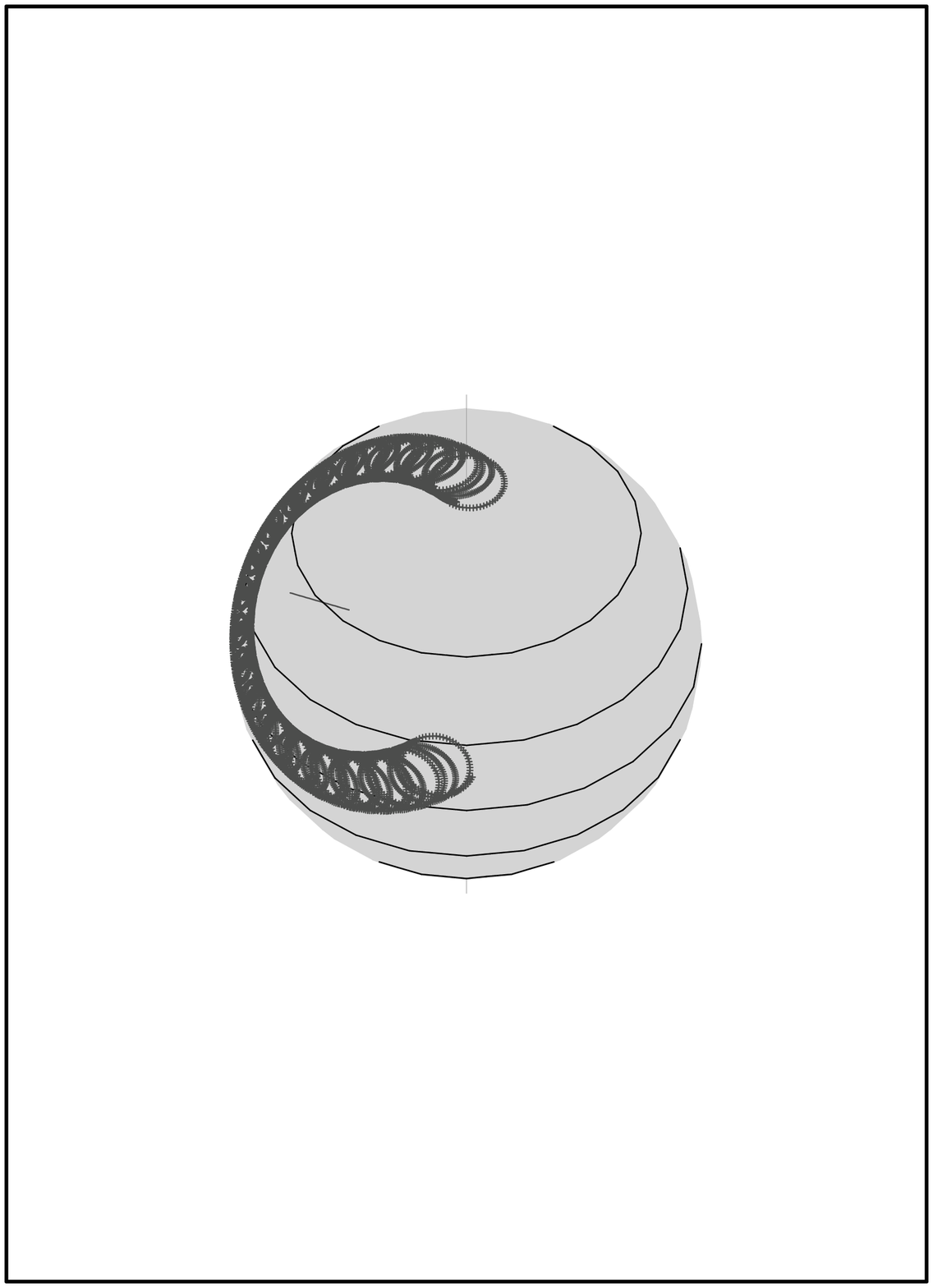,width=2in }   \\
[0.4cm] \mbox{$x_{tip}(\mathbf{\Phi}(t,u_{0.25},0.25))$}
\end{array}$
\end{center}
\caption{Case 3, $\lambda=0.25$ }\label{figure:figtest-fig6}
\end{figure}

We have seen that the study of the reaction-diffusion system (\ref{E:RDparam}) on $\mathbf{Y}^{\alpha}$ reduces to
the study of the finite-dimensional system (\ref{E:reduceddiffeqparam1}) on the center manifold $M^{cu}_{u_{0}}(\lambda)$ of the relative equilibrium $SO(3)u_{0}$, with $\mathbf{\Phi}(t,u_{0},0)=e^{X_{0}t}u_{0}$.
Let us denote $X_{G}(q,\lambda)= F^{x}(q,\lambda)L_{x}+F^{y}(q,\lambda)L_{y}+F^{z}(q,\lambda)L_{z}$ for any $q \in V_{*}$ and $\left | \lambda \right | $ small.\\
If we parameterize $SO(3)$ by Euler angles, that is $A=R_{z}(\psi)R_{x}(\theta)R_{z}(\phi)$, where
$\phi \in [0,2\pi)$, $\theta \in [0,\pi]$, $\psi \in [0,2\pi)$, the finite dimensional system (\ref{E:reduceddiffeqparam1}) becomes
\begin{equation}\label{E:Numerical1}
\begin{array}{lll}
\dot{\phi} &=& F^{z}(q,\lambda)-\cot \theta \left [F^{y}(q,\lambda) \cos \phi +F^{x}(q,\lambda) \sin \phi \right ],\\
\dot{\theta} &=& -F^{y}(q,\lambda)\sin \phi+ F^{x}(q, \lambda)\cos \phi,\\
\dot{\psi} &=& \frac{1}{\sin \theta}\left [F^{y}(q,\lambda)\cos \phi + F^{x}(q,\lambda)\sin \phi \right ],\\
\dot{q} &=& X_{N}(q,\lambda).
\end{array}
\end{equation}
We consider the following initial value problem associated with the system (\ref{E:Numerical1}):
\begin{equation}\label{E:Numerical11}
\begin{array}{lll}
\dot{\phi} &=& F^{z}(q,\lambda)-\cot \theta [F^{y}(q,\lambda) \cos \phi + F^{x}(q,\lambda) \sin \phi],\\
\dot{\theta} &=& -F^{y}(q,\lambda)\sin \phi+ F^{x}(q, \lambda)\cos \phi,\\
\dot{\psi} &=& \frac{1}{\sin \theta}[F^{y}(q,\lambda)\cos \phi +F^{x}(q,\lambda)\sin \phi],\\
\phi(0) &=& 0,\\
\theta(0) & \neq & 0,\\
\psi(0) &=& 0,\\
\dot{q} &=& X_{N}(q,\lambda).
\end{array}
\end{equation}
We choose $\theta(0) \neq 0$ near $0$ in (\ref{E:Numerical11}).\\
We consider $V_{*} \simeq \mathbb{C}$ and a supercritical Hopf bifurcation takes place in $\dot{q}=X_{N}(q,\lambda)$ at $q=0$ for $\lambda=0$ with eigenvalues $\pm i \omega_{bif}$. Let $q(t,\lambda)$ be the periodic solution of the second differential equation in (\ref{E:reduceddiffeqparam}), that appears by a supercritical Hopf bifurcation for $\lambda > 0$ small. Let $T(\lambda)=\frac{2\pi}{\left | \omega_{\lambda}\right | }$ be its period, where $\omega_{\lambda}=\omega_{bif}+O(\lambda)$ for $\lambda \geq 0$ small. Let $q(t,0)=0$ for any $t \in \mathbb{R}$. If we substitute $q(t,\lambda)$ in $X_{G}(q,\lambda)$ and in the first three equations of the system (\ref{E:Numerical11}), we get $X^{G}(t,\lambda)=X_{G}(q(t,\lambda),\lambda)= F^{x}(q(t,\lambda), \lambda)L_{x}+ F^{y}(q(t,\lambda),\lambda)L_{y}+F^{z}(q(t,\lambda),\lambda)L_{z}=\left | X_{0} \right | X_{0}^{1}+\sqrt{\lambda}H(t,\lambda)$ for $\lambda \geq 0$ small and $t \in \mathbb{R}$.
If we write $F^{xx}(t,\lambda)= F^{x}(q(t,\lambda),\lambda)$, $F^{yy}(t,\lambda)= F^{y}(q(t,\lambda),\lambda)$ and $F^{zz}(t,\lambda)= F^{z}(q(t,\lambda),\lambda)$, we get that $X^{G}(t,\lambda)= F^{xx}(t,\lambda)L_{x}+ F^{yy}(t,\lambda)L_{y}+F^{zz}(t,\lambda)L_{z}$ is $\frac{2\pi}{\left | \omega_{\lambda} \right |}$-periodic.\\
The following initial value problem is obtained
\begin{equation}\label{E:Numerical2}
\begin{array}{lll}
\dot{\phi} &=& F^{zz}(t,\lambda)-\cot \theta [F^{yy}(t,\lambda)\cos \phi + F^{xx}(t,\lambda)\sin \phi],\\
\dot{\theta} &=& -F^{yy}(t,\lambda)\sin \phi+ F^{xx}(t, \lambda)\cos \phi,\\
\dot{\psi} &=& \frac{1}{\sin \theta}[F^{yx}(t,\lambda)\cos \phi + F^{xx}(t,\lambda)\sin \phi],\\
\phi(0) &=& 0,\\
\theta(0) & \neq & 0,\\
\psi(0) &=& 0.
\end{array}
\end{equation}
We choose $\theta(0) \neq 0$ near $0$ in (\ref{E:Numerical2}).\\
We use a Maple program that integrates numerically the system (\ref{E:Numerical2}) and finds numerically a
primary frequency vector (\cite{Se}, if its norm is not $0$ or $\pi$) of the modulated rotating wave that appears by a
supercritical Hopf bifurcation, as discussed in Section \ref{S:HopfBif}. Then, for a choice of the point
$x_{0} \in r\mathbf{S}^{2}$ near $\frac{r}{\left | X_{0} \right | }\overrightarrow X_{0}$, we represent the
tip motion
\begin{equation}\label{E:motion1}
x_{tip}(\mathbf{\Phi}(t,u_{\lambda},\lambda))= A(0,\lambda)^{-1}A(t,\lambda)x_{0} \mbox{ on } r\mathbf{S}^{2}
\end{equation}
for $\lambda \geq 0$ small. Recall that $\mathbf{\Phi}(t,u_{\lambda},\lambda)= A(0,\lambda)^{-1}A(t,\lambda)\mathbf{\Psi}(q(t,\lambda))$.\\
We consider the following cases:
\begin{enumerate}
\item [Case 1.] $X(t,\lambda)$ is given in Example (\ref{ex:exHoph1}) from Section \ref{S:Examples_Transition};
$X_{0}^{1}= L_{z}, X_{1}=L_{x}, X_{2}=L_{y}$; $\omega_{bif}=20$, $\left | X_{0} \right |=2$, $g(t,\lambda)=
\sin((\omega_{bif}+\lambda)t)$, $r=3$, $\theta(0)=0.01$,
$x_{0}=\left (\begin{array}{c}
         0\\ 0.92\\2.85\\
              \end{array} \right )$.
We get Figures \ref{figure:figtest-fig1} to \ref{figure:figtest-fig2} for $\lambda=0.01$, $0.05$.\\
Since $\left | X_{0} \right | \neq k \omega_{bif}$ for any $k \in \mathbb{Z}$, we have the nonresonant case.
On Figures \ref{figure:figtest-fig1} and  \ref{figure:figtest-fig2}, it is visualized the tip motion given by (\ref{E:motion1}) of the modulated
rotating waves $\mathbf{\Phi}(t,u_{\lambda},\lambda)$ obtained by a supercritical Hopf bifurcation that takes place in
$\dot{q}=X_{N}(q,\lambda)$ at $q=0$ for $\lambda=0$. On Figure \ref{figure:figtest-fig1} we consider $\lambda=0.01$.
On Figure \ref{figure:figtest-fig2} we consider $\lambda=0.05$.\\
On the second of Figure \ref{figure:figtest-fig1} we plot the points $x_{tip}(\mathbf{\Phi}(i\frac{2\pi}{\omega_{\lambda}},u_{\lambda},\lambda))$
for some $i=0$,~$1$,~$2$,\ldots \, and we can see that they are points of a circle on the sphere $r\mathbf{S}^{2}$ with
the center on the line having the direction of the primary frequency vector of $\mathbf{\Phi}(t, u_{\lambda},\lambda)$,
where $\lambda=0.01$. \\
The grey line is the line containing the frequency of the rotating wave undergoing the Hopf bifurcation,
$\overrightarrow{X_{0}}$. The black lines are the lines corresponding to the primary frequency vector associated to the
modulated rotating wave. One of them is computed numerically and the other is the exact one. We can see that they are
very close.

\item [Case 2.] $X(t,\lambda)$ is given in Example (\ref{ex:exHoph2}) from Section \ref{S:Examples_Transition};
$X_{0}^{1}= L_{z},X_{1}=L_{x},X_{2}=L_{y}$; $\omega_{bif}= \left | X_{0} \right |=20$,
$g(t,\lambda)= \sin((\omega_{bif}+\lambda)t)$, $r=3$, $\theta(0)=0.02$,
$x_{0}=\left (\begin{array}{c}
         0\\ 0.92\\2.85\\
       \end{array} \right )$.\\
We get Figures \ref{figure:figtest-fig3} to \ref{figure:figtest-fig4} for $\lambda=0.05$, $0.1$.\\
Since $\left | X_{0} \right | = k \omega_{bif}$ for $k=1 \in \mathbb{Z}$, we have the resonant case. On Figures
\ref{figure:figtest-fig3} and \ref{figure:figtest-fig4}, it is visualized the tip motion given by (\ref{E:motion1}) of
the modulated rotating waves $\mathbf{\Phi}(t,u_{\lambda},\lambda)$ obtained by a supercritical Hopf bifurcation that
takes place in $\dot{q}=X_{N}(q,\lambda)$ at $q=0$ for $\lambda=0$. On Figure \ref{figure:figtest-fig3} we consider
$\lambda=0.05$. On Figure \ref{figure:figtest-fig4} we consider $\lambda=0.1$. We can see that the primary frequency
vectors associated to the modulated rotating waves are orthogonal to $\overrightarrow{X_{0}}$. We call this
phenomenon resonant drift.\\
On the second of Figure \ref{figure:figtest-fig3} we plot the points $x_{tip}(\mathbf{\Phi}(i\frac{2\pi}{\omega_{\lambda}},
u_{\lambda},\lambda))$ for some $i=0$,~$1$,~$2$,\ldots \, and we can see that they are points of a circle
 on the sphere $r\mathbf{S}^{2}$ with the center on the line having the direction of the primary frequency
 vector of $\mathbf{\Phi}(t, u_{\lambda},\lambda)$, where $\lambda=0.05$.\\
The grey line is the line containing the frequency of the rotating wave undergoing the Hopf bifurcation,
$\overrightarrow{X_{0}}$. The black line is the line corresponding to the primary frequency
vector associated to the modulated rotating wave.

\item [Case 3.] $X(t,\lambda)$ is given in Example (\ref{ex:exHoph3}) from Section \ref{S:Examples_Transition};
$X_{0}^{1}= L_{z}, X_{1}=L_{x}, X_{2}=L_{y}$; $\omega_{bif}=\left | X_{0} \right | =20$, $g(t,\lambda)=
\sin((\omega_{bif}+\lambda)t)$, $r=3$, $\theta(0)=0.5$,
$x_{0}=\left (\begin{array}{c}
         0.44\\ 0.14\\2.96\\
       \end{array} \right ).$\\
We get Figures \ref{figure:figtest-fig5} to \ref{figure:figtest-fig6} for $\lambda=0.05$, $0.25$.\\
Since $\left | X_{0} \right | = k \omega_{bif}$ for $k=1 \in \mathbb{Z}$, we have the resonant case. On Figures
\ref{figure:figtest-fig5} and \ref{figure:figtest-fig6}, it is visualized the tip motion given by (\ref{E:motion1}) of
the modulated rotating waves $\mathbf{\Phi}(t,u_{\lambda},\lambda)$ obtained by a supercritical Hopf bifurcation that
takes place in $\dot{q}=X_{N}(q,\lambda)$ at $q=0$ for $\lambda=0$. On Figure \ref{figure:figtest-fig5} we consider
$\lambda=0.05$. On Figure \ref{figure:figtest-fig6} we consider $\lambda=0.25$. We can see that the primary frequency
vectors associated to the modulated rotating waves are not orthogonal to $\overrightarrow{X_{0}}$. This can
happen in the resonant case if we consider only one parameter $\lambda$.\\
On the second of Figure \ref{figure:figtest-fig5} we plot the points
$x_{tip}(\mathbf{\Phi}(i\frac{2\pi}{\omega_{\lambda}},u_{\lambda},\lambda))$for some $i=0$,~$1$,~$2$,\ldots \, and
we can see that they are points of  a circle on the sphere $r\mathbf{S}^{2}$ with the center on the line having
the direction of the primary frequency vector of $\mathbf{\Phi}(t, u_{\lambda},\lambda)$, where $\lambda=0.05$.\\
The grey line is the line containing the frequency of the rotating wave undergoing the Hopf bifurcation,
$\overrightarrow{X_{0}}$. The black line is the line corresponding to the primary frequency vector associated to
the modulated rotating wave.
\end{enumerate}
\newpage
The nonuniformly rigid rotation shown in Example \ref{ex:exHoph4} is not generically. Therefore, we do not present it here.

\section{Proofs of Theorems}\label{S:Proof_Transition}

\begin{proof}[\textbf{Proof of Lemma \ref{lem:Hopf_group}}]
Let us omit the parameter $\lambda$. Since $SO(3)$ is a compact manifold and $X^{G}$ is a sufficiently $CS$ function, the initial value problem
(\ref{E:auxHopf}) has a unique sufficiently $CS$ solution that is globally defined. Let $A^{*}(t)$ be this solution. Since the exponential map $exp \colon so(3) \rightarrow SO(3)$ is surjective, there exists a matrix $X \in so(3)$ such
that $A^{*}(T)= e^{XT}.$\\
Let us make the following change of variable $B= e^{-Xt}A$. \\
Then the first equation of (\ref{E:auxHopf}) becomes:
\begin{equation}
\dot{B}= -e^{-Xt}XA(t)+e^{-Xt}\dot{A}(t),
\end{equation}
or since $\dot{A}= AX^{G}(t)$, we get
\begin{equation}
\dot{B}= -e^{-Xt}XA(t)+e^{-Xt}AX^{G}(t),
\end{equation}
or since $A=e^{Xt}B$ it follows that
\begin{equation}
\dot{B}= -XB(t)+B(t)X^{G}(t).
\end{equation}
The initial condition $A(0)= I_{3}$ becomes $B(0)= I_{3}$.\\
Let us consider now the initial value problem
\begin{equation}\label{E:auxHopf_B}
\begin{array}{lll}
\dot{B} &=& -XB+BX^{G}(t),\\
 B(0)   &=& I_{3}.
\end{array}
\end{equation}
The initial value problem (\ref{E:auxHopf_B}) has a unique global solution since $SO(3)$ is a compact manifold.\\
Since $A^{*}(t)$ is the solution of the initial value problem (\ref{E:auxHopf}), it follows that $B^{*}(t)= e^{-Xt}A^{*}(t)$ is a solution of (\ref{E:auxHopf_B}) and  $B^{*}(T)= I_{3}$.\\
Let us define \[C^{*}(t)= B^{*}(t+T).\] We will check that $C^{*}$ is another solution of (\ref{E:auxHopf_B}). We have $C^{*}(0)= B^{*}(T)= I_{3}$ and for any $t \in [0,\infty)$,
\begin{equation}
\begin{split}
\frac{d}{d t}\left( C^{*}(t) \right )&= \frac{d}{d t}\left ( B^{*}(t+T) \right )= \frac{d B^{*}}{d t}(t+T)= -XB^{*}(t+T)+B^{*}(t+T)X^{G}(t+T)\\
&= -XB^{*}(t+T)+B^{*}(t+T)X^{G}(t)= -XC^{*}(t)+C^{*}(t)X^{G}(t),
\end{split}
\end{equation}
where we have used the fact that $X^{G}$ is a $T$-periodic function.\\
Therefore, $B^{*}(t)= C^{*}(t)$ for $t \in [0,\infty)$, that is, $B^{*}$ is $T$-periodic.\\
So, $A^{*}(t)= e^{Xt}B^{*}(t)$ for $t \in [0,\infty)$.\\
Taking into account the parameter $\lambda > 0$, since $X^{G}(t,\lambda)$ is sufficiently $CS$ for $t \in [0,\infty)$ and $\lambda > 0$, we have that $A(t,\lambda)$ is sufficiently $CS$ for $t \in [0, \infty)$ and $\lambda > 0$ small. If we work with $\lambda \geq 0$, then $A(t, \lambda)$ is sufficiently $CS$ for $t \in [0, \infty)$ and $\lambda \geq 0$ small.
\end{proof}

\begin{proof}[\textbf{Proof of Theorem \ref{thm:Hopf_RW}}]

We apply Theorems \ref{thm:CMR_REparam}. Then, it is sufficient to study the differential equations (\ref{E:reduceddiffeqparam1}).\\
The supercritical Hopf bifurcation in $V^{*}$ for the second differential equation in (\ref{E:reduceddiffeqparam1}) at $q= 0$ for $\lambda= 0$ implies that there exists a unique sufficiently $CS$ branch of periodic solutions $q(t, \lambda)$ for $\lambda >0$ near $q= 0$.\\
Recall that we have denoted by $T(\lambda)= \frac{2\pi}{\left | \omega_{\lambda}\right | }$ the period of $q(t, \lambda)$, where $\omega_{\lambda}= \omega_{bif}+O(\lambda)$ is a sufficiently smooth function for $\lambda \geq 0$ small.\\
If we substitute $q(t, \lambda)$ in the first of the differential equations (\ref{E:reduceddiffeqparam}), we get the differential equation
\begin{equation}\label{E:eq1}
\dot{A} = AX^{G}(t, \lambda),
\end{equation}
where we recall that the function $X^{G}(t,\lambda)$ is sufficiently $CS$ and $\frac{2\pi}{\left | \omega_{\lambda}\right |}$-periodic for $t \in [0,\infty)$ and $\lambda \geq 0$ small (and it is defined by (\ref{E:InitialData})).\\
We consider the following initial value problem on $SO(3)$:
\begin{equation}\label{E:PeriodicEq}
\begin{array}{lll}
\dot{A} &=& AX^{G}(t, \lambda),\\
 A(0)   &=& I_{3},\\
\end{array}
\end{equation}
because any solution of the differential equation (\ref{E:eq1}) with $A(0)= A_{0}$, where $A_{0} \in SO(3)$ is given by $A_{0}A^{*}(t,\lambda)$ and $A^{*}(t,\lambda)$ is the solution of the initial value problem (\ref{E:PeriodicEq}).\\
For $\lambda > 0$ small, by applying Lemma \ref{lem:Hopf_group}, we get
\[A(t,\lambda)= e^{X(\lambda)t}B(t, \lambda) \mbox{ for } t \in [0,\infty) \mbox{ and } \lambda \geq 0 \mbox{ small },\] where $B(t, \lambda)$ has period $T(\lambda)$ and $B(0, \lambda)= I_{3}$.\\
Consequently, by using Theorems \ref{thm:CMR_REparam}, there exists a solution
$\mathbf{\Phi}(t,u_{\lambda},\lambda)$ of the reaction-diffusion system (\ref{E:RDparam}) of the form:
\[\mathbf{\Phi}(t, u_{\lambda}, \lambda)= A(t, \lambda)\Psi(q(t, \lambda))=
e^{X(\lambda)t}B(t, \lambda)\Psi(q(t, \lambda)), \mbox { where } t \in [0,\infty) \mbox{ and } \lambda > 0 \,
\mbox{ small },\] with $\Psi(q(0, \lambda))= u_{\lambda}$ and $\Psi$ appears in the statement of
Theorem \ref{thm:CMR_REparam} (see  Appendix \ref{S:app2_CMR}).
If we write \[q_{1}(t, \lambda)= B(t, \lambda)\Psi(q(t, \lambda)),\] then the function $q_{1}(t, \lambda)$ is $T(\lambda)$-periodic in $t$, since $B(t,\lambda)$ and $q(t,\lambda)$ are $T(\lambda)$-periodic in $t$. Also, we get
\[\mathbf{\Phi}(t,u_{\lambda},\lambda)=e^{X(\lambda)t}q_{1}(t, \lambda) \mbox{ and } \mathbf{\Phi}(T(\lambda), u_{\lambda}, \lambda)= e^{X(\lambda)T(\lambda)}u_{\lambda}.\] We check that, at least generically, $\mathbf{\Phi}(t,u_{\lambda},\lambda) \notin SO(3)u_{\lambda}$ for any $t \in (0,T(\lambda))$ and $\lambda >0$ small. \\
We prove that $\overrightarrow{X(\lambda)}$ is a primary frequency vector of $\mathbf{\Phi}(t,u_{\lambda},\lambda)$.\\
Suppose that there exists a $T_{1}(\lambda) \in (0,T(\lambda))$ such that $\mathbf{\Phi}(T_{1}(\lambda),u_{\lambda},\lambda) \in SO(3)u_{\lambda}$ and for $t \in (0,T_{1}(\lambda))$ we have $\mathbf{\Phi}(t,u_{\lambda},\lambda) \notin SO(3)u_{\lambda}$.
Let $T(\lambda)=lT_{1}(\lambda)+ r$, where $r \in [0,T_{1}(\lambda))$ and  $l \in \mathbb{Z}$.
Let $\mathbf{\Phi}(T_{1}(\lambda),u_{\lambda},\lambda)= e^{Y(\lambda)T_{1}(\lambda)}u_{\lambda}$,
where $Y(\lambda) \in so(3)$.
Using the definition and $SO(3)$-equivariance of the semiflow $\mathbf{\Phi}$, we get
$\mathbf{\Phi}(T(\lambda),u_{\lambda},\lambda)=\mathbf{\Phi}(r,\mathbf{\Phi}(lT_{1}(\lambda), u_{\lambda},\lambda))=
e^{lY(\lambda)T_{1}(\lambda)}\mathbf{\Phi}(r, u_{\lambda},\lambda)$.\\
Therefore, $\mathbf{\Phi}(r, u_{\lambda},\lambda)= e^{-lY(\lambda)T_{1}(\lambda)}e^{X(\lambda)T(\lambda)}u_{\lambda}
\in SO(3)u_{\lambda}$ and $r \in [0,T(\lambda))$. It follows that $r=0$ and so $T(\lambda)=lT_{1}(\lambda)$. Thus, $\mathbf{\Phi}(T(\lambda), u_{\lambda}, \lambda)= e^{X(\lambda)T(\lambda)}u_{\lambda}=e^{Y(\lambda)T(\lambda)}u_{\lambda}$ and $e^{X(\lambda)T(\lambda)}=e^{Y(\lambda)T(\lambda)}$ because $\Sigma_{u_{\lambda}}=I_{3}$.\\
Also, we get $\mathbf{\Phi}(t,u_{\lambda},\lambda)=e^{X(\lambda)t}q_{1}(t, \lambda)=e^{Y(\lambda)t}e^{-Y(\lambda)t}
e^{X(\lambda)t} B(t, \lambda)\Psi(q(t, \lambda))$.\\
Let $B_{1}(t, \lambda)=e^{-Y(\lambda)t}e^{X(\lambda)t}B(t,\lambda)$. $B_{1}(t, \lambda)$ is $T(\lambda)$-periodic, because $B(t,\lambda)$ is $T(\lambda)$-periodic and $e^{X(\lambda)T(\lambda)}=e^{Y(\lambda)T(\lambda)}$. Thus,
\[\mathbf{\Phi}(t,u_{\lambda},\lambda)=e^{Y(\lambda)t}B_{1}(t, \lambda)\Psi(q(t,\lambda)).\]
Then
$\mathbf{\Phi}(T_{1}(\lambda),u_{\lambda},\lambda)=e^{Y(\lambda)T_{1}(\lambda)}u_{\lambda}=
e^{Y(\lambda)T_{1}(\lambda)}B_{1}(T_{1}(\lambda),\lambda)\Psi(q(T_{1}(\lambda),\lambda))$,
and it follows that $\Psi(q(T_{1}(\lambda),\lambda))=(B(T_{1}(\lambda),\lambda))^{-1}u_{\lambda}$.\\
We have that $u_{\lambda}$ is close to $u_{0}$ and $SO(3)u_{0} \cap \Psi(V_{*})=\{u_{0}\}$, therefore, generically,
$SO(3)u_{\lambda} \cap \Psi(V^{*})=\{u_{\lambda}\}$. Since $\Psi(q(T_{1}(\lambda),\lambda))=(B(T_{1}(\lambda),\lambda))^{-1}u_{\lambda}
\in \Psi(V_{*}) \cap SO(3)u_{\lambda}$, it follows that $B(T_{1}(\lambda),\lambda)u_{\lambda}= u_{\lambda}$ and, since
$\Sigma_{u_{\lambda}}=I_{3}$, $B_{1}(T_{1}(\lambda),\lambda)=I_{3}$ and $\Psi(q(T_{1}(\lambda),\lambda))= u_{\lambda}=
\Psi(q(0,\lambda))$. Because $\Psi$ is a local diffeomorphism (see \cite{FSSW, SSW}) and
$q(t,\lambda)=\sqrt{\lambda}r(t,\lambda)$, it results that $q(T_{1}(\lambda), \lambda)= q(0,\lambda)$. We know that $q(t,\lambda)$ is a solution of the differential equation $\dot{q}= X_{N}(q,\lambda)$. Therefore, it follows that
$q(t,\lambda)$ is $T_{1}(\lambda)$-periodic. But it has the period $T(\lambda)$. That is, $T(\lambda)=T_{1}(\lambda)$ contradicting the fact that $T_{1}(\lambda) \in (0,T(\lambda))$.\\
If $\left | X(\lambda)\right | T(\lambda)= 2k\pi$ for some $k \in \mathbb{Z}$, then $\mathbf{\Phi}(t, u_{\lambda}, \lambda)$ is a periodic solution with the period $T(\lambda)$.\\
Otherwise, $\mathbf{\Phi}(t, u_{\lambda}, \lambda)$ is a modulated rotating wave with a primary frequency vector $\overrightarrow {X(\lambda)}$.\\
By Proposition \ref{prop:Stability_RW}, the orbital stability of the modulated rotating wave or of the periodic solution obtained above is the same as the stability of the periodic solution $q(t, \lambda)$.\\
Therefore, since the periodic solution $q(t, \lambda)$ for $\lambda > 0$ is stable, then the corresponding modulated rotating wave or periodic solution is orbitally stable as well.\\
Since $A(t,\lambda)$ and $q(t,\lambda)$ are sufficiently $CS$ and $\Psi$ is a local diffeomorphism, we get that $\mathbf{\Phi}(t, u_{\lambda}, \lambda)$ is sufficiently $CS$. The same result can be inferred if we take into account that $u_{\lambda}=\Psi(q(0,\lambda))$.
\end{proof}

\begin{proof}[\textbf{Proof of Theorem \ref{thm:Aux1}}]
We will use the following result (see \cite{Ha}):
\begin{lemma}\label{lem:scalar}
Let $g(t,u)$ be a continuous function on an open connected set $[a_{1},b_{1}) \times \mathbb{R^{+}} \subset \Omega \subset \mathbb{R}^{2}$ and such that the initial value problem for the scalar equation $\dot{u}=g(t,u)$ has a unique solution $u(t) \geq 0$ on $t \in [a_{1},b_{1})$. If $f \colon [a_{1},b_{1}) \times \mathbb{R}^{n} \rightarrow \mathbb{R}^{n}$ is continuous and $\left \| f(t, x) \right \| \leq g(t, \left \| x \right \|)$ for $t \in [a_{1},b_{1})$ and $x \in \mathbb{R}^{n}$, then the solutions of $\dot{x}= f(t, x)$, $\left \| x(a_{1}) \right \| \leq u(a_{1})$ exists on $[a_{1},b_{1})$ and $\left \| x(t) \right \| \leq u(t)$ for $t \in [a_{1},b_{1})$.
\end{lemma}
There exists a constant $M > 0$ independent of $\lambda \geq 0$ small such that $\left | X^{G}(t,\lambda)\right | < M$ for $t \in [0,\infty)$ and $\lambda \geq 0$ small. For $\lambda \geq 0$ small and any $\left | Z \right | \leq \pi$, we have
\begin{eqnarray}
\left \| \left [ I_{3}+  \frac{1}{2}Z+ \left (\frac{1}{ \left | Z \right | ^{2}}-\frac{\cos \frac{\left | Z \right |}{2}} {2\sin \frac{\left | Z \right |}{2}\left| Z \right | }\right )Z^{2} \right ]\overrightarrow {X^{G}(t, \lambda)}\right \| & \leq & \nonumber \\
\left [ \left \| I_{3} \right \| +\frac{1}{2} \left \| Z \right \| +\left( \frac{1}{\left | Z \right | ^{2}}+\frac{\cos \frac{\left | Z \right |}{2}} {2\sin \frac{\left | Z \right |}{2}\left| Z \right |}\right ) \left \| Z \right \| ^{2} \right] M,
\end{eqnarray}

\begin{eqnarray}
\left \| \left [I_{3}+\frac{1}{2}Z+ \left (\frac{1}{\left | Z \right | ^{2}}-\frac{\cos \frac{\left | Z \right |}{2}} {2\sin \frac{\left | Z \right |}{2}\left| Z \right | }\right )Z^{2} \right ]\overrightarrow {X^{G}(t, \lambda)} \right \| & \leq & \nonumber \\
M \left [\sqrt{3}+\frac{1}{2}\sqrt{2}\left | Z \right |+ \left (\frac{1}{\left | Z \right | ^{2}}+\frac{\cos \frac{\left | Z \right |}{2}} {2\sin \frac{\left | Z \right |}{2}\left| Z \right|}\right )2\left | Z \right | ^{2}\right ]
\end{eqnarray}
or, if we use the fact that the function $x \rightarrow \frac{x\cos x}{\sin x}$ is decreasing on $[0,\frac{\pi}{2}]$,
\begin{eqnarray}
\left \| \left [I_{3}+\frac{1}{2}Z+ \left (\frac{1}{\left | Z \right | ^{2}}-\frac{\cos \frac{\left | Z \right |}{2}} {2\sin \frac{\left | Z \right |}{2}\left| Z \right | }\right )Z^{2} \right ]\overrightarrow {X^{G}(t, \lambda)} \right \| & \leq & \nonumber \\
M\left [\sqrt{3}+\frac{1}{2}\sqrt{2}\left | Z \right |+(2+2\cdot 1)\right ]\leq M(6+\left | Z \right |).
\end{eqnarray}
The initial value problem
\begin{equation}\label{E:zetaaux}
\begin{array}{lll}
\dot{a} &=& M(6+a),\\
a(0) &=& 0\\
\end{array}
\end{equation}
has solution on an maximal interval $[0, b_{max}]$, where $b_{max}$ is  defined such that $a(t) \leq \pi$ on $[0,b_{max}]$; $b_{max}$ is independent of $\lambda$.\\
Therefore, using Lemma \ref{lem:scalar} with $a_{1}=0$, it follows that the sufficiently $CS$ solution $Z_{1}(t, \lambda)$ of the initial value problem (\ref{E:ODEs_frequency}) is defined for any $t \in [0, b_{max}]$ and $\lambda \geq 0$ small, and that  $\left | Z(t, \lambda) \right | \leq \pi $ for any $t \in [0, b_{max}]$ and $\lambda \geq 0$ small.\\
If we choose a positive integer $n > 0$ such that $n > \frac{T(\lambda)}{ b_{max}}$, then the initial value problem
(\ref{E:ODEs_frequency}) has a solution $Z_{1}(t, \lambda)$ defined and sufficiently $CS$ on $t \in [0, \frac{T(\lambda)}{n}]$ and $\lambda \geq 0$ small. Since $T(\lambda)= T(0)+O(\lambda)$, then we can choose $n$ independent of $\lambda$ for $\lambda \geq 0$ small, where $T(0)=\frac{2\pi}{\left | \omega_{bif} \right |}$.\\
This ends the proof of (1) of Theorem \ref{thm:Aux1}.\\
To prove the formula (\ref{E:Zeta0}), we write the Taylor expansion of \\ $Z^{*}(t, \epsilon) \overset{\text{def}}=Z_{1}(t, \epsilon^{2})$:
\begin{eqnarray}
Z^{*}(t, \epsilon) &=& z_{00}(t, \epsilon)X_{0}^{1}+z_{11}(t,\epsilon)X_{1}+z_{22}(t, \epsilon)X_{2},\label{E:ZetaTaylor1}\\
z_{00}(t, \epsilon) &=& Z_{0}(t)+Z_{1}(t)\epsilon+ \epsilon^{2}Z_{2}(t, \epsilon). \label{E:ZetaTaylor2}
\end{eqnarray}
If we substitute (\ref{E:ZetaTaylor1}), (\ref{E:ZetaTaylor2}) and (\ref{E:InitialData}) in the first of the differential
equations given by (\ref{E:ODEs_frequency}) with $\lambda= \epsilon^{2}$, identify the orders of $\epsilon$ in both
sides and then take $\epsilon= \sqrt{\lambda}$, we get the formula (\ref{E:Zeta0}) for $Z_{1}(t, \lambda)$. For the sake of simplicity, we omit $\epsilon$.\\
The result is
\begin{eqnarray}\label{E:Tayloraux}
\dot{z}_{00}X_{0}^{1}+\dot{z}_{11}X_{1}+\dot{z}_{22}X_{2}= [I_{3}+\frac{z_{00}}{2}X_{0}^{1}+\frac{z_{11}}{2}X_{1}+\frac{z_{22}}{2}X_{2} \nonumber \\ +\left (\frac{1}{z_{00}^{2}+z_{11}^{2}+z_{22}^{2}}-\frac{\cos \frac{\sqrt{z_{00}^{2}+ z_{11}^{2}+z_{22}^{2}}}{2}}{2\sqrt{z_{00}^{2}+z_{11}^{2}+z_{22}^{2}}\sin \frac{\sqrt{z_{00}^{2}+ z_{11}^{2}+z_{22}^{2}}}{2}}\right ) \nonumber \\ \cdot (z_{00}X_{0}^{1}+z_{11}X_{1}+z_{22}X_{2})^{2} ]\overrightarrow{x_{0}X^{1}_{0}+ \epsilon x_{1}X_{1}+\epsilon x_{2}X_{2}}.\nonumber\\
\end{eqnarray}
Since $Z^{*}$ is sufficiently smooth, we get $Z^{*}(t, \epsilon)= Z^{*}(t,0)+ \epsilon H(t, \epsilon)= \left | X_{0} \right | t X_{0}^{1} + \epsilon H(t, \epsilon)$.\\
Therefore, $z_{11}(t,\epsilon)=\epsilon z_{1}(t, \epsilon)$ and $z_{22}(t,\epsilon)= \epsilon z_{2}(t, \epsilon)$. In the right-hand side of (\ref{E:Tayloraux}), the coefficient of $X_{0}^{1}$ is $x_{0}(t,\epsilon)+\epsilon^{2} k(t, \epsilon) =
\left | X_{0} \right |+\epsilon x_{01}(t)+\epsilon^{2} x_{02}(t, \epsilon)+ \epsilon^{2} k(t, \epsilon)$ and thus
$\dot{z}_{00}= \left | X_{0} \right |+\epsilon x_{01}(t)+\epsilon^{2} k_{1}(t, \epsilon)$, where $x_{0}$, $x_{01}$,
$x_{02}$ are defined in (\ref{E:InitialData}). Therefore, \[z_{00}(t,\lambda)= \left | X_{0} \right | t+
\epsilon \int_{0}^{t} x_{01}(s)\, ds +\epsilon^{2} x_{03}(t, \epsilon).\] This gives the formula (\ref{E:Zeta0}) after
we relabel $z_{1}(t,\sqrt{\lambda})$, $z_{2}(t,\sqrt{\lambda})$ and $x_{03}(t,\sqrt{\lambda})$
to $z_{1}(t,\lambda)$, $z_{1}(t,\lambda)$ and $x_{03}(t,\lambda)$. This ends the proof of (2) of Theorem \ref{thm:Aux1}.
\end{proof}

\begin{proof}[\textbf{Proof of Corollary \ref{cor:Aux2}}]
We have the same notations as in the proof of Theorem \ref{thm:Aux1}.
From the proof of Theorem \ref{thm:Aux1}, we know that the initial value problem
\begin{equation}\label{E:zetaaux1}
\begin{array}{lll}
\dot{a} &=& M(6+a),\\
a(0) &=& 0\\
\end{array}
\end{equation}
has the solution $a$ defined on $[0, b_{max}]$ with $a(t) \leq \pi$ on $[0,b_{max}]$.
Let $t_{0}=\frac{T(\lambda)}{n} \in [0,b_{max}]$. Consider the initial value problem
\begin{equation}\label{E:zetaaux2}
\begin{array}{lll}
\dot{a} &=& M(6+a),\\
a(t_{0}) &=& 0\\
\end{array}
\end{equation}
The solution to this problem satisfies $a(t) \leq \pi$ for any $t \in [t_{0}, 2t_{0}]$.\\
We repeat the previous argument for $t_{0}=i\frac{T(\lambda)}{n }$ for $i=2$,~$3$, \ldots, ~$n-1$. From this and Lemma
\ref{lem:scalar} with $a_{1}=t_{0}= i\frac{T(\lambda)}{n }$ for $i=2$,~$3$, \ldots, ~$n-1$, we get the conclusion
of Corollary \ref{cor:Aux2}.
\end{proof}

\begin{proof}[\textbf{Proof of Theorem \ref{thm:GlobalZeta}}]
Recall that $T(\lambda)=\frac{2\pi}{\left | \omega_{\lambda} \right |}$ for $\lambda \geq 0$ small.
\begin{enumerate}
\item [Part 1]
On the interval $[0, T(\lambda)]$:\\
Consider the initial value problem
\begin{equation}\label{E:auxHopf1}
\begin{array}{lll}
\dot{A} &=& AX^{G}(t, \lambda),\\
 A(0)   &=& I_{3}.
\end{array}
\end{equation}
Let us make the change of variable $A= e^{Z}$ \, near $I_{3}$ in the initial value problem (\ref{E:auxHopf1}).\\
Then, using Proposition \ref{prop:properties_Exp} (2) we get the following initial value problem in $Z$:\\
\begin{equation}\label{E:Zeta10}
\begin{array}{rll}
e^{Z}\sum_{n= 0}^{\infty}\frac{(-1)^{n}}{(n+1)!}(ad Z)^{n}\dot{Z} &= & e^{Z}X^{G}(t, \lambda),\\
Z(0) &=& O_{3},
\end{array}
\end{equation}
\begin{equation}\label{E:Zeta11}
\begin{array}{rll}
\sum_{n= 0}^{\infty}\frac{(-1)^{n}}{(n+1)!}(ad Z)^{n}\dot{Z} &=& X^{G}(t, \lambda),\\
Z(0) &=& O_{3},
\end{array}
\end{equation}
\begin{equation}\label{E:Zeta12}
\begin{array}{rll}
\overrightarrow{\sum_{n= 0}^{\infty}\frac{(-1)^{n}}{(n+1)!}(ad Z)^{n}\dot{Z}} &=& \overrightarrow{X^{G}(t, \lambda)},\\
Z(0) &=& O_{3},
\end{array}
\end{equation}
or, using again Proposition \ref{prop:geometrical_interp} (1) we get
\begin{equation}\label{E:Zeta13}
\begin{array}{rll}
\sum_{n= 0}^{\infty}\frac{(-1)^{n}}{(n+1)!}Z^{n}\overrightarrow{\dot{Z}} &=& \overrightarrow{X^{G}(t, \lambda)},\\
Z(0) &=& O_{3}.
\end{array}
\end{equation}
Using Proposition \ref{prop:geometrical_interp} (3), we have that
\begin{equation}\label{E:auxadjoint}
\sum_{n= 0}^{\infty}\frac{(-1)^{n}}{(n+1)!}Z^{n}= \left (I_{3}+\frac{\cos\left | Z \right | -1}{ \left | Z \right |
^{2}}Z+ \frac{\left | Z \right | -\sin\left | Z \right | } {\left | Z \right | ^{3}}Z^{2}\right ).
\end{equation}
(\ref{E:auxadjoint}) will proved later. \\
If we put (\ref{E:auxadjoint}) into (\ref{E:Zeta13}), we get
\begin{equation}\label{E:Zeta2}
\begin{array}{rll}
\left (I_{3}+\frac{\cos\left | Z \right | -1}{ \left | Z \right | ^{2}}Z+ \frac{\left | Z \right | -\sin\left | Z \right | } {\left | Z \right | ^{3}}Z^{2} \right )\overrightarrow{\dot{Z}} &=& \overrightarrow{X^{G}(t, \lambda)},\\
Z(0) &=& O_{3}.
\end{array}
\end{equation}
We will prove later that for any $\left | Z \right | < 2\pi$,
\begin{equation}\label{E:check1}
\left ( \frac{\cos\left | Z \right | -1}{ \left | Z \right | ^{2}}Z+ \frac{\left | Z \right | -\sin\left | Z \right | } {\left | Z \right | ^{3}}Z^{2}\right ) ^{-1}= I_{3}+\frac{1}{2}Z+\left (\frac{1}{\left | Z \right | ^{2}}-\frac{\cos \frac{\left | Z \right |}{2}}{2\sin \frac{\left | Z \right |}{2}\left| Z \right |} \right )Z^{2}.
\end{equation}
Since we are looking for $\left | Z \right | \leq \pi < 2\pi$, the system (\ref{E:Zeta2}) becomes
\begin{equation}\label{E:Zeta4}
\begin{array}{lll}
\overrightarrow{\dot{Z}} &=& \left[ I_{3}+\frac{1}{2}Z+ \left (\frac{1}{\left | Z \right | ^{2}}-\frac{\cos \frac{\left | Z \right |}{2}}{2\sin \frac{\left | Z \right |}{2}\left| Z \right |}\right )Z^{2}\right ] \overrightarrow{X^{G}(t, \lambda)},\\
Z(0) &=& O_{3},
\end{array}
\end{equation}
where we are looking for a solution $Z$ such that $\left | Z \right | \leq \pi$ on some interval $[0,t_{0}(\lambda)]$.\\
We now prove (\ref{E:check1}). It is enough to check that
\begin{equation}\label{E:c1}
\left (I_{3}+\frac{\cos\left | Z \right | -1}{ \left | Z \right | ^{2}}Z+\frac{\left | Z \right | -\sin\left | Z \right | } {\left | Z \right | ^{3}}Z^{2}\right )\left [I_{3}+\frac{1}{2}Z+\left (\frac{1}{\left | Z \right | ^{2}}-
 \frac{\cos \frac{\left | Z \right |}{2}}{2\sin \frac{\left | Z \right |}{2}\left| Z \right |} \right )Z^{2} \right ]=I_{3}.
\end{equation}
Using Proposition \ref{prop:geometrical_interp} (3), that is $Z^{3}=-\left | Z \right | ^{2} Z$ and $Z^{4}=-\left | Z \right | ^{2} Z^{2}$, we have
\begin{multline}
\left (I _{3}+\frac{\cos \left | Z \right | -1}{ \left | Z \right | ^{2}}Z+\frac{\left | Z \right | -\sin\left | Z \right | } {\left | Z \right | ^{3}}Z^{2}\right )\left [I_{3}+\frac{1}{2}Z+\left (\frac{1}{\left | Z \right | ^{2}}-
 \frac{\cos \frac{\left | Z \right |}{2}}{2\sin \frac{\left | Z \right |}{2}\left| Z \right |} \right )Z^{2} \right ]\\
\shoveleft{ =I_{3}+\frac{1}{2}Z+\left (\frac{1}{\left | Z \right | ^{2}}-\frac{\cos \frac{\left | Z \right |}{2}}{2\sin \frac{\left | Z \right |}{2}\left | Z \right |} \right )Z^{2}+ \frac{\cos\left | Z \right | -1}{ 2\left | Z \right | ^{2}}Z^{2}}\\
\shoveleft{+\frac{\cos\left | Z \right | -1}{ \left | Z \right | ^{2}}Z+\frac{\cos \left | Z \right | -1}{ \left | Z \right | ^{2}} \left (\frac{1}{\left | Z \right | ^{2}}-\frac{\cos \frac{\left | Z \right |}{2}}{2\sin \frac{\left | Z \right |}{2}\left | Z \right |}\right )Z^{3} +\frac{\left | Z \right | -\sin\left | Z \right | }{\left | Z \right | ^{3}}Z^{2}+ \frac{\left | Z \right | -\sin\left | Z \right | }{2\left | Z \right | ^{3}}Z^{3}}\\ \shoveleft{+\frac{\left | Z \right | -\sin\left | Z \right | }{\left | Z \right | ^{3}} \left (\frac{1}{\left | Z \right | ^{2}}-\frac{\cos \frac{\left | Z \right |}{2}}{2\sin \frac{\left |
Z \right |} {2}\left| Z \right |}\right )Z^{4}}\\
\shoveleft{=I_{3}+\left [ \frac{1}{2}+\frac{\cos\left | Z \right |-1}{ \left | Z \right | ^{2}}- \frac{\cos\left | Z \right | -1}{ \left | Z \right | ^{2}}+\frac{\cos \left | Z \right | -1}{ \left | Z \right | ^{2}} \frac{\left | Z \right |\cos \frac{\left | Z \right |}{2}}{2\sin \frac{\left | Z \right |}{2}}-\frac{1}{2}+ \frac{\sin \left | Z \right |}{2\left | Z \right |} \right ] Z}\\
\shoveleft{\left [ \frac{1}{\left | Z \right | ^{2}}-\frac{\cos \frac{\left | Z \right |}{2}}{2\sin \frac{\left | Z \right |}{2}\left | Z \right |}+ \frac{\cos\left | Z \right | -1}{ 2\left | Z \right | ^{2}}+ \frac{\left | Z \right | -\sin\left | Z \right | }{\left | Z \right | ^{3}}- \frac{\left | Z \right | -\sin\left | Z \right | }{\left | Z \right | } \frac{1}{\left | Z \right | ^{2}}+\right . }\\
\shoveleft {\left . + \frac{\left | Z \right | -\sin\left | Z \right | }{\left | Z \right | }\frac{\cos \frac{\left | Z \right |}{2}}{2\sin \frac{\left | Z \right |}{2}\left | Z \right | } \right ]Z^{2}} \\
\shoveleft{ =I_{3}+\left [ \frac{-\sin^{2} \frac{\left | Z \right |}{2}}{\left | Z \right |}\frac{\cos\frac{\left | Z \right |}{2}}{\sin \frac{\left | Z \right |}{2}}+ \frac{\sin \left | Z \right |}{2\left | Z \right |} \right ]Z +\left [ \frac{\cos ^{2}\frac{\left | Z \right |}{2}}{\left | Z \right | ^{2}}- \frac{2\sin \frac{\left | Z \right |}{2}\cos \frac{\left | Z\right |}{2}}{2\left | Z \right |^{2}} \frac{\cos \frac{\left | Z \right |}{2}}{\sin \frac{\left | Z \right |}{2}}\right ]Z^{2}=I_{3}}.
\end{multline}
The proof of (\ref{E:auxadjoint}) follows by using Proposition \ref{prop:geometrical_interp} (3):
\begin{equation}\label{E:aux11}
\begin{split}
\sum_{n=0}^{\infty}\frac{(-1)^{n}}{(n+1)!}Z^{n} &=I_{3}+\sum_{k= 1}^{\infty}\frac{(-1)^{2k}}{(2k+1)!}Z^{2k}+\sum_{k= 0}^{\infty}\frac{(-1)^{2k+1}}{(2k+2)!}Z^{2k+1}\\ &=I_{3}+\left (-\frac{1}{2!}+\frac{1}{4!}\left | Z \right |^{2}-\frac{1}{6!}\left | Z \right | ^{4}+\ldots \right )Z\\
&+\left (\frac{1}{3!}-\frac{1}{5!}\left | Z \right |^{2}+\frac{1}{7!}\left | Z \right |^{4}-\ldots \right )Z^{2}\\
&=I_{3}+\frac{\cos\left | Z \right | -1}{ \left | Z \right | ^{2}}Z+\frac{\left | Z \right | -\sin\left | Z \right | } {\left | Z \right | ^{3}}Z^{2},
\end{split}
\end{equation}
where we have used the Taylor expansions for \emph{sine} and \emph{cosine}.\\
We apply Theorem \ref{thm:Aux1} to get $A(t, \lambda)= e^{Z_{1}(t,\lambda)}$, for any $t \in [0, \frac{T(\lambda)}{n}]$
and $\lambda \geq 0$ small, where $Z_{1}(t, \lambda)$ is the sufficiently $CS$ solution of (\ref{E:ODEs_frequency}).
We make the change of variable $B= (A(\frac{T(\lambda)}{n}, \lambda))^{-1}A$. We get that the solution of the initial
value problem (\ref{E:auxHopf1}) is given by \[A(t, \lambda)= e^{Z_{1}(\frac{T(\lambda)}{n}, \lambda)}B(t, \lambda),\]
where $B(t, \lambda)$ is the solution of the initial value problem
\begin{equation}\label{E:auxHopf2}
\begin{array}{lll}
\dot{B} &=& BX^{G}(t,\lambda),\\
 B(\frac{1}{n}T(\lambda)) &=& I_{3}.
\end{array}
\end{equation}
Using the above argument and Corollary \ref{cor:Aux2}, we have that $B(t, \lambda)= e^{Z_{2}(t, \lambda)}$ for any
$t \in [\frac{T(\lambda)}{n}, 2\frac{T(\lambda)}{n}]$ and $\lambda \geq 0$ small, where $Z_{2}(t, \lambda)$ is the
sufficiently $CS$ solution of (\ref{E:ODEs_frequency1}) for $i= 2$.\\
If we continue this, we get for $\lambda \geq 0$ small and any $t \in [(i-1)\frac{T(\lambda)}{n},i\frac{T(\lambda)}{n}]$ , \[A(t, \lambda)= e^{Z^{0}(t, \lambda)},\] where  $Z^{0}(t, \lambda)$ is given by
\begin{equation}\label{E:ZetaFormula0}
e^{Z^{0}(t, \lambda)}= e^{Z_{1}(\frac{T(\lambda)}{n}, \lambda)}e^{Z_{2}(2\frac{T(\lambda)}{n}, \lambda)} \ldots
e^{Z_{i}(t, \lambda)}
\end{equation}
for any $ t \in [(i-1)\frac{T(\lambda)}{n},i\frac{T(\lambda)}{n}]$ and small $\lambda \geq 0$ for $i= 1$,~$2$, \ldots,
~$n$, where $Z_{i}$ is the solution of the initial value problem (\ref{E:ODEs_frequency1}) on the interval $[(i-1)\frac{T(\lambda)}{n}, i\frac{T(\lambda)}{n}]$.\\
Then $[Z(t, \lambda)]$ is defined by formula (\ref{E:ZetaFormula00}).\\
Also, from (\ref{E:Zeta0}) in Theorem \ref{thm:Aux1} we get
\begin{multline}\label{E:ZetaIt}
Z_{i}(t, \lambda)= \left( \left | X_{0} \right | (t-(i-1)\frac{T(\lambda)}{n})+\sqrt{\lambda} \int_{(i-1)\frac{T(\lambda)}{n}}^{t}x_{01}(s) \, d s + \lambda x_{03}^{i}(t,\lambda)\right )X_{0}^{1}\\+\sqrt{\lambda} z_{1}^{i}(t, \lambda)X_{1}+\sqrt{\lambda}z_{2}^{i}(t, \lambda)X_{2}
\end{multline}
for any $t \in [(i-1)\frac{T(\lambda)}{n}, i\frac{T(\lambda)}{n}]$ and small $\lambda \geq 0$ for $i= 1$,~$2$, \ldots, ~$n$.\\
Because the $BCH$ formula in $so(3)$ is smooth from $so(3) \times so(3)$ into $D$ (see Theorem \ref{thm:BCH_formula}) and $Z_{i}(t,\lambda)$ for $i= 1$,~$2$, \ldots, ~$n$ are sufficiently $CS$, then  $[Z(t, \lambda)]$ is sufficiently $CS$ on $[0,T(\lambda)]$ and $\lambda \geq 0$ small.
\item [Part 2]
On the interval $[iT(\lambda), (i+1)T(\lambda)]$ for any integer $i \geq 0$.\\
Since \[A(t+T(\lambda), \lambda)= e^{X(\lambda)T(\lambda)}A(t, \lambda)= A(T(\lambda),\lambda)A(t,\lambda)=
e^{Z^{0}(T(\lambda), \lambda)}e^{Z^{0}(t, \lambda)},\] it is easy to see that we can define for any $t \in [T(\lambda), 2T(\lambda)]$ and $ \lambda \geq 0$ small, \[[Z^{1}(t,\lambda)]=BCH(Z^{0}(T(\lambda),\lambda), Z^{0}(t-T(\lambda), \lambda)).\] Because the $BCH$ formula in $so(3)$ is smooth from $so(3) \times so(3)$ into $D$, $Z^{0}(t,\lambda)$ is sufficiently $CS$ and $T(\lambda)$ is sufficiently smooth, then $[Z(t, \lambda)]$ is sufficiently $CS$ for any $t \in [T(\lambda), 2T(\lambda)]$ and for $\lambda \geq 0$ small.\\
We then repeat the above argument.
\end{enumerate}
It is clear that the function $[Z(t,\lambda)]$ is sufficiently $CS$, and $A(t,\lambda)= e^{Z(t,\lambda)}$ for
$t \in [0,\infty)$ and $\lambda \geq 0$ small. Let us define $Z_{1}(t, \lambda)=q([Z(t, \lambda)])$.
Therefore, $Z_{1}(t, \lambda)$ is sufficiently $CS$ and $A(t,\lambda)= e^{Z(t,\lambda)}=e^{[Z(t,\lambda)]}=
e^{Z_{1}(t,\lambda)}$.\\
Thus, $Z_{1}(t, \lambda)$ is the solution of the initial value problem (\ref{E:ODEs_frequency_final}).

\end{proof}

\begin{proof}[\textbf{Proof of Proposition \ref{prop:Frequency1}}]
If we define
\begin{equation}\label{E:freqfreq}
X(\lambda)= \frac{1}{T(\lambda)}Z(T(\lambda), \lambda) \mbox{ for } \lambda \geq 0 \mbox{ small },
\end{equation}
then $A(T(\lambda), \lambda)= e^{X(\lambda)T(\lambda)}$.\\
Since $Z(t,\lambda)$ is sufficiently $CS$ for $ t \in [0,\infty)$ and $\lambda \geq 0$ small and $T(\lambda) > 0$ for
$\lambda \geq 0$ is sufficiently smooth, we get that $X(\lambda)$ is sufficiently $CS$ for $\lambda \geq 0$ small.
The branch $X(\lambda)$ defined by (\ref{E:freqfreq}) is such that  $X(\lambda)T(\lambda) \in q(D)$ for $\lambda \geq 0$
and it is the branch whose existence is stated in Theorem \ref{thm:FrequencyGeneral1}.\\
\end{proof}

\begin{proof}[\textbf{Proof of Corollary \ref{cor:Frequency2}}]
Using Proposition \ref{prop:Frequency1} we know that $X(\lambda)$ is sufficiently $CS$ for $\lambda \geq o$ small.
Therefore, $\left | X(\lambda) \right | ^{2} $ for $\lambda \geq 0$ small is sufficiently $CS$.
Using Taylor formula for $\epsilon \rightarrow \left | X(\epsilon^{2}) \right | ^{2} $ around $\epsilon= 0$ up to the
term of order $\epsilon$, and then taking $\epsilon= \sqrt{\lambda}$, we get $\left | X(\lambda) \right | ^{2}=
\left | X(0) \right | ^{2} + O(\epsilon)= \left | X(0) \right | ^{2}+ O(\sqrt{\lambda})$ for $\lambda \geq 0$ small.
Also, $\left | X(\lambda) \right |$ is continuous for $\lambda \geq 0$.\\
The branch $X^{f}(\lambda)$ is constructed below. Let $\left | X_{0} \right |T(0)= \alpha_{0}+2k\pi$, where $\alpha_{0}
\in [0,2\pi)$, $k \in \mathbb{Z}$, $k \geq 0$, and $T(0)=\frac{2\pi}{\left | \omega_{bif}\right |}$.\\
\begin{enumerate}
\item If $\left | X_{0} \right | \neq k\omega_{bif}$ for all $k \in \mathbb{Z}$, then $\alpha_{0} \neq 0$.\\
Since $A(T(0),0)=e^{X_{0}T(0)}=e^{X(0)T(0)}=e^{Z(t, \lambda)}$, we get that\\
$X(0)T(0)=Z(t,0)=q(BCH(0, X_{0}T(0)))=\alpha_{0}X_{0}^{1}$.
Therefore, $\left | X(0) \right | =\alpha_{0} >0$. Since $X(\lambda)$ is sufficiently $CS$, then
$\left | X(\lambda) \right | >0$ for $\lambda \geq 0$ small.\\
Then $X^{f}(\lambda)=\frac{\left | X(\lambda) \right |+k\left |\omega_{\lambda}\right | }{\left | X(\lambda) \right | }
X(\lambda)$ is well-defined for $\lambda \geq 0$ small.\\
We check that $X^{f}(0)=X_{0}$. Using $X(0)=\left | X(0) \right | X_{0}^{1}$, it yields
$X^{f}(0)=\frac{\left | X(0) \right |T(0)+2k\pi}{T(0)}X_{0}^{1}=\frac{\alpha_{0}+2k\pi}{T(0)}X_{0}^{1}=
\left | X_{0} \right | X_{0}^{1}=X_{0}$.\\
Since $\left | X(\lambda) \right | >0$ and $X(\lambda)$, $\omega_{\lambda}$ are sufficiently $CS$, then $X^{f}(\lambda)$
is sufficiently $CS$.\\
Using the Taylor formula for $\epsilon \rightarrow X^{f}(\epsilon^{2})$ around $\epsilon= 0$ up to the term of order
$\epsilon$ and then taking $\epsilon= \lambda^{\frac{1}{2}}$, we get $X^{f}(\lambda)=
X^{f}(0)+O(\epsilon)X_{0}^{1}+O(\epsilon)X_{1}+O(\epsilon)X_{2}=
X_{0}+O(\lambda^{\frac{1}{2}})X_{0}^{1}+O(\lambda^{\frac{1}{2}})X_{1}+
O(\lambda^{\frac{1}{2}})X_{2}$ for $\lambda \geq 0$ small, thus we get the formula (\ref{E:freqTaylor01}).\\
We have only modulated rotating waves because $\left | X_{0} \right | T(0) \neq 2k\pi$ for any $k \in \mathbb{Z}$ implies
$\left | X^{f}(\lambda) \right | T(\lambda) \neq 2k\pi$ for any $k \in \mathbb{Z}$ and for $\lambda \geq 0$ small by the
continuity of $\left | X^{f}(\lambda) \right | T(\lambda)$.\\
\item If $\left | X_{0}\right | = k\omega_{bif}$ for some $ k \in \mathbb{Z}$, $k \neq 0$, then $\alpha_{0} = 0$.\\
Let $X^{f}(0)=X_{0}$.\\
If $\left | X(\lambda) \right | \neq 0$, then we define $X^{f}(\lambda)= \frac{\left | X(\lambda) \right |+
 k\left |\omega_{\lambda}\right |}{\left | X(\lambda) \right | }X(\lambda)$ for $\lambda > 0$ small.\\
If $\left | X(\lambda) \right | = 0$, then we define $X^{f}(\lambda)= X(\lambda)+ k\left | \omega_{\lambda}\right |
Q(\lambda)$, where $Q(\lambda) \in so(3)$, $\left | Q(\lambda)\right | =1$ for $\lambda > 0$ small.\\
If $\left | X_{0} \right | = k\omega_{bif}$ for some $k \in \mathbb{Z}$, $k \neq 0$, then we check that
$\left | X^{f}(\lambda) \right | = \left | X(\lambda) \right | + \frac{2k\pi}{T(\lambda)}$ for all $\lambda \geq 0$ small.
For $\lambda > 0$ small such that $\left | X(\lambda) \right | > 0$, we have $\left | X^{f}(\lambda) \right |=
\left | X(\lambda) \right |+ 2k\pi \frac{\left | X(\lambda) \right |}{T(\lambda)\left | X(\lambda)\right |}=
\left | X(\lambda) \right |+ \frac{2k\pi}{T(\lambda)}$.\\
For $\lambda > 0$ small such that $X(\lambda) = 0$, we have $\left | X^{f}(\lambda) \right |=
\frac{2k\pi}{T(\lambda)}\left | Q(\lambda) \right |= \frac{2k\pi}{T(\lambda)}$.\\
Hence, for $\lambda \geq 0$ small, we get $\left | X^{f}(\lambda) \right |= \left | X(\lambda) \right |+
\frac{2k\pi}{T(\lambda)}$.\\
Since $\left | X(\lambda)\right |^{2}= \left | X(0) \right |^{2}+O(\sqrt{\lambda})$ and $\left | X(0) \right | =0$,
we get  $\left | X(\lambda)\right | = O(\lambda^{\frac{1}{4}}$) and then $\left | X^{f}(\lambda) \right |=
O(\lambda^{\frac{1}{4}})+\frac{2 k \pi}{T(0)}+O(\lambda)= O(\lambda^{\frac{1}{4}})+\left | X_{0} \right |$ for
$\lambda \geq 0$ small, thus we get the formula (\ref{E:freqTaylor02}). \\
Also, $\left | X^{f}(\lambda) \right | $ is continuous since $X(\lambda)$ and $T(\lambda)$ are continuous for
$\lambda \geq 0$ small.
\end{enumerate}

By using Proposition \ref{prop:properties_Exp} ((6) and (7)), we have $e^{X^{f}(\lambda)T(\lambda)}=
e^{X(\lambda)T(\lambda)}$ for $\lambda \geq 0$ small.
\end{proof}

\begin{proof}[\textbf{Proof of Theorem \ref{thm:Periodic_part}}]
\begin{enumerate}
\item Suppose that $\left | X_{0} \right | \neq k\omega_{bif}$ for any $k \in \mathbb{Z}$.\\
We have $B^{f}(t,\lambda)= e^{-X^{f}(\lambda)t}A(t, \lambda)= e^{-X^{f}(\lambda)t}e^{Z(t, \lambda)}= e^{q(BCH(-X^{f}(\lambda)t,
Z(t, \lambda)))}$ for $t \in [0, \infty)$ and $\lambda \geq 0$ small.\\
Since $X^{f}(\lambda)$ and $Z(t, \lambda)$ are sufficiently $CS$ for $t \in [0,\infty)$ and $\lambda \geq 0$ small
and the $BCH$ formula in $so(3)$ is smooth, it follows that $Per^{f}(t,\lambda)=q(BCH(-X^{f}(\lambda)t, Z(t, \lambda)))$
is sufficiently $CS$ for $t \in [0, \infty)$ and $\lambda \geq 0$ small.\\
Since $e^{-X^{f}(0)t}e^{Z(t,0)}= A(t)^{-1}A(t)= I_{3}= e^{BCH(-X^{f}(0)t,Z(t,0))}$, it follows that
$Per^{f}(t,0)= q(BCH(-X^{f}(0)t,Z(t,0)))= O_{3}$. Then, $Per^{f}(t,\lambda)=\sqrt{\lambda}Y(t, \lambda)$ for all
$t \in \mathbb{R}$ and $\lambda \geq 0$ small.
The periodicity of $Y(.,\lambda)$ results from the following Remark:
\begin{remark}
If $P(t,\lambda)=e^{G(t,\lambda)}$, $ P \colon \mathbb{R} \times \mathbb{R} \rightarrow SO(3)$ is a continuous $T$-periodic
function in $t$ for $\lambda \geq 0$ small and $G \colon \mathbb{R} \times \mathbb{R} \rightarrow so(3)$ is a continuous
function such that $G(t,\lambda)=\lambda H(t,\lambda)$, then $H(t,\lambda)$ is a $T$-periodic function in $t$ for $\lambda \geq
0$ small.
\begin{proof}
We have $P(t+T,\lambda)=P(t,\lambda) \Longrightarrow e^{G(t+T,\lambda)}=e^{G(t,\lambda)} \Longrightarrow e^{\lambda
H(t+T)}=e^{\lambda H(t)}$ and $\lambda H(t,\lambda)$, $\lambda H(t+T,\lambda)$ are in a neighborhood of $O_{3}$ for $\lambda \geq
0$ small. Then, by using the fact that the exponential map $exp$ is a local diffeomorphism at $O_{3}$, it follows that $\lambda
H(t,\lambda)= \lambda H(t+T,\lambda)$ or $G(t+T,\lambda)=G(t,\lambda)$ for $\lambda \geq 0$ small.
\end{proof}
\end{remark}
Writing the Taylor formula for $\epsilon \rightarrow Per^{f}(t,\epsilon^{2})$ at $\epsilon= 0$ up to order of
$\epsilon$ and taking $\epsilon=\sqrt{\lambda}$, we get $Per^{f}(t,\lambda)= O_{3}+\epsilon Y_{1}(t,\epsilon)=
\sqrt{\lambda}Y(t,\lambda))$ for $t \in [0, \infty)$ and $\lambda \geq 0$ small.
\item Suppose that $\left | X_{0} \right | = k\omega_{bif}$ for some $k \in \mathbb{Z}$, $k \neq 0$.\\
We have $B(t,\lambda)= e^{-X(\lambda)t}A(t, \lambda)= e^{-X(\lambda)t}e^{Z(t, \lambda)}=
e^{q(BCH(-X(\lambda)t, Z(t,\lambda)))}$ for $t \in [0, \infty)$ and $\lambda \geq 0$ small.
Since $X(\lambda)$ and $Z(t, \lambda)$ are sufficiently $CS$ for $t \in [0, \infty)$ and $\lambda \geq 0$ small,
and $BCH$ formula in $so(3)$ is smooth, it is clear that $Per(t,\lambda)=q(BCH(-X(\lambda)t, Z(t,\lambda)))$ is sufficiently $CS$ for $t \in [0, \infty)$
and $\lambda \geq 0$ small.\\
Since $X(0)=O_{3}$, we get $X(\lambda)=O(\sqrt{\lambda})X_{0}^{1}+O(\sqrt{\lambda})X_{1}+O(\sqrt{\lambda})X_{2}$.
Since $A(t,0)=e^{X_{0}t}=e^{Z(t,0)}$ and $Per(t,0)=q(BCH(O_{3},Z(t,0)))=q(BCH(O_{3}, X_{0}t))$, it yields that
$Per(t,0)=X_{0}t (\mbox{ mod } 2\pi)$. We denote $r(t)=t (\mbox { mod  } \frac{2\pi}{\left | X_{0} \right |})$.
Therefore, $Per(t, \lambda)= X_{0}r(t)+\sqrt{\lambda}H_{1}(t, \lambda)$ for
$t \in [0, \infty)$ and $\lambda \geq 0$ small. Also, it is clear that $B(t, \lambda)=e^{X_{0}t+\sqrt{\lambda}H(t, \lambda)}$.
Since $B(t, \lambda)$ is $\frac{2\pi}{\left | \omega_{\lambda} \right |} $-periodic in $t$ and $B(0, \lambda)=I_{3}$,
it follows that
$e^{X_{0}\frac{2\pi}{\left | \omega_{\lambda} \right | }+\sqrt{\lambda}H(\frac{2\pi}{\left | \omega_{\lambda}\right | },\lambda)}=I_{3}$ for $\lambda \geq 0$ small.
\end{enumerate}
\end{proof}

\begin{proof}[\textbf{Proof of Theorem \ref{thm:MRW}}]
Since $\mathbf{\Phi}(t, u_{\lambda}, \lambda)=e^{X(\lambda)}B(t,\lambda)\Psi(q(t, \lambda))=
e^{X^{f}(\lambda)}B^{f}(t,\lambda)\Psi(q(t, \lambda))$, the conclusion of Theorem \ref{thm:MRW} result by applying
Theorems \ref{prop:Frequency1}, \ref{cor:Frequency2} and \ref{thm:Periodic_part}.
\end{proof}

\begin{proof}[\textbf{Proof of Theorem \ref{thm:Resonant_drift1}}]
Let us define $Y(\alpha, \beta)= e^{\alpha X_{1}+\beta X_{2}}$.\\
Let us denote $\lambda= \epsilon^{2}$, $\widetilde{X}(\epsilon, \mu)= X(\epsilon^{2}, \mu)$ and $\widetilde{T}(\epsilon, \mu)= T(\epsilon^{2}, \mu)= \frac{2\pi}{\left | \omega_{\epsilon^{2}, \mu}\right |}$. Let $\widetilde{T}(0, \mu)= \frac{2\pi}{\left | \omega_{0,\mu}\right | }$.\\
We have to prove that there exist sufficiently smooth curves
\[\mu= \mu(\epsilon), \alpha= \alpha(\epsilon) \mbox{ and } \beta= \beta(\epsilon)\mbox{ such that } \mu(0)= 0,
\alpha(0)= 0 \mbox{ and } \beta(0)= 0\] satisfying
\begin{equation}\label{E:alphabeta1}
Y(\alpha, \beta)= e^{\widetilde{X}(\epsilon, \mu)\widetilde{T}(\epsilon, \mu)} \mbox{ or } Y(\alpha, \beta)^{-1}e^{\widetilde{X}(\epsilon, \mu)\widetilde{T}(\epsilon,\mu)}= I_{3}.
\end{equation}
Let us define the function $F \colon \mathbb{R}^{4} \rightarrow SO(3)$ by
\[F(\alpha, \beta, \epsilon, \mu)= e^{-\alpha X_{1}-\beta X_{2}}e^{\widetilde{X}(\epsilon, \mu)
\widetilde{T}(\epsilon, \mu)}.\] We will prove the existence of the smooth curves using the implicit
function theorem for the function $F$ at the point $(\alpha, \beta, \epsilon, \mu)= (0, 0, 0, 0)$.\\
We now prove that:
\begin{enumerate}
\item $F(0,0,0,0)=I_{3}$\\ and \item
$[(DF)_{(0,0,0,0)}]_{X_{0}^{1},X_{1},X_{2}}$ has rank 3.
\end{enumerate}
We have that $F(0, 0, 0, 0)= e^{-0 \cdot X_{1}-0 \cdot X_{2}}e^{\widetilde{X}(0, 0) \widetilde{T}(0,0)}=
e^{X(0,0)T(0,0)}= e^{X_{0}\frac{2\pi}{\left | \omega_{bif} \right | }}= I_{3}$ since $\left | X_{0} \right |=
k\omega_{bif}$. \\
Recall that we have \[X_{G}(q, \lambda, \mu)= x_{0}(q, \lambda, \mu)X_{0}^{1}+ x_{1}(q, \lambda, \mu)X_{1}+
_{2}(q, \lambda, \mu)X_{2}.\] Then $X_{G}(0,0, \mu)= a(\mu)X_{0}^{1}+b(\mu)X_{1}+c(\mu)X_{2}$, where
\[a(\mu)= x_{0}(0, 0, \mu)=\left | X_{0} \right |+O(\mu), b(\mu)= x_{1}(0, 0, \mu)=O(\mu) \mbox{ and } c(\mu)=
x_{2}(0, 0, \mu)=O(\mu).\] Since $A(t,\lambda,\mu)=e^{X(\lambda, \mu)T(\lambda, \mu)}$, then
$A(t,0,\mu)=e^{X(0, \mu)T(0, \mu)}= e^{\widetilde{X}(0, \mu)\widetilde{T}(0, \mu)}$.
Also $A(t,0,\mu)=e^{X_{G}(0,0,\mu)T(0,\mu)}$. Using Proposition \ref{prop:properties_Exp} (2), we get
\begin{multline}\label{E:eq22} \frac{\partial}{\partial{\mu}} e^{X(0, \mu)T(0, \mu)} |_{\mu= 0}=
\frac{\partial}{\partial{\mu}} e^{X_{G}(0,0, \mu)T(0, \mu)} |_{\mu= 0} \\
= e^{X_{G}(0,0,0)T(0,0)} \sum_{n= 0}^{\infty}\left [\frac{(-1)^{n}}{(n+1)!} \left ( ad(X_{G}(0,0)T(0,0))\right )^{n}\left ( (X_{G}(0,0,\mu)T(0,\mu))^{'} |_{\mu= 0} \right )\right ].
\end{multline}
We have $T(0,0)=\frac{2\pi}{\left | \omega_{bif}\right |}$, $X_{G}(0,0,0)=X_{0}$ and
\begin{equation}
\begin{split}
(X_{G}(0,0,\mu)T(0,\mu))^{'}|_{\mu= 0} &= (X_{G}(0,0,\mu))^{'}|_{\mu=0} T(0,0)+ X_{G}(0,0,0)(T(0,\mu))^{'}|_{\mu=0}\\
&=\frac{2\pi}{\left | \omega_{bif}\right |} \left [ (a^{'}(0) X_{0}^{1}+b^{'}(0) X_{1}+c^{'}(0) X_{2}) \right ]\\&+
\left | X_{0} \right | \left [ -2\pi \cdot sgn(\omega_{bif})\right ]\frac{(\omega_{0,\mu})^{'}|_{\mu= 0}}{\omega_{bif}^{2}}X_{0}^{1}\\
&= \frac{2\pi}{\left | \omega_{bif} \right |} \left [(a^{'}(0)-k\cdot (sgn(\omega_{bif}))^{2}(\omega_{0,\mu})
^{'}|_{\mu= 0}) X_{0}^{1}+b^{'}(0) X_{1}+c^{'}(0) X_{2} \right ],
\end{split}
\end{equation}
where we have used $\left | X_{0} \right |=k\omega_{bif}$. Therefore, (\ref{E:eq22}) becomes
\begin{equation}
\begin{split}
& e^{\frac{2\pi \left | X_{0} \right |}{\left | \omega_{bif} \right | }X_{0}^{1}} \sum_{n= 0}^{\infty}\frac{(-1)^{n}}
{(n+1)!}\left ( ad(\frac{2\pi\left | X_{0} \right |}{\left | \omega_{bif}\right |}X_{0}^{1})\right ) ^{n}
\left (\frac{2\pi} {\left | \omega_{bif}\right |} \left [(a^{'}(0)- k (\omega_{0,\mu})^{'}|_{\mu= 0})X_{0}^{1}+
b^{'}(0)X_{1}+ c^{'}(0)X_{2}\right ]\right )\\ &=\frac{2\pi}{\left | \omega_{bif} \right | }
\left [(a^{'}(0)- k (\omega_{0,\mu})^{'}|_{\mu= 0})X_{0}^{1}+b^{'}(0)X_{1}+c^{'}(0)X_{2}\right ]\\
&+\sum_{n= 1}^{\infty}{\frac{(-1)^{n}}{(n+1)!} \left( ad(2| k |\pi X_{0}^{1})\right )^{n}\left (\frac{2\pi}{\left |
\omega_{bif}\right |} \left [(a^{'}(0)- k (\omega_{0,\mu})^{'}|_{\mu= 0})X_{0}^{1}+b^{'}(0)X_{1}+ c^{'}(0)X_{2}\right ]
\right )}\\ &=\frac{2\pi}{\left | \omega_{bif} \right | }\left [(a^{'}(0)- k (\omega_{0,\mu})^{'}|_{\mu= 0})X_{0}^{1}+
(b^{'}(0)+*) X_{1}+(c^{'}(0)+*)X_{2}\right ],
\end{split}
\end{equation}
where we have used $e^{\frac{2\pi\left | X_{0} \right |}{\left | \omega_{bif}\right |}X_{0}^{1}}=I_{3}$ and where
$*$ denotes other terms. These term appear only in the coefficients of $X_{1}$ and $X_{2}$ since
\begin{equation}
\begin{split}
& \left( ad(2| k | \pi X_{0}^{1})\right ) \left (\frac{2\pi}{\left | \omega_{bif}\right |} \left [(a^{'}(0)-
k (\omega_{0,\mu})^{'}|_{\mu= 0})X_{0}^{1}+b^{'}(0)X_{1}+ c^{'}(0)X_{2}\right ]\right )\\ &=
\frac{2\pi}{\left | \omega_{bif} \right | }\left [2| k| \pi X_{0}^{1}, (a^{'}(0)-k(\omega_{0,\mu})^{'}|_{\mu= 0})
X_{0}^{1} +b^{'}(0)X_{1}+c^{'}(0)X_{2}\right ] \\ &= \frac{4| k |\pi^{2}}{\left | \omega_{bif}\right |}(-c^{'}(0)X_{1}+
b^{'}(0)X_{2}),
\end{split}
\end{equation}
and for any integer $n \geq 2$,
\begin{multline}
\left( ad(2| k |\pi X_{0}^{1})\right )^{n}\left (\frac{2\pi}{\left | \omega_{bif}\right |} \left [(a^{'}(0)-
k(\omega_{0,\mu})^{'}|_{\mu= 0})X_{0}^{1}+b^{'}(0)X_{1}+c^{'}(0)X_{2}\right ]\right )\\=
\left [2| k |\pi X_{0}^{1},\left( ad(2| k |\pi X_{0}^{1})\right )^{n-1}\left (\frac{2\pi}{\left | \omega_{bif}\right |}
\left [(a^{'}(0)- k (\omega_{0,\mu})^{'}|_{\mu= 0})X_{0}^{1}+b^{'}(0)X_{1}+c^{'}(0)X_{2}\right ]\right )\right ]
\end{multline}
will contain only linear combinations of $X_{1}$ and $X_{2}$.\\
Then with respect to the basis $\{X_{0}^{1}, X_{1}, X_{2}\}$,\\
\[(DF)_{(0,0,0,0)}= \left ( \begin{array}{ccc}
                                   0 & -1 & 0\\
                                   0 & 0 & -1\\
                                   0 & 0 & 0 \\
                             \frac{2\pi}{\left | \omega_{bif}\right |}\left(a^{'}(0)- k (\omega_{0,\mu})^{'}|_{\mu= 0}\right)
                             &\frac{2\pi}{\left | \omega_{bif}\right | }b^{'}(0)+*&\frac{2\pi}{\left |\omega_{bif}\right | }c^{'}(0)+*
                             \end{array}
                    \right ),
\] where $*$ denotes other terms (Below we show that in fact $\frac{2\pi}{\left | \omega_{bif}\right | }b^{'}(0)+*=0$
and $\frac{2\pi}{\left | \omega_{bif}\right | }c^{'}(0)+*=0$).\\
Let $U=\sum_{n= 0}^{\infty} \frac{(-1)^{n}}{(n+1)!}\left ( ad(\frac{2\pi\left | X_{0} \right |}
{\left | \omega_{bif}\right |}X_{0}^{1})\right ) ^{n}\left (\frac{2\pi}{\left | \omega_{bif}\right |}
\left [(a^{'}(0)- k (\omega_{0,\mu})^{'}|_{\mu= 0})X_{0}^{1}+b^{'}(0)X_{1}+c^{'}(0)X_{2}\right ]\right )$.
Using Proposition \ref{prop:geometrical_interp} (1), we get
\begin{equation}\label{E:mike1}
\begin{split}
\overrightarrow U &= \overrightarrow{\sum_{n= 0}^{\infty} \frac{(-1)^{n}}{(n+1)!} \left (
ad(2\pi |k| X_{0}^{1})\right ) ^{n}\left (\frac{2\pi}{\left | \omega_{bif}\right |}
\left [(a^{'}(0)- k (\omega_{0,\mu})^{'}|_{\mu= 0})X_{0}^{1}+b^{'}(0)X_{1}+ c^{'}(0)X_{2}\right ]\right )}\\
&=\frac{2\pi}{\left | \omega_{bif}\right |}\sum_{n= 0}^{\infty} \frac{(-1)^{n}}{(n+1)!} (2\pi |k| X_{0}^{1})^{n}
\left [(a^{'}(0)-k(\omega_{0,\mu})^{'}|_{\mu= 0})\overrightarrow X_{0}^{1}+b^{'}(0)\overrightarrow X_{1}+ c^{'}(0)
\overrightarrow X_{2}\right ].
\end{split}
\end{equation}
If we apply the relation (\ref{E:auxadjoint}) for $Z=2\pi |k| X_{0}^{1}$, then taking into account that
$\left | Z \right |= 2\pi |k|$, we get
\begin{equation}\label{E:mike2}
\sum_{n= 0}^{\infty} \frac{(-1)^{n}}{(n+1)!}(2\pi |k| X_{0}^{1})^{n}= I_{3}+(X_{0}^{1})^{2}.
\end{equation}
If we substitute (\ref{E:mike2}) into (\ref{E:mike1}), it follows that
\begin{equation}
\begin{split}
\overrightarrow{U}&=\frac{2\pi}{\left | \omega_{bif}\right |}[I_{3}+(X_{0}^{1})^{2}] \left [(a^{'}(0)- k
(\omega_{0,\mu})^{'}|_{\mu= 0})\overrightarrow X_{0}^{1}+ b^{'}(0)\overrightarrow X_{1}+c^{'}(0)\overrightarrow X_{2}\right ]\\
&=\frac{2\pi}{\left | \omega_{bif}\right |}\left [ (a^{'}(0)- k(\omega_{0,\mu})^{'}|_{\mu= 0}) \overrightarrow
X_{0}^{1}+b^{'}(0)\overrightarrow X_{1}+ c^{'}(0)\overrightarrow X_{2}\right .  \\  &  \left . +(a^{'}(0)- k(\omega_{0,\mu})^{'}|_{\mu= 0})(X_{0}^{1})^{2}\overrightarrow X_{0}^{1}+ b^{'}(0)(X_{0}^{1})^{2} \overrightarrow X_{1}+ c^{'}(0)(X_{0}^{1})^{2} \overrightarrow X_{2}\right ].
\end{split}
\end{equation}
Since $(X_{0}^{1})^{2}\overrightarrow X_{0}^{1}=\overrightarrow 0$, $(X_{0}^{1})^{2}\overrightarrow X_{1}=-\overrightarrow X_{1}$ and $(X_{0}^{1})^{2}\overrightarrow X_{2}=-\overrightarrow X_{2}$, we get $\overrightarrow{U}=\frac{2\pi}{\left | \omega_{bif}\right |}(a^{'}(0)- |k |(\omega_{0,\mu})^{'}|_{\mu= 0})\overrightarrow X_{0}^{1}$.\\
In order to have the rank of $(DF)_{(0,0,0,0)}$ equal to 3, it is necessary and sufficient that \[a^{'}(0) \neq k (\omega_{0, \mu})^{'}|_{\mu= 0}.\]
We note that depending on the function $X^{G}(t, \lambda, \mu)$ we can establish if $\alpha(\epsilon)= 0$ and/or $\beta(\epsilon)= 0$.\\
\emph{ \textbf{Another proof of Theorem \ref{thm:Resonant_drift1}}} We have that $A(T(\lambda,\mu),\lambda,\mu)=
e^{X(\lambda,\mu)T(\lambda,\mu)}$, where $X(\lambda,\mu)$ is defined as in Proposition \ref{prop:Frequency1}.\\
Let $X(\lambda,\mu)=a_{1}(\lambda,\mu)X_{0}^{1}+b_{1}(\lambda,\mu)X_{1}+ c_{1}(\lambda,\mu)X_{2}$. Then, since $X(\lambda,\mu)$ is sufficiently $CS$, we get that $a_{1}$, $b_{1}$, $c_{1}$ are sufficiently $CS$ and $b_{1}(0,\mu)=O(\mu),c_{1}(0,\mu)=O(\mu)$.\\
Since $e^{X(0,0)T(0,0)}= A(T(0,0),0,0)= I_{3}$, then $X(0,0)=O_{3}$, that is $a_{1}(0,0)=0$.\\
We want to find a sufficiently $CS$ branch $\mu=\mu(\lambda)$ such that $\mu(0)=0$ and $a_{1}(\lambda,\mu(\lambda))=0$, for $\lambda \geq 0$ small. If $(a_{1})_{\mu}(0,0) \neq 0$, then by applying the implicit function theorem we get the existence of the required sufficiently $CS$ branch.\\
Now we show that $(a_{1})_{\mu}(0,0) \neq 0$ is equivalent with $a^{'}(0) \neq k (\omega_{0, \mu})^{'}|_{\mu= 0}$.\\
We have that $A(T(0,\mu),0,\mu)= e^{X(0,\mu)T(0,\mu)}= e^{X_{G}(0,0,\mu)T(0,\mu)}$ for $\left | \mu \right |$ small, then we get $X(0,\mu)T(\mu,0)= X_{G}(0,0,\mu)T(\mu,0)+\frac{2l(\mu)\pi}{\left | X_{G}(0,0,\mu)
\right | } X_{G}(0,0,\mu)$ for some $l(\mu) \in \mathbf{Z}$ for $\left | \mu \right |$ small.\\
Thus, $\left [ a_{1}(0,\mu)X_{0}^{1}+b_{1}(0,\mu)X_{1}+c_{1}(0,\mu)X_{2}\right ]T(0,\mu)= \\
\left [a(\mu)X_{0}^{1}+b(\mu)X_{1}+c(\mu)X_{2}\right ]T(0,\mu)+ \frac{2l(\mu)\pi}{\sqrt{a(\mu)^{2}+O(\mu^{2})}}\left
[a(\mu)X_{0}^{1}+ b(\mu)X_{1}+c(\mu)X_{2}\right ]$ for $\left | \mu \right | \geq 0$ small.\\
This implies that
\begin{equation}\label{E:eqauxres}
a_{1}(0,\mu)=a(\mu)+\frac{a(\mu)}{\sqrt{a(\mu)^{2}+O(\mu^{2})}}l(\mu)\omega_{0,\mu}\cdot sgn(\omega_{bif}),
\end{equation}
which for $\mu=0$ gives $l(0)=-|k|$, where we have used the fact that $\left | X_{0} \right | =k \omega_{bif}>0$ implies $sgn(k)= sgn(\omega_{bif})$, and thus, $k \cdot sgn(\omega_{bif}) =|k|$. Since $l(\mu) \in \mathbb{Z}$ is continuous, we get $l(\mu)=-| k |$ for $\left | \mu \right | \geq 0$ small. \\
Then by taking into account that $sgn(\omega_{bif})=sgn(k)$, (\ref{E:eqauxres}) becomes
\begin{equation}\label{E:eqauxres1}
a_{1}(0,\mu)=a(\mu)- k \frac{a(\mu)}{\sqrt{a(\mu)^{2}+ O(\mu^{2})}}\omega_{0,\mu}.
\end{equation}
By differentiation of (\ref{E:eqauxres1}), it follows that $(a_{1})_{\mu}(0,0)= a^{'}(0) - k (\omega_{0, \mu})^{'}|_{\mu= 0}\neq 0$.\\\\
\end{proof}

\begin{proof}[\textbf{Proof of Corollary \ref{cor: Resonant_drift2}}]
The conclusions of Corollary \ref{cor: Resonant_drift2} result in the same way as in the proof of Theorems
\ref{prop:Frequency1} and \ref{thm:Periodic_part}, except the scaling of the primary
frequencies $ \left | X(\lambda, \mu(\lambda)) \right |$ that results from Remark 3.5 in \cite{Wu2}.
\end{proof}

\appendix
\section{BCH formula in $so(3)$}\label{S:app1_BCH}
It is known that $SO(3)$ is diffeomorphic as a manifold to the real projective space $\mathbb{R} P^{3}$ (see \cite{Ro}).
If $Y=\left(\begin{array}{ccc}
                  0 & a & -b \\
                  -a & 0 & c \\
                  b & -c & 0\\
             \end{array}
       \right) \in so(3)$,
then we define $\overrightarrow {Y}= \left( \begin{array}{c} c\\ b\\ a\\
\end{array} \right)$. Also, we define $\left | Y \right | =\left \| \overrightarrow{Y} \right \| $.
A model for the space $\mathbb{R} P^{3}$ is the set
$D= \{ \overrightarrow y \in \mathbb{R}^{3} \mid \left \| \overrightarrow y \right \| \leq \pi ,
\mbox{ whith the antipodals points of the norm } \left | y \right | = \pi \mbox{ identified }\}$. In fact, D is the
quotient set $E/\sim$, where $ \sim$ is the equivalence relation $\vec y \sim \vec z \mbox { iff
} z=-y, \left | y \right | =\pi $ and $E= \{ \overrightarrow y \in \mathbb{R}^{3} \mid \left \| \overrightarrow y
\right \| \leq \pi \}$. The set $D$ is considered with the quotient topology. Sometimes,
instead of $\overrightarrow y \in \mathbb{R}^{3}$, $\left | \overrightarrow y \right | \leq \pi$,
we use $y \in so(3)$, in which case we denote the equivalence class $[y]=[\overrightarrow y]$.
The projection map $p \colon E \rightarrow D$ , $p(y)=[y]$ is smooth.
There exists a unique smooth function $d^{*} \colon SO(3) \rightarrow D$ such that $e^{d^{*}(A)}=A$.

\begin{definition}\label{def:BCH1}
For any $X$, $Y \in so(3)$, we define $BCH(X,Y)= d^{*}(e^{X}e^{Y})$, where $d^{*}$ is defined above.
Clearly, we have $e^{BCH(X,Y)}=e^{X}e^{Y}$ and $BCH(X,Y) \in D$, for any $X$, $Y \in so(3)$.
\end{definition}
\begin{theorem}[$BCH$ Formula in $so(3)$, \cite{Co, En, Se}]\label{thm:BCH_formula}
The $BCH$ formula in $so(3)$ has the form
\begin{equation}\label{E:BCH_formula}
BCH(X,Y)= [\alpha X+ \beta Y+ \gamma [X,Y]] \mbox{ for } X, \, Y \in so(3),
\end{equation}
where
\[\alpha=k(X,Y)h_{\alpha}(X,Y), \, \beta=k(X,Y)h_{\beta}(X,Y),\, \gamma=k(X,Y)h_{\gamma}(X,Y),\]
and
\[e= \cos \frac{\left | X \right |}{2}\cos\frac{\left | Y \right |}{2}-  \sin\frac{\left | X \right |}{2}\sin\frac{\left | Y \right |}{2}\cos(\angle(\vec X,\vec Y)),\]
\[a_{1}= \sin \frac{\left | X \right |}{2}\cos \frac{\left | Y \right |}{2},\, a=a_{1}e ,\]
\[b_{1}=  \sin \frac{\left | Y \right |}{2} \cos \frac{\left | X \right |}{2},\, b=b_{1}e ,\]
\[c_{1}= \sin \frac{\left | X \right |}{2}\sin \frac{\left | Y \right |}{2},\, c=c_{1}e ,\]
\[d_{1}= \sqrt{a_{1}^{2}+b_{1}^{2}+2a_{1}b_{1}\cos(\angle(\vec X,\vec Y))+
c_{1}^{2}(\sin(\angle(\vec X,\vec Y))^{2}},\]
\[d= d_{1}\left | e \right | ,\]
where $\angle(\vec X, \vec Y)$ is the angle between the two vectors $\overrightarrow X$ and $\overrightarrow Y$,
\[h_{\alpha}(X,Y) = \begin{cases}
                     \frac{a_{1}}{\left | X \right | } &\text{if $X \neq O_{3}$;}\\
                     \cos \frac{\left | Y \right | }{2} &\text{if $X= O_{3}$,}\\
                    \end{cases}\]
\[h_{\beta}(X,Y) = \begin{cases}
                      \frac{b_{1}}{\left | Y \right | } &\text{if $Y \neq O_{3}$;}\\
                      \cos \frac{\left | X \right | }{2} &\text{if $Y= O_{3}$,}\\
                    \end{cases}\]
\[h_{\gamma}(X,Y) = \begin{cases}
                        \frac{c_{1}}{\left | X \right | \left | Y \right |} &\text{if $X \neq O_{3},Y \neq O_{3}$;}\\
                        \frac{\sin \frac{\left | Y \right | }{2}}{\left | Y \right |} &\text{if $X= O_{3},Y \neq O_{3}$;}\\
                        \frac{\sin \frac{\left | X \right | }{2}}{\left | X \right |} &\text{if $Y= O_{3},X \neq O_{3}$;}\\
                        1 &\text{if $Y= O_{3},X = O_{3}$,}
                     \end{cases}\]
\[k(X,Y)=\begin{cases}
           s\frac{\arcsin(d)}{d_{1}} &\text{if $(e^{X}e^{Y})^{2} \neq I_{3}$,\, $e^{X}e^{Y}$  has eigenvalues with positive real parts;}\\
           s\frac{\pi-\arcsin(d)}{d_{1}} &\text{if $(e^{X}e^{Y})^{2} \neq I_{3}$,\, $e^{X}e^{Y}$  has two eigenvalues with negative } \\ &\text{ or zero real parts;}\\
       \pi &\text{if $(e^{X}e^{Y})^{2} = I_{3}$,\, $e^{X}e^{Y} \neq I_{3}$;}\\
           s &\text{if $ e^{X}e^{Y}= I_{3}$,}
        \end{cases}
\]
where $s= sgn(e)=\begin{cases}
                  1, &\text{if $e >0$;}\\
                 -1, &\text{if $e < 0$,}
                 \end{cases}$ for any $(X,Y) \in so(3) \times so(3)$, such that
$(e^{X}e^{Y})^{2} \neq I_{3}$ or $e^{X}e^{Y}= I_{3}$. The functions $\alpha^{2}$, $\beta^{2}$, $\gamma ^{2}$ are
smooth on $so(3) \times so(3)$ and $\left | \alpha \right | $, $\left | \beta \right | $, $\left | \gamma \right | $
are continuous on $so(3) \times so(3)$. Also, the function $BCH$ is smooth from $so(3) \times so(3)$ into $D$.
\end{theorem}
\begin{proposition}\label{prop:geometrical_interp}
\begin{enumerate}
\item \cite{Ro} For any $X $, $Y \in so(3)$,
\begin{equation}\label{E:adj1}
\begin{array}{lll}
\overrightarrow {ad(X)Y} &=& X\overrightarrow Y \mbox{ and }\\
\overrightarrow {(ad(X))^{n}Y} &=& X^{n}\overrightarrow Y \mbox{ for any integer } n >0.\\
\end{array}
\end{equation}
\item For any $A \in SO(3)$ and any $X \in so(3)$, we have
\begin{enumerate}
\item
\begin{equation}\label{E:adj2}
\left |  AXA^{-1} \right | = \left | X \right | ;
\end{equation}
\item \cite{Ro} \[AXA^{-1}= B \mbox{ if and only if } A\overrightarrow X= \overrightarrow B.\]
\end{enumerate}
\item \cite{MR} For any $X \in so(3)$, we have:
\begin{equation}\label{E:adj3}
\begin{array}{ccc}
X^{2n} &=& (-1)^{n-1}\left | X \right |^{2(n-1)}X^{2} \mbox{ for any } n \geq 1,\\
X^{2n+1} &=& (-1)^{n}\left | X \right |^{2n}X \mbox{ for any } n \geq 0.\\
\end{array}
\end{equation}
\end{enumerate}
\end{proposition}
The properties of the exponential map $e \colon so(3) \rightarrow SO(3)$ are:
\begin{proposition}\label{prop:properties_Exp}
\begin{enumerate}
\item \cite{BD, On} The exponential map $exp \colon so(3) \rightarrow SO(3)$ is surjective.
\item \cite{On, Ro} The exponential map $exp$ is a smooth function on $so(3)$ and its differential is
given by:
\[(d(exp))_{X}(Y)= e^{X}\sum_{n= 0}^{\infty}\frac{(-1)^{n}}{(n+1)!}(ad (X))^{n}(Y) \mbox { for any } X,\, Y \in so(3)\] or
\[(d(exp))_{X}(Y)= \sum_{n= 0}^{\infty}\frac{1}{(n+1)!}(ad (X))^{n}(Y)e^{X} \mbox { for any } X,\, Y \in so(3)\]
( we will use the first formula later ). Moreover, it is a local diffeomorphism
near any $X \in so(3)$ if and only if the operator $ad(X)$ has no eigenvalues of the form $2\pi ik$ with $ k \neq 0$, that is if and
only if $\left | X \right | \neq 2k\pi$  for $k \in \mathbb{Z}$, $k \neq 0$. \item \cite{Ro} $\frac{d }{d t} \left ( e^{X(t)}
\right )= \dot X(t)e^{X(t)}$ if and only if $X(t)= Xg(t)$, where $X \in so(3)$ and $g \colon \mathbb{R} \rightarrow
\mathbb{R}$ is a $C^{1}$ function.
\item \cite{En, Se} The Rodrigues' formula holds:
\[e^{X}= I_{3}+\frac{\sin\left | X \right | }{\left | X \right | }X+ 2\frac{\sin^{2}\frac{\left | X \right| }{2}} {\left | X \right
|^{2}}X^{2} \mbox{ for any }  X \in so(3),\] where we take the limit when $X= O_{3}$.
\item \cite{Ro} The exponential map $exp$ maps  any $ X \in so(3)$, $X \neq O_{3}$ to the right-handed rotation $A \in SO(3)$ with angle $ \left | X \right | $ around
$\overrightarrow X$.
\item \cite{Ro} $e^{X}= I_{3}$ if and only if $\left | X \right | = 2k\pi$ for some $k \in \mathbb{Z}$.
\item \cite{Ro} $e^{X}= e^{Y}$ if and only if either $\left | X \right | = 2k\pi= \left | Y \right | $ or $\overrightarrow Y=
\overrightarrow X +\frac{2k\pi}{\left | X \right | } \overrightarrow X$ for some $k \in \mathbb{Z}$.\\
\end{enumerate}
\end{proposition}

\section{Equivariant Center Manifold Reduction for Reaction-Diffusion Systems on $r\mathbf{S^{2}}$}\label{S:app2_CMR}

We consider the reaction-diffusion system (\ref{E:RDparam}).
\begin{theorem}[\cite{FSSW, SSW}]\label{thm:CMR_REparam}
Let $u_{0} \in \mathbf{Y^{\alpha}}$ be a relative equilibrium that is not an equilibrium for (\ref{E:RDparam}) at
$\lambda= 0$ and such that the stabilizer of $u_{0}$ is $\Sigma_{u_{0}}= {I_{3}}$. Let $L$ be  the linearization of
the right-hand side of (\ref{E:RDparam}) with respect to the rotating wave
$\mathbf{\Phi}(t, u_{0}, 0)= e^{X_{0}t}u_{0}$  at $\lambda=0$ in the co-rotating frame, that is
\[L= D\Delta_{S}+D_{u}F(u_{0},0)-X_{0}.\]
Suppose that:
\begin{enumerate}
\item $\sigma(L) \cap \{z \in \mathbb{C} \mid Re \, (z) \geq 0 \}$ is a spectral set with spectral projection $P_{*}$, and dim$(R(P_{*})) < \infty ;$ \item the semigroup $e^{Lt}$ satisfies $\left | e^{Lt}|_{R(1-P_{*})} \right | \leq Ce^{-\beta _{0}t}$ for some $\beta _{0} > 0$ and $ C > 0$.
\end{enumerate}
Then there exists  a sufficiently smooth parameter-dependent center manifold $M^{cu}_{u_{0}}(\lambda)$ for the relative equilibrium $u_{0}$.
Let $V_{*}$ be the orthogonal complement of $T_{u_{0}}(SO(3)u_{0})$ in $E^{cu}= R(P_{*})$.\\
Then the center manifold $M^{cu}_{u_{0}}(\lambda)$ is diffeomorphic to $SO(3) \times V_{*}$ for $\left | \lambda \right | $ small. Furthermore, there exist sufficiently smooth functions $X_{G} \colon V_{*} \times \mathbb{R} \rightarrow so(3)$ and $X_{N} \colon V_{*} \times \mathbb{R} \rightarrow V_{*}$ such that any solution of
\begin{equation}\label{E:reduceddiffeqparam}
    \begin{array}{lll}
    \dot{A} &=& AX_{G}(q, \lambda),\\
    \dot{q} &=& X_{N}(q, \lambda),\\
    \end{array}
\end{equation}
on $SO(3) \times V_{*}$ corresponds to a solution of the reaction-diffusion system (\ref{E:RDparam}) on
$M_{u_{0}}^{cu}(\lambda)$ under the diffeomorphic identification  $(A,q) \rightarrow A\Psi(q)$, with $\Psi$ a local chart
from $V$ to $M$, for $\left | \lambda  \right | $ small.
Also, $X_{G}(0, 0)= X_{0}$, $X_{N}(0, 0)= 0$ and $\sigma(D_{u}X_{N}(0,0))= \sigma(Q_{*}L|_{V_{*}})$, where $Q_{*}$ is the projection onto $V_{*}$ along $T_{u_{0}}(SO(3)u_{0})$.
\end{theorem}
\begin{remark}
Theorems \ref{thm:CMR_REparam} is valid for more than one parameter.
\end{remark}

\section{Some Computations}

We present the proof of (\ref{E:Numerical1}): \\
If we parameterize $SO(3)$ by Euler angles, that is $A= R_{z}(\psi)R_{x}(\theta)R_{z}(\phi)$, where $\phi \in [0,2\pi)$, $\theta \in [0,\pi]$, $\psi \in [0,2\pi)$, and we substitute this into the equation $\dot{A}=AX_{G}(q,\lambda)$, we get
\begin{multline}
\dot{\psi}\left(\frac{d R_{z}}{d \psi}(\psi)\right)R_{x}(\theta)R_{z}(\phi)+ \dot{\theta} R_{z}(\psi)\left(\frac{d
R_{x}}{d\theta}(\theta)\right)R_{z}(\phi)+\dot{\phi}R_{z}(\psi)R_{x}(\theta)\left(\frac{d R_{z}} {d \phi}(\phi)\right )\\ = R_{z}(\psi)R_{x}(\theta)R_{z}(\phi)\left [ F^{x}(q,\lambda)L_{x}+F^{y}(q,\lambda)L_{y}+F^{z}(q,\lambda)L_{z}
\right ]
\end{multline}
or
\begin{equation}\label{E:auxa}
\begin{split}
& \dot{\psi}R_{z}(-\phi)R_{x}(-\theta)R_{z}(-\psi)\left (\frac{d R_{z}}{d \psi}(\psi)\right )R_{z}(\phi)R_{x}(\theta)R_{z}(\phi) \\ & +\dot{\theta}R_{z}(-\phi)R_{x}(-\theta)\left ( \frac{d R_{x}}{d \theta}(\theta)\right )R_{z}(\phi)+\dot{\phi}R_{z}(-\phi)\left ( \frac{d R_{z}}{d \phi}(\phi)\right ) \\
& =F^{x}(q,\lambda)L_{x}+F^{y}(q,\lambda)L_{y}+F^{z}(q,\lambda)L_{z}.
\end{split}
\end{equation}
By computation, it follows that
\begin{equation}\label{E:auxb}
\begin{array}{rrl}
R_{z}(-\phi)\left ( \frac{d R_{z}}{d \phi}(\phi)\right ) &=& L_{z},\\
R_{z}(-\phi)R_{x}(-\theta)\left (\frac{d R_{x}}{d \theta}(\theta)\right )R_{z}(\phi)&=& (-\sin \phi)L_{y}+ \cos \phi L_{x},\\
R_{z}(-\phi)R_{x}(-\theta)R_{z}(-\psi)\left (\frac{d R_{z}}{d \psi}(\psi)\right )R_{z}(\phi)R_{x}(\theta)R_{z}(\phi) &=&
\cos \theta L_{z}+\cos \phi \sin \theta L_{y}+ \sin \phi \sin \theta L_{x}.\\
\end{array}
\end{equation}
If we substitute (\ref{E:auxb}) into (\ref{E:auxa}), we get
\begin{multline}\label{E:auxc}
\dot{\psi}(\cos \theta L_{z}+\cos \phi \sin \theta L_{y}+ \sin \phi \sin \theta L_{x}) +\dot{\theta}((-\sin \phi)L_{y}+ \cos \phi L_{x})+\dot{\phi}L_{z}\\= F^{x}(q,\lambda)L_{x}+F^{y}(q,\lambda)L_{y}+F^{z}(q,\lambda)L_{z}.
\end{multline}
If we identify the corresponding coefficients of $L_{x}$, $L_{y}$ and $L_{z}$ from both sides of the equation
(\ref{E:auxc}), we get the system
\begin{equation}\label{E:auxd}
\begin{array}{rrl}
\dot{\psi}\cos \theta+\dot{\phi} &=& F^{z}(q,\lambda),\\
\dot{\psi}\cos \phi \sin \theta -\dot{\theta}\sin \phi &=& F^{y}(q,\lambda),\\
\dot{\psi} \sin \phi \sin \theta+\dot{\theta}\cos \phi &=& F^{x}(q,\lambda).\\
\end{array}
\end{equation}
Then, by multiplying the second equation in (\ref{E:auxd}) by $-\sin \phi$ and the third equation in (\ref{E:auxd}) by $\cos \phi$ and adding them, we get $\dot{\theta} = -F^{y}(q,\lambda)\sin \phi+ F^{x}(q, \lambda)\cos \phi$.
Similarly, by multiplying the second equation in (\ref{E:auxd}) by $\cos \phi$ and the third equation in (\ref{E:auxd}) by $\sin \phi$ and adding them, we get $\dot{\psi} = \frac{1}{\sin \theta}\left [F^{y}(q,\lambda)\cos \phi +F^{x}(q,\lambda)\sin \phi \right ]$. If we substitute $\dot{\psi}$ given by the above relation into the first equation of the system (\ref{E:auxd}), we get $\dot{\phi} = F^{z}(q,\lambda)-\cot \theta \left [F^{y}(q,\lambda) \cos \phi +F^{x}(q,\lambda) \sin \phi \right ]$.

\section*{Acknowledgments}
I wish to thank you to my Ph.D. supervisor, Victor G. LeBlanc, for his help to write this article. I am also grateful to the
Natural Sciences and Engineering Research Council of Canada for a PGS-B Scholarship.

\end{document}